\documentclass[twoside]{article}%
\usepackage{amssymb}
\usepackage{amsmath}
\usepackage{fancyhdr}
\usepackage{amsfonts}
\usepackage{graphicx}%
\newtheorem{theorem}{Theorem}[section]

\newtheorem{corollary}[theorem]{Corollary}

\newtheorem{definition}[theorem]{Definition}

\newtheorem{lemma}[theorem]{Lemma}

\newtheorem{proposition}[theorem]{Proposition}
\newtheorem{remark}[theorem]{Remark}

\newenvironment{proof}[1][Proof]{\textbf{#1.} }{\ \rule{0.5em}{0.5em}}
\begin{document}

\title{A Mean Value Theorem for Closed Geodesics on Congruence Surfaces}
\author{Vladimir Lukianov\\Tel-Aviv University}
\maketitle

\paragraph{\thispagestyle{empty}Abstract.}

We define a weighted multiplicity function for closed geodesics of given
length on a finite area Riemann surface. These weighted multiplicities appear
naturally in the Selberg trace formula, and in particular their mean square
plays an important role in the study of statistics of the eigenvalues of the
Laplacian on the surface.

In the case of the modular domain, E. Bogomolny, F. Leyvraz and C. Schmit gave
a formula for the mean square, which was rigorously proved by M. Peter. In
this paper we calculate the mean square of weighted multiplicities for some
surfaces associated to congruence subgroups of the unit group of a rational
quaternion algebra, in particular for congruence subgroups of the modular
group. Remarkably, the result turns out to be a rational multiple of the mean
square for the modular domain.

\section{Introduction}

The aim of the present work is to give a tool for finding the eigenvalue
statistics of the hyperbolic Laplacian%
\[
\Delta=y^{2}\left(  \dfrac{\partial^{2}}{\partial x^{2}}+\dfrac{\partial^{2}%
}{\partial y^{2}}\right)
\]
defined on a Riemann surface $\Omega=\Gamma\backslash\mathbb{H}$, for a
discrete cofinite subgroup $\Gamma$ of $PSL_{2}(\mathbb{R})$ of "congruence
type" (see below).

\bigskip

Let $\lambda_{0}=0<\lambda_{1}\leq\lambda_{2}\leq...$ be the discrete spectrum
of the hyperbolic Laplacian on $\Omega$. It is known \cite{hejhal} that it
satisfies the Weyl's law:%
\[
N(T):=\#\{0\leq r_{j}\leq T\}=\dfrac{vol(\Omega)}{4\pi}T^{2}+O(T\ln T)\text{,}%
\]
where $\lambda_{j}=\dfrac{1}{4}+r_{j}^{2}$. We define a smooth version of the
counting function as follows. Let $f$ be an even test function with the
compactly supported\ smooth Fourier transform $\widehat{f}$. Define%
\[
N_{f}(\tau)=%
%TCIMACRO{\dsum \limits_{j\geq0}}%
%BeginExpansion
{\displaystyle\sum\limits_{j\geq0}}
%EndExpansion
f(L(r_{j}-\tau))+f(L(-r_{j}-\tau))\text{.}%
\]
For example, if $f$ is the characteristic function of the interval
$(-1/2,1/2)$, then $N_{f}(\tau)$ counts the number of $r_{j}$ lying in the
intervals $\pm(\tau-1/2L,\tau+1/2L)$.

\bigskip

\bigskip

Put
\[
h(r)=f(L(r-\tau))+f(L(-r-\tau))
\]
then we can use the Selberg Trace Formula to express $N_{f}(\tau)$ as:%
\begin{align}
N_{f}(\tau)  &  =%
%TCIMACRO{\dsum \limits_{j\geq0}}%
%BeginExpansion
{\displaystyle\sum\limits_{j\geq0}}
%EndExpansion
h(r_{j})=\{identity~contribution\}+\nonumber\\
&  +\{hyperbolic~contribution\}+\label{stf}\\
&  +\{elliptic~contribution\}+\nonumber\\
&  +\{parabolic~and~continuous~spectrum~contribution\}\nonumber
\end{align}
We say that the element $T$ of $\Gamma$ is hyperbolic if $\left\vert
trT\right\vert >2$. Such an element is conjugate in $PSL_{2}(\mathbb{R})$ to a
diagonal matrix $\left[
\begin{array}
[c]{cc}%
\lambda & 0\\
0 & \lambda^{-1}%
\end{array}
\right]  $, where $\lambda>1$ is real. For the element $T$, such that
$\left\vert trT\right\vert =t$ we define the norm\footnote{The number
$\lambda^{2}$ is called also a multiplier of $T$.} of $T$ to be
\[
\mathcal{N}(T)=\lambda^{2}=\left(  \dfrac{t+\sqrt{t^{2}-4}}{2}\right)
^{2}\text{.}%
\]
For us the most important term in formula (\ref{stf}) is the hyperbolic
contribution term, which is defined explicitly as:%
\[%
%TCIMACRO{\dsum \limits_{t>2}}%
%BeginExpansion
{\displaystyle\sum\limits_{t>2}}
%EndExpansion%
%TCIMACRO{\dsum \limits_{\substack{\{T\}~hyperbolic\\\left\vert trT\right\vert
%=t}}}%
%BeginExpansion
{\displaystyle\sum\limits_{\substack{\{T\}~hyperbolic\\\left\vert
trT\right\vert =t}}}
%EndExpansion
\dfrac{\ln\mathcal{N}(T_{0})}{\mathcal{N}(T)^{1/2}-\mathcal{N}(T)^{-1/2}}%
g(\ln\mathcal{N}(T))\text{.}%
\]
There are in general several conjugacy classes with $\left\vert trT\right\vert
=t$. Define the weighted multiplicity function $\beta_{\Gamma}(t)$ by:%
\[
\beta_{\Gamma}(t)=\dfrac{1}{4}%
%TCIMACRO{\dsum \limits_{\substack{\{T\}~\subset\Gamma\\T=T_{0}^{k}%
%~is~hyperbolic\\\left\vert trT\right\vert =t}}}%
%BeginExpansion
{\displaystyle\sum\limits_{\substack{\{T\}~\subset\Gamma\\T=T_{0}%
^{k}~is~hyperbolic\\\left\vert trT\right\vert =t}}}
%EndExpansion
\dfrac{\ln\mathcal{N}(T_{0})}{\mathcal{N}(T)^{1/2}-\mathcal{N}(T)^{-1/2}%
}\text{,}%
\]
where $T_{0}$ is a primitive hyperbolic element, that is $T_{0}$ is not a
power of any other hyperbolic element. In this notation we can rewrite the
hyperbolic contribution term as follows%
\[%
%TCIMACRO{\dsum \limits_{t>2}}%
%BeginExpansion
{\displaystyle\sum\limits_{t>2}}
%EndExpansion
\beta_{\Gamma}(t)g(\ln\mathcal{N}(T))\text{,}%
\]
so the information about the weighted multiplicities will be useful for
understanding the behavior of this term.

From the prime geodesic theorem it follows that the mean value of
$\beta_{\Gamma}(t)$ is unity:%
\[
\lim\limits_{N\rightarrow\infty}\dfrac{1}{N}%
%TCIMACRO{\dsum \limits_{2<t\leq N}}%
%BeginExpansion
{\displaystyle\sum\limits_{2<t\leq N}}
%EndExpansion
\beta_{\Gamma}(t)=1\text{.}%
\]
Our main result is a computation of the mean square of the weighted
multiplicities for the case that $\Gamma$ is congruence subgroup of the unit
group of a quaternion algebra over the rationals.

Let $B$ be a quaternion algebra over $\mathbb{Q}$. Let%
\[
N_{B}:B\rightarrow\mathbb{Q}%
\]%
\[
Tr_{B}:B\rightarrow\mathbb{Q}%
\]
denote respectively the reduced norm and the reduced trace of the elements of
$B$. If $\alpha\in B$, then
\[
N_{B}(\alpha)=\alpha+\overline{\alpha}%
\]%
\[
Tr_{B}(\alpha)=\alpha\overline{\alpha}%
\]
where $\overline{\alpha}$ is the conjugate of $\alpha$.

We recall that quaternion algebra over field $F$ is either division algebra,
or is isomorphic to $M_{2}(F)$.

Let $d_{B}$ be the reduced discriminant of $B$, that is $d_{B}$ is the product
of all primes $p$, such that $B_{p}=B\otimes_{\mathbb{Q}}\mathbb{Q}_{p}$ is a
division algebra. Note that $d_{B}>1$ if and only if $B$ is a division
algebra. We assume that $B$ is indefinite, that is $B_{\infty}:=$%
\ $B\otimes_{\mathbb{Q}}\mathbb{R\cong}M_{2}(\mathbb{R})$.

let $R$ be an order in $B$, then for $\alpha\in R$ we have that both the
reduced trace and the reduced norm of $\alpha$ are integers.

Fix now the isomorphism $B\otimes_{\mathbb{Q}}\mathbb{R\cong M}_{2}%
(\mathbb{R})$, and set%
\[
\Gamma_{R}=\{\alpha\in R\mid N_{B}(\alpha)=1\}\text{.}%
\]
Under above isomorphism $\Gamma_{R}$ is identified with a subgroup of
$SL_{2}(\mathbb{R)}$, which is a cofinite Fuchsian group, (and is co-compact
if $d_{B}>1$). Moreover, the traces of its elements are \textit{integers}. In
particular, one can consider $B=M_{2}(\mathbb{Q)}$, and for the given natural
$Q$ the set
\[
R=\left\{  \left[
\begin{array}
[c]{cc}%
a & b\\
Qc & d
\end{array}
\right]  \in M_{2}(\mathbb{Z)}\right\}  \text{.}%
\]
Then $B$ is indefinite quaternion algebra over $\mathbb{Q}$ and $R$ is an
order in $B$. Here
\begin{equation}
\Gamma_{R}=\Gamma_{0}(Q)=\left\{  \left[
\begin{array}
[c]{cc}%
a & b\\
c & d
\end{array}
\right]  \in SL_{2}(\mathbb{Z})\mid c\equiv0(\operatorname{mod}Q)\right\}
\text{.} \label{gammaRQ}%
\end{equation}

\begin{remark}
We will often identify $\Gamma_{R}$ with its image in the projective group
$PSL_{2}(\mathbb{R)}$.
\end{remark}

We are now able to state our main theorems:

\begin{theorem}
\label{mean-sq-theorem(congr)}Let
\[
\beta_{Q}(t)=\dfrac{1}{4}%
%TCIMACRO{\dsum \limits_{\substack{\{T\}~\subset\Gamma_{0}(Q)\\T=T_{0}%
%^{k}~is~hyperbolic\\\left\vert trT\right\vert =t}}}%
%BeginExpansion
{\displaystyle\sum\limits_{\substack{\{T\}~\subset\Gamma_{0}(Q)\\T=T_{0}%
^{k}~is~hyperbolic\\\left\vert trT\right\vert =t}}}
%EndExpansion
\dfrac{\ln\mathcal{N}(T_{0})}{\mathcal{N}(T)^{1/2}-\mathcal{N}(T)^{-1/2}%
}\text{,}%
\]
where $\Gamma_{0}(Q)$ is the congruence subgroup of $SL_{2}(\mathbb{Z)}$,
defined in (\ref{gammaRQ}). Then for the squarefree odd $Q$%
\[
\lim\limits_{N\rightarrow\infty}\dfrac{1}{N}%
%TCIMACRO{\dsum \limits_{2<t\leq N}}%
%BeginExpansion
{\displaystyle\sum\limits_{2<t\leq N}}
%EndExpansion
\beta_{Q}^{2}(t)=C_{1}%
%TCIMACRO{\dprod \limits_{q\mid Q}}%
%BeginExpansion
{\displaystyle\prod\limits_{q\mid Q}}
%EndExpansion
\dfrac{2(q^{2}-q-1)(q+1)^{2}}{q(q^{3}+q^{2}-q-3)}=\kappa_{Q}\text{.}%
\]

\end{theorem}%

\[
\]
Here
\[
C_{1}=\lim\limits_{N\rightarrow\infty}\dfrac{1}{N}%
%TCIMACRO{\dsum \limits_{2<t\leq N}}%
%BeginExpansion
{\displaystyle\sum\limits_{2<t\leq N}}
%EndExpansion
\beta_{1}^{2}(t)=\dfrac{1015}{864}%
%TCIMACRO{\dprod \limits_{p\neq2}}%
%BeginExpansion
{\displaystyle\prod\limits_{p\neq2}}
%EndExpansion
\dfrac{p^{2}(p^{3}+p^{2}-p-3)}{(p^{2}-1)^{2}(p+1)}=1.328...
\]
which is the mean square of the weighted multiplicities for the modular domain
($Q=1$), proved by M. Peter \cite{man}, following a conjecture of Bogomolny et
al \cite{bogom}.

\begin{theorem}
\label{mean-sq-theorem(quater)}Let $B$ be an indefinite division quaternion
algebra over $\mathbb{Q}$ with the maximal order $R$, and discriminant $d_{B}%
$. Then for
\[
\beta_{R}(t)=\dfrac{1}{4}%
%TCIMACRO{\dsum \limits_{\substack{\{T\}~\subset\Gamma_{R}\\T=T_{0}%
%^{k}~is~hyperbolic\\\left\vert trT\right\vert =t}}}%
%BeginExpansion
{\displaystyle\sum\limits_{\substack{\{T\}~\subset\Gamma_{R}\\T=T_{0}%
^{k}~is~hyperbolic\\\left\vert trT\right\vert =t}}}
%EndExpansion
\dfrac{\ln\mathcal{N}(T_{0})}{\mathcal{N}(T)^{1/2}-\mathcal{N}(T)^{-1/2}}%
\]
the mean square of weighted multiplicities is
\[
\lim\limits_{N\rightarrow\infty}\dfrac{1}{N}%
%TCIMACRO{\dsum \limits_{2\leq n\leq N}}%
%BeginExpansion
{\displaystyle\sum\limits_{2\leq n\leq N}}
%EndExpansion
\beta_{R}^{2}(n)=C_{1}\cdot%
%TCIMACRO{\dprod \limits_{p\mid d_{B}}}%
%BeginExpansion
{\displaystyle\prod\limits_{p\mid d_{B}}}
%EndExpansion
\dfrac{2(p^{3}-1)(p-1)}{p(p^{3}+p^{2}-p-3)}\text{.}%
\]

\end{theorem}%

\[
\]
We now explain the method of proof. As a first step, we express the weighted
multiplicities in terms of Dirichlet's $L$-functions. For $\Gamma_{0}(Q)$ the
expression is
\begin{equation}
\beta_{Q}(t)=\sum\limits_{\substack{D,v\geq1\\D\text{ is a discriminant}%
\\Dv^{2}=t^{2}-4}}\frac{1}{v}L(1,\chi_{D})\cdot\prod\limits_{q\mid Q}\left\{
\begin{array}
[c]{cc}%
2, & q^{2}\mid D\\
1+\left(  \dfrac{D}{q}\right)  & q^{2}\nmid D
\end{array}
\right\}  \text{.} \label{L-exp-for-beta}%
\end{equation}
This is proved in section \ref{L-func-expr} by connecting the weighted
multiplicities with class numbers and using Dirichlet's class number formula.

The principal tool in our approach, following Peter, is that the formula
(\ref{L-exp-for-beta}) displays the weighted multiplicity $\beta_{Q}(t)$ as a
"limit periodic function" in a suitable sense (see section
\ref{LimitPeriodBckgr} for background on these). To show this, we approximate
the $L$-functions by a finite Euler product using a zero-density theorem, in a
certain semi-norm coming from the theory of limit periodic functions (section
\ref{Limit-Period-sect}). For computing mean squares, this suffices and allows
us to use Parseval's equality in this setting to express the mean square as%
\[
\lim\limits_{N\rightarrow\infty}\frac{1}{N}\sum\limits_{2\leq t\leq N}%
\beta_{Q}^{2}(t)=\sum\limits_{b\geq1}\sum\limits_{\substack{1\leq a\leq
b\\\gcd(a,b)=1}}\left\vert \widehat{\beta_{Q}}\left(  \frac{a}{b}\right)
\right\vert ^{2}=\text{ \ \ \ \ \ \ \ \ \ \ \ \ \ \ \ \ \ \ \ \ \ \ \ \ \ }%
\]%
\[
\text{ \ \ \ \ \ \ \ \ \ \ \ \ \ \ \ \ \ \ \ \ \ \ \ \ \ \ }=\prod
\limits_{p\text{ - prime}}\left(  1+\sum\limits_{c\geq1}\sum
\limits_{\substack{1\leq a\leq p^{c}\\a\neq0(\operatorname{mod}p)}}\left\vert
\widehat{\beta_{(p,Q)}}\left(  \frac{a}{p^{c}}\right)  \right\vert
^{2}\right)  ,
\]
where $\widehat{\beta}$ are the Fourier coefficients of $\beta$, defined in
section \ref{LimitPeriodBckgr}.

We then carry out a length calculation of the Fourier coefficients in section
\ref{Fur-coef-comp}, finally ending up with rather complicated expressions
described in Theorem \ref{Fur-coef-summary-th}.

The result is that the mean square is an Euler product%
\[
\lim\limits_{N\rightarrow\infty}\frac{1}{N}\sum\limits_{2\leq n\leq N}%
\beta_{Q}^{2}(n)=\prod\limits_{p\text{ - prime}}M_{p}(Q),
\]
with
\[
M_{p}(Q)=1+\sum\limits_{c\geq1}A_{Q}(p^{c}),
\]
and $A_{Q}(p^{c})$ is given by (\ref{A-Q}), section \ref{mean-qs-calc-sect}.
We evaluate the sum $M_{p}(Q)$ over prime powers as a rational function of $p$
and find that it depends on divisibility of $Q$ by $p$, in particular, for
$p\nmid Q$, $p\neq2$,
\[
M_{p}(Q)=M_{p}(1)=\dfrac{p^{2}(p^{3}+p^{2}-p-3)}{(p^{2}-1)^{2}(p+1)}.
\]
This will prove Theorem \ref{mean-sq-theorem(congr)}.

In section \ref{Quater-sect} we sketch this procedure for the case of the unit
group of the maximal order of a quaternion algebra (theorem
\ref{mean-sq-theorem(quater)}).

\subsection*{\bigskip An application:}

The computation of the mean square is used to study the statistics of
$N_{f}(\tau)$: by Weyl's law we expect $N_{f}(\tau)$ to be asymptotically
equal to a multiple of $\tau/L$. To study the fluctuations around this
expectation one may consider $\tau$ as a random variable.

Then define an averaging operator:
\[
\left\langle F\right\rangle _{T}=\frac{1}{T}%
%TCIMACRO{\dint \limits_{T}^{2T}}%
%BeginExpansion
{\displaystyle\int\limits_{T}^{2T}}
%EndExpansion
F(\tau)d\tau
\]

Studying the moments of $N_{f}$ shows that if $L\rightarrow\infty$, but
$L=o(\ln T)$, then its limiting value distribution is Gaussian with the mean
\[
\left\langle N_{f}(\tau)\right\rangle _{T}=\frac{\tau}{L}\frac{2vol(\Omega
)}{4\pi}%
%TCIMACRO{\dint \limits_{-\infty}^{\infty}}%
%BeginExpansion
{\displaystyle\int\limits_{-\infty}^{\infty}}
%EndExpansion
f(x)dx
\]
and the variance%
\[
\left\langle \left(  N_{f}(\tau)-\left\langle N_{f}(\tau)\right\rangle
_{T}\right)  ^{2}\right\rangle _{T}=\frac{4\kappa}{\pi L}%
%TCIMACRO{\dint \limits_{0}^{\infty}}%
%BeginExpansion
{\displaystyle\int\limits_{0}^{\infty}}
%EndExpansion
\widehat{f}^{2}(u)e^{\pi Lu}du\text{,}%
\]
where $\kappa$ is the mean square of the weighted multiplicities. See
\cite{rud}, \cite{luk}. Thus our computation of $\kappa$ for $\Gamma$ of
congruence type yields the value distribution of $N_{f}$.\ 

\bigskip

\paragraph{Acknowledgement.}

This work was supported in part by Gutwirth Scholarship and by the Israel
Science Foundation founded by the Israel Academy of Sciences and Humanities.
It was carried out as part of the author's PHD thesis at Tel Aviv University,
under the supervision of prof. Ze\'{e}v Rudnick.

\section{Preliminaries on orders in quaternion algebras}

First of all we recall that a ring $B$ with unity is called an algebra of
dimension $n$ over a\ field $F$, if the following three conditions are satisfied:

\mathstrut

\ \ \ \ \ \ $1^{\circ}$ $B\supset F,$and $1_{B}=1_{F}$;

\ \ \ \ \ \ $2^{\circ}$ any element of $F$ commutes with all elements of $B$;

\ \ \ \ \ \ $3^{\circ}$ $B$ is a vector space over $F$ of dimension $n$.

\begin{definition}
Let $B$ be a finite dimensional algebra over $\mathbb{Q}$ with identity
element $1_{B}$. A subset $R\subset B$ is called an order in $B$ if the
following conditions are satisfied:

$\ \ \ \ \ \ \ \ \ 1^{\circ}$ $R$ is a finitely generated $\mathbb{Z}$-module;

$\ \ \ \ \ \ \ \ \ 2^{\circ}$ $R$ contains basis of $B$ over $\mathbb{Q}$;

$\ \ \ \ \ \ \ \ \ 3^{\circ}$ $R$ is a subring of $B$ and $1_{B}\in R$.
\end{definition}

For any finite dimensional algebra $B$ over $\mathbb{Q}$ and for each place
$\upsilon$ of $\mathbb{Q}$ (i.e. $\upsilon=\infty$ or $\upsilon=p$ a prime) we
define $B_{\upsilon}:=B\otimes_{\mathbb{Q}}\mathbb{Q}_{\upsilon}$.
$B_{\upsilon}$ is an algebra over $\mathbb{Q}_{\upsilon}$ of dimension
$\dim_{\mathbb{Q}}B$. For $\upsilon=p$ a prime, orders in $B_{p}$ defined as
in the previous definition, replacing $\mathbb{Q}$ and $\mathbb{Z}$ by
$\mathbb{Q}_{p}$ and $\mathbb{Z}_{p}$. If $R$ is an order of $B$, we get an
order in $B_{p}$ by the definition$\ \ \ \ \ \ \ \ \ \ \ \ \ \ \ \ \ $%
\[
\ R_{p}:=R\otimes_{\mathbb{Z}}\mathbb{Z}_{p}=[\text{the closure of }R~\text{in
}B_{p}].
\]
We also recall \ the definition of adele ring $B_{A}$ \cite[pp. 197-198]{miy}
and \cite[p. 62]{wei}. As a set this is$\ \ \ \ \ \ \ \ \ \ \ $%
\[
\ B_{A}=\{\left(  b_{\upsilon}\right)  \in\prod\limits_{\upsilon}B_{\upsilon
}\mid b_{p}\in R_{p}~\text{for almost all primes }p\},
\]
where $R$ is any order in $B$.\ \ \ \ \ \ \ \ \ \ \ \ \ \ \ \ \ \ \ \ \ \ \ \ \ \ \ \ \ \ \ \ \ \ \ \ 

\begin{definition}
\bigskip A nonzero integer $d$ is called a fundamental discriminant $\ $if
$d\equiv1(\operatorname{mod}4),d$ square free, $d\neq1$; or $d\equiv
0(\operatorname{mod}4)$, $\dfrac{d}{4}\neq1(\operatorname{mod}4)$, $\dfrac
{d}{4}$ squarefree.
\end{definition}

Consider now a quadratic extension of $\mathbb{Q}$ arising as the splitting
field over $\mathbb{Q}$ of a polynomial $P(X)=X^{2}-tX+1,$ $t\in\mathbb{Z}$.
Notice that $P(X)$ is irreducible iff $\sqrt{t^{2}-4}\notin\mathbb{Q}$. In
this case we can write $t^{2}-4=l^{2}d,$ $l\in\mathbb{Z}^{+},$ $d$ is a
fundamental discriminant. The splitting field of $P(X)$ is $\mathbb{Q(}%
\sqrt{d})$ and the zeros of $P(X)$ are%

\[
x_{1}=\frac{t}{2}+l\frac{\sqrt{d}}{2},\ \ \ \ \ \ \ \ \ \ \ \ \ \ \ \ \ x_{2}%
=\frac{t}{2}-l\frac{\sqrt{d}}{2}.
\]
For each fixed fundamental discriminant $d$ we let for each $f\in
\mathbb{Z}^{+}$%

\[
\mathfrak{r}[f]:=\mathbb{Z+}f\omega\mathbb{Z},~~\text{where \ }\omega
=\dfrac{d+\sqrt{d}}{2}.
\]
Then orders in $\mathbb{Q(}\sqrt{d})$ are precisely $\mathfrak{r}[f],$
$f=1,2,3,...$ \cite[pp. 48-49]{cohn}, and in particular $\mathfrak{r}[1]$ is
the unique maximal order.\cite[pp. 146-147]{ded}. The index $\left[
\mathfrak{r}[1]:\mathfrak{r}[f]\right]  =f$, and we call $\mathfrak{r}[f]$ the
order of index $f$ in $\mathbb{Q(}\sqrt{d})$. Note that $\mathfrak{r}%
[f_{1}]\subset\mathfrak{r}[f_{2}]$ iff \ $f_{2}\mid f_{1}$.

For $P(X)$ as above, the roots $x_{1},x_{2}$ of $P(X)$\ satisfy%

\[
\mathbb{Z+}x_{1}\mathbb{Z=Z+}x_{2}\mathbb{Z=}\mathfrak{r}[l]\ \ \ \ \ (t^{2}%
-4=l^{2}d).
\]
Let $K=\mathbb{Q(}\sqrt{d}),$ $d$ as above, then

\begin{center}%
\[
K_{\infty}=\mathbb{Q(}\sqrt{d})\otimes_{\mathbb{Q}}\mathbb{R\simeq}\left\{
\begin{array}
[c]{cc}%
\mathbb{R}\oplus\mathbb{R} & ,d>0\\
\mathbb{C} & ,d<0
\end{array}
\right\}  ,
\]

\end{center}

\[
K_{p}=\mathbb{Q(}\sqrt{d})\otimes_{\mathbb{Q}}\mathbb{Q}_{p}\mathbb{\simeq
}\left\{
\begin{array}
[c]{cc}%
\mathbb{Q}_{p}\oplus\mathbb{Q}_{p} & ,\left(  \frac{d}{p}\right)  =1\\
\mathbb{Q}_{p}\mathbb{(}\sqrt{d}) & ,\left(  \frac{d}{p}\right)  =-1
\end{array}
\right\}  .
\]
The distinct orders in $K_{p}$ are%

\[
\mathfrak{o}[p^{n}]:=\mathbb{Z}_{p}\mathbb{+}p^{n}\omega\mathbb{Z}%
_{p},n=0,1,2,...~\text{and }\omega=\dfrac{d+\sqrt{d}}{2}~\ \text{as above.}%
\]
We have%

\[
\mathfrak{o}[1]\supset\mathfrak{o}[p]\supset\mathfrak{o}[p^{2}]\supset
...~\text{and }\left[  \mathfrak{o}[1]:\mathfrak{o}[p^{n}]\right]  =p^{n}.
\]
Combining the definitions of order in $K$ and in $K_{p}$ we get%

\[
\left(  \mathfrak{r}[f]\right)  _{p}=\mathfrak{o}[p^{n}],~\text{where
}n=ord_{p}f.
\]
For any order $\mathfrak{r}$ in $K=\mathbb{Q(}\sqrt{d})$ we define%
\begin{align*}
\mathfrak{r}_{\infty}  &  :=K_{\infty}\\
\mathfrak{r}_{A}  &  :=\prod\limits_{\upsilon}\mathfrak{r}_{\upsilon}\subset
K_{A}\\
\mathfrak{r}_{\infty^{+}}^{\times}  &  :=\{x\in K_{\infty}^{\times}\mid
N(x)>0\}\\
\mathfrak{r}_{A^{+}}^{\times}  &  :=\{(\alpha_{\upsilon})\in\mathfrak{r}%
_{A}^{\times}\mid N(\alpha_{\infty})>0\}=\mathfrak{r}_{\infty^{+}}^{\times
}\times\prod\limits_{p}\mathfrak{r}_{p}^{\times}\subset\mathfrak{r}%
_{A}^{\times}%
\end{align*}
Here $N$ denotes the extension of the norm $N:K\longrightarrow\mathbb{Q}$ to
the regular norm in the algebra $K_{\infty}$ over $\mathbb{Q}_{\infty
}=\mathbb{R}$. \cite[p.53]{wei}.

The subgroup $\mathfrak{r}_{A^{+}}^{\times}\cdot K^{\times}$ has finite index
in $K_{A}^{\times}$, and we set%

\[
h(\mathfrak{r}):=\left[  K_{A}^{\times}:(\mathfrak{r}_{A^{+}}^{\times}\cdot
K^{\times})\right]  ,
\]
the class number of $\mathfrak{r}$.

For $K=\mathbb{Q(}\sqrt{d})$ as above one can also show \cite[Chapter
5.2]{cohen} that $h(\mathfrak{r}[f])$ is equal to the number of inequivalent
primitive (and if $d<0$, positive) quadratic forms $ax^{2}+bxy+cy^{2}$ with
discriminant $b^{2}-4ac=df^{2}$. In other words $h(\mathfrak{r}[f])=h(df^{2}%
)$, where the right-hand side is as in \cite[vol I, pp. 127-]{land}.

\section{\label{L-func-expr}Expression of weighted multiplicities function in
the terms of Dirichlet's $L$-functions}

In this section we define the weighted multiplicities function $\beta_{Q}(t)$
and express it in the terms of Dirichlet's $L$-functions. It is a necessary
step, which allows to analyze the behavior of \ $\beta_{Q}(t)$. The result is
stated in Theorem \ref{theorem}, proof of which takes the rest of the section.

Let $T\in SL_{2}(\mathbb{Z})$, such that \ $\left\vert t\right\vert
=\left\vert tr(T)\right\vert >2$ . For such $T$ the equation $\det(T-\lambda
I)=0$ has two different solutions: $\lambda_{1,2}=\frac{t\pm\sqrt{t^{2}-4}}%
{2}$, and so $T$ is diagonalizable in $\mathbb{R}$. Such $T$ is called
hyperbolic. Define $\mathcal{N}(T):=\lambda^{2}$, where $\left\vert
\lambda\right\vert >1$. Hence $\mathcal{N}(T)=\frac{\left(  \left\vert
t\right\vert +\sqrt{t^{2}-4}\right)  ^{2}}{4}$.%
\[
\]
Let $Q$ be an odd squarefree integer, and let
\[
\Gamma_{0}(Q)=\left\{  \left[
\begin{array}
[c]{cc}%
a & b\\
c & d
\end{array}
\right]  \in SL_{2}(\mathbb{Z})\mid c\equiv0(\operatorname{mod}Q)\right\}
\]
be the congruence subgroup of $SL_{2}(\mathbb{Z})$. Define$\ $the weighted
multiplicities function%

\[
\ \beta_{Q}(t):=\frac{1}{4}\sum\limits_{\substack{\{T\}\in\Gamma_{0}(Q)
\\\left\vert tr(T)\right\vert =t \\hyperbolic,~doesn^{\prime}t \\fix~cusp}%
}\frac{\ln\mathcal{N}(T_{0})}{\mathcal{N}(T)^{\frac{1}{2}}-\mathcal{N}%
(T)^{-\frac{1}{2}}},
\]
where \ the sum is taken over conjugacy classes of $\ \Gamma_{0}(Q)$.

\begin{theorem}
\label{theorem}$\ $For $2<\left\vert t\right\vert \in\mathbb{Z}$
\[
\beta_{Q}(t)=\sum\limits_{\substack{D,v\geq1\\D\text{ is a discriminant}%
\\Dv^{2}=t^{2}-4}}\frac{1}{v}L(1,\chi_{D})\cdot\prod\limits_{q\mid Q}\left\{
\begin{array}
[c]{cc}%
2, & q^{2}\mid D\\
1+\left(  \dfrac{D}{q}\right)  & q^{2}\nmid D
\end{array}
\right\}  ,
\]
and $\beta_{Q}(t)=0$, otherwise. Here $\left(  \frac{\cdot}{\cdot}\right)  $
is a Legendre's symbol and $\chi_{D}$ is a quadratic character.
\end{theorem}

\[
\]
We start from noting that
\[
\Gamma_{0}(Q)=\bigsqcup\limits_{t\in\mathbb{Z}}H_{t}\text{,}%
\]
where
\[
H_{t}=\left\{  T\in\ \Gamma_{0}(Q)\mid tr(T)=t\right\}  \text{.}%
\]
So we can write ( the change from $1/4$ to $1/2$ comes from not collecting
together $t$ and $-t$ )%

\begin{equation}
\beta_{Q}(t)=\frac{1}{2}\sum\limits_{\substack{T\in H_{t}//\Gamma_{0}(Q)
\\hyperbolic,~doesn^{\prime}t \\fix~cusp}}\frac{\ln\mathcal{N}(T_{0}%
)}{\mathcal{N}(T)^{\frac{1}{2}}-\mathcal{N}(T)^{-\frac{1}{2}}} \label{main12}%
\end{equation}
Here $H_{t}//\Gamma_{0}(Q)$ is the set of $\Gamma_{0}(Q)$-conjugacy classes of
$H_{t}$.

\begin{lemma}
\bigskip\ For a hyperbolic element $T\in H_{t}$ (i.e. $t^{2}-4>0$) no fixed
point of $T$ is a cusp of \ $\Gamma_{0}(Q).$
\end{lemma}

\begin{proof}
Note, that the cusps of $\Gamma_{0}(Q)$ are exactly the points in $\left\{
i\infty\right\}  \cup\mathbb{Q}$, by \cite[cor. 1.5.5, th. 4.1.3(2)]{miy}.
Assume $T=\left[
\begin{array}
[c]{cc}%
a & b\\
c & d
\end{array}
\right]  $. First, $c\neq0$, since we will get $T$ parabolic. So, the
fixpoints of $T$ are the two real solutions to $\dfrac{ax+b}{cx+d}%
=x\Leftrightarrow x=\dfrac{a-d}{2c}\pm\dfrac{\sqrt{t^{2}-4}}{2\left\vert
c\right\vert }$, which are both irrational, since $\sqrt{t^{2}-4}%
\notin\mathbb{Q}$.
\end{proof}

\[
\]
Now we want to know \ when in\ the sum (\ref{main12}), $H_{t}$ is not the
empty set.

\begin{lemma}
\bigskip\ Let \ $t\in\mathbb{Z}$, such that $\sqrt{t^{2}-4}\notin\mathbb{Q}$ ,
and write $t^{2}-4=l^{2}d$, $l\in\mathbb{Z}^{+},~d$ is a fundamental
discriminant. Then $H_{t}$ is not empty iff for all primes $q$ divide $Q$ we
have $q\mid l$ or $\left(  \dfrac{d}{q}\right)  \neq-1.$
\end{lemma}

\begin{proof}
$H_{t}\neq\emptyset$ iff there some $a,b,c\in\mathbb{Z}$, such that
$\det\left[
\begin{array}
[c]{cc}%
a & b\\
Qc & t-a
\end{array}
\right]  =1.$ Hence $H_{t}\neq\emptyset$ iff there is some $a\in\mathbb{Z}$,
such that $a(t-a)\equiv1(\operatorname{mod}Q)\Leftrightarrow(2a-t)^{2}\equiv
t^{2}-4(\operatorname{mod}Q)\Leftrightarrow(2a-t)^{2}\equiv l^{2}%
d(\operatorname{mod}Q).$ Since $Q$ is odd squarefree, the last congruence is
solvable iff $(2a-t)^{2}\equiv l^{2}d(\operatorname{mod}q)$ for all primes
$q\mid Q.$ Or equivalently: the last congruence is solvable iff $q\mid l$ or
$\left(  \dfrac{d}{q}\right)  \neq-1,$ for all $q\mid Q.$
\end{proof}

\[
\]
Note, that for $T\in H_{t}$
\[
\mathcal{N}(T)^{\frac{1}{2}}-\mathcal{N}(T)^{-\frac{1}{2}}=\dfrac{\left\vert
t\right\vert +\sqrt{t^{2}-4}}{2}-\dfrac{2}{\left\vert t\right\vert
+\sqrt{t^{2}-4}}=\sqrt{t^{2}-4}.
\]
So, for fixed $t\in\mathbb{Z}$, such that $\sqrt{t^{2}-4}\notin\mathbb{Q}$,
$t^{2}-4>0$, $H_{t}\neq\emptyset$ we can write (\ref{main12}) in the form%

\begin{equation}
\beta_{Q}(t)=\frac{1}{2}\sum\limits_{T\in H_{t}//\Gamma_{0}(Q)}\frac
{\ln\mathcal{N}(T_{0})}{\sqrt{t^{2}-4}} \label{main2}%
\end{equation}
\ We will enumerate the $\Gamma_{0}(Q)$-conjugacy classes in $H_{t}$, and for
each conjugacy class will get the number $\mathcal{N}(T_{0}).$

Fix some $T_{t}\in H_{t}$, and take the polynomial \ $P(X)=X^{2}-tX+1$, i.e.
the characteristic polynomial of $T_{t}$, and so $P(T_{t})=0.$ $P(X)$ is
irreducible, since $\sqrt{t^{2}-4}\notin\mathbb{Q}$. Hence $\mathbb{Q}\left[
T_{t}\right]  \cong\mathbb{Q(}\sqrt{d})$, where $t^{2}-4=l^{2}d.$

Define the set $C(T_{t}):=\{\delta T_{t}\delta^{-1}\mid\delta\in
GL_{2}(\mathbb{Q)\}}$. Then $H_{t}=R\cap C(T_{t})$, where $R=\left\{  \left[
\begin{array}
[c]{cc}%
a & b\\
c & d
\end{array}
\right]  \in M_{2}(\mathbb{Z)\mid}\text{ \ \ }Q\mid c\right\}  $ is an order
in $M_{2}(\mathbb{Q})$.

\bigskip For any $\delta\in GL_{2}(\mathbb{Q)}$, \ $\left(  \mathbb{Q[}%
T_{t}]\cap\delta R\delta^{-1}\right)  $ is an order in $\mathbb{Q[}T_{t}],$ by
the fact, that if $A^{\prime}$ is a subalgebra of $A$, and $\mathcal{O}$ is an
order in $A$, then $A^{\prime}\cap\mathcal{O}$ is an order in $A^{\prime}$.
For any order $\mathfrak{r}$ in $\mathbb{Q[}T_{t}]$ define $C(T_{t}%
,\mathfrak{r}):=\{\delta T_{t}\delta^{-1}\mid\delta\in GL_{2}(\mathbb{Q)}$,
$\mathbb{Q[}T_{t}]\cap\delta^{-1}R\delta=\mathfrak{r}\mathbb{\}}$. Following
\cite{stromb}, \cite{miy} we have a

\begin{lemma}
\mathstrut%
\[
\bigskip C(T_{t})=\bigsqcup\limits_{f=1}^{\infty}C(T_{t},\mathfrak{r}[f]),
\]%
\[
H_{t}=\bigsqcup\limits_{f\mid l}C(T_{t},\mathfrak{r}[f]),
\]
and each $C(T_{t},\mathfrak{r}[f])$ is closed under $\Gamma_{0}(Q)$%
-conjugation.\footnote{$\mathfrak{r}[f]:=\mathbb{Z+}fw\mathbb{Z}$, is an all
orders in $\mathbb{Q(}\sqrt{d}\mathbb{)}$, where $f\in\mathbb{Z}^{+},$
$w=\frac{d+\sqrt{d}}{2}.$ Since $\mathbb{Q(}\sqrt{d}\mathbb{)\cong Q[}T_{t}]$
we use the symbol $\mathfrak{r}[f]$ for the corresponding order in
$\mathbb{Q[}T_{t}]$ .}
\end{lemma}

\begin{proof}
Clearly, $C(T_{t})=\bigcup\limits_{f=1}^{\infty}C(T_{t},\mathfrak{r}[f])$,
since there are no other orders than \ $\mathfrak{r}[1]$, $\mathfrak{r}[2]$,
... in $\mathbb{Q[}T_{t}]$. To prove disjointedness, we must show that if
$\delta_{1}T_{t}\delta_{1}^{-1}=\delta_{2}T_{t}\delta_{2}^{-1}$ $(\delta
_{1},\delta_{2}\in GL_{2}(\mathbb{Q))}$, then $\delta_{1}^{-1}R\delta_{1}$,
and $\delta_{2}^{-1}R\delta_{2}$ have the same intersection with
$\mathbb{Q[}T_{t}]$. From $\delta_{1}T_{t}\delta_{1}^{-1}=\delta_{2}%
T_{t}\delta_{2}^{-1}$ we get $\delta_{1}^{-1}\delta_{2}\in\mathbb{Q[}T_{t}]$,
since if an element of $M_{2}(\mathbb{Q)}$ commutes with $T_{t}$, it is in
$\mathbb{Q[}T_{t}]$, by \cite[Lemma 5.2.2(3)]{miy}, and thus
\[
\mathbb{Q[}T_{t}]\cap\delta_{2}^{-1}R\delta_{2}=(\delta_{1}^{-1}\delta
_{2})\left(  \mathbb{Q[}T_{t}]\cap\delta_{2}^{-1}R\delta_{2}\right)
(\delta_{1}^{-1}\delta_{2})^{-1}=\mathbb{Q[}T_{t}]\cap\delta_{1}^{-1}%
R\delta_{1}.
\]
This prove the first relation.

Next, for any order $\mathfrak{r}$ in $\mathbb{Q[}T_{t}]$ and any $\delta\in
GL_{2}(\mathbb{Q)}$, such that $\delta T_{t}\delta^{-1}\in C(T_{t}%
,\mathfrak{r}) $, we have :
\[
T_{t}\in\mathfrak{r}\Leftrightarrow T_{t}\in\mathbb{Q[}T_{t}]\cap\delta
^{-1}R\delta\Leftrightarrow T_{t}\in\delta^{-1}R\delta\Leftrightarrow\delta
T_{t}\delta^{-1}\in R\Leftrightarrow\delta T_{t}\delta^{-1}\in H_{t}.
\]

In other words :

\hspace{1in}if $T_{t}\in\mathfrak{r}$, then $C(T_{t},\mathfrak{r})\subset
H_{t};$

\hspace{1in}if $T_{t}\notin\mathfrak{r}$, then $C(T_{t},\mathfrak{r})\cap
H_{t}=\emptyset.$

\mathstrut

Note that $T_{t}\in\mathfrak{r}[l],$ where $t^{2}-4=l^{2}d$, $l\in
\mathbb{Z}^{+}.$So the orders of $\mathfrak{r}$ in $\mathbb{Q[}T_{t}]$ that
contain $T_{t}$ are exactly $\mathfrak{r}[f],$ such that $\mathfrak{r}%
[f]\supset\mathfrak{r}[l],$ that is, such that $f\mid l$. By definition
$H_{t}\subset C(T_{t})$, \ and so from the first relation the second relation follows.

Finally, $\gamma^{-1}R\gamma=R$, for any $\gamma\in R^{\times}$. So we get
that each $C(T_{t},\mathfrak{r})$ is closed under $R^{\times}$-conjugation,
and in particular under $\Gamma_{0}(Q)$-conjugation.
\end{proof}

\mathstrut

Let $\varepsilon_{d}=\dfrac{x+y\sqrt{d}}{2}$ be the proper fundamental unit in
$\mathbb{Q(}\sqrt{d})$ ( $(x,y)$ is the positive integer solution to
$x^{2}-dy^{2}=4,$ for which $y>0$ is minimal). Define
\[
\mathfrak{r}[f]^{1}:=\{\alpha\in\mathfrak{r}[f]\mid N(\alpha)=1\},
\]
the units of the order $\mathfrak{r}[f].$ Since
\[
\mathfrak{r}[1]^{1}=\left\{  \pm\varepsilon_{d}^{k}\mid k\in\mathbb{Z}%
\right\}  ,
\]
we have
\[
\mathfrak{r}[f]^{1}=\left\{  \pm\left(  \varepsilon_{d}^{\left[
\mathfrak{r}[1]^{1}:\mathfrak{r}[f]^{1}\right]  }\right)  ^{k}\mid
k\in\mathbb{Z}\right\}  .
\]

\begin{lemma}
If $T$ is hyperbolic, i.e. $d>0$, then $\mathcal{N}(T_{0})=\varepsilon
_{d}^{2\left[  \mathfrak{r}[1]^{1}:\mathfrak{r}[f]^{1}\right]  }.$
\end{lemma}

\begin{proof}
Choose some $\delta\in GL_{2}(\mathbb{Q)}$, such that $T=\delta T_{t}%
\delta^{-1}$, and $\mathbb{Q[}T_{t}]\cap\delta R\delta^{-1}=\mathfrak{r} $.
Then $Z_{\Gamma_{0}(Q)}(T),$ the centralizer of $T$ in $\Gamma_{0}(Q)$ is
\[
Z_{\Gamma_{0}(Q)}(T)=\mathbb{Q[}T]\cap\Gamma_{0}(Q)=\delta(\mathbb{Q[}%
T_{t}]\cap\delta^{-1}\Gamma_{0}(Q)\delta)\delta^{-1}=\delta\mathfrak{r}%
[f]\delta^{-1}%
\]
by \cite[lemma 5.2.2(3)]{miy}.Since
\[
\mathfrak{r}[f]^{1}=\left\{  \pm\left(  \varepsilon_{d}^{\left[
\mathfrak{r}[1]^{1}:\mathfrak{r}[f]^{1}\right]  }\right)  ^{k}\mid
k\in\mathbb{Z}\right\}  ,
\]
for $d>0$ we can take for $T_{0\text{ }}$ the image of $\varepsilon
_{d}^{\left[  \mathfrak{r}[1]^{1}:\mathfrak{r}[f]^{1}\right]  }$ under \ the
composite of the two isomorphisms $\mathbb{Q}\left[  T_{t}\right]
\cong\mathbb{Q(}\sqrt{d})$ and $\mathbb{Q}\left[  T_{t}\right]  \ni
S\longmapsto\delta S\delta^{-1}$\bigskip$\in\mathbb{Q}\left[  T\right]  .$
Then
\[
tr_{M_{2}(\mathbb{Q)}}(T_{0})=tr\left(  \varepsilon_{d}^{\left[
\mathfrak{r}[1]^{1}:\mathfrak{r}[f]^{1}\right]  }\right)  .
\]
If we will denote $\varepsilon_{d}^{\left[  \mathfrak{r}[1]^{1}:\mathfrak{r}%
[f]^{1}\right]  }=\alpha+\beta\sqrt{d}\in\mathbb{Q(}\sqrt{d})$, note, that
$\alpha,\beta\geq0,$ and $N(\alpha+\beta\sqrt{d})=\alpha^{2}-d\beta^{2}=1.$ We
have $tr(\alpha+\beta\sqrt{d})=t,$ and so
\begin{align*}
\mathcal{N}(T_{0})  &  =\frac{\left(  \left\vert t\right\vert +\sqrt{t^{2}%
-4}\right)  ^{2}}{4}=\frac{\left(  2\left\vert \alpha\right\vert
+\sqrt{4\alpha^{2}-4}\right)  ^{2}}{4}=\\
&  =\frac{\left(  2\left\vert \alpha\right\vert +\sqrt{4(1+d\beta^{2}%
)-4}\right)  ^{2}}{4}=\frac{\left(  2\left\vert \alpha\right\vert +2\left\vert
\beta\right\vert \sqrt{d}\right)  ^{2}}{4}=\\
&  =\left(  \left\vert \alpha\right\vert +\left\vert \beta\right\vert \sqrt
{d}\right)  ^{2}=\left(  \alpha+\beta\sqrt{d}\right)  ^{2}=\varepsilon
_{d}^{2\left[  \mathfrak{r}[1]^{1}:\mathfrak{r}[f]^{1}\right]  }.
\end{align*}

\end{proof}%

\[
\]
Now we can rewrite (\ref{main2}) in the form%
\[
\beta_{Q}(t)=\sum\limits_{f\mid l}\left\vert C(T_{t},\mathfrak{r}%
[f])//\Gamma_{0}(Q)\right\vert \frac{\ln\varepsilon_{d}^{\left[
\mathfrak{r}[1]^{1}:\mathfrak{r}[f]^{1}\right]  }}{\sqrt{t^{2}-4}},
\]
where we write $t^{2}-4=dl^{2}$ with $d$ a fundamental discriminant and
$l\geq1$.

For the quantities $\left\vert C(T_{t},\mathfrak{r}[f])//\Gamma_{0}%
(Q)\right\vert $ \ we refer to the \cite[\S 6.6]{miy}, \cite{stromb}:

\begin{proposition}
\label{prop}Let $Q$ be squarefree and $h(df^{2})$ the narrow class number of
$\mathbb{Q(}\sqrt{d})$. Then
\[
\left\vert C(T_{t},\mathfrak{r}[f])//\Gamma_{0}(Q)\right\vert =h(df^{2}%
)\cdot\left\{
\begin{array}
[c]{cc}%
2, & \text{if }d<0\\
1, & \text{if }d>0
\end{array}
\right\}  \cdot\prod\limits_{p\mid Q}\left\{
\begin{array}
[c]{cc}%
2, & p\mid f\\
1+\left(  \frac{d}{p}\right)  & p\nmid f
\end{array}
\right\}  .
\]

\end{proposition}

\begin{corollary}
$\ \ $\ $\ \ \ \ \ \ \ \ \ \ \ $%
\[
\beta_{Q}(t)=%
%TCIMACRO{\dsum \limits_{f\mid l}}%
%BeginExpansion
{\displaystyle\sum\limits_{f\mid l}}
%EndExpansion
h(df^{2})\dfrac{\ln\varepsilon_{d}^{\left[  \mathfrak{r}[1]^{1}:\mathfrak{r}%
[f]^{1}\right]  }}{\sqrt{t^{2}-4}}\cdot%
%TCIMACRO{\dprod \limits_{q\mid Q}}%
%BeginExpansion
{\displaystyle\prod\limits_{q\mid Q}}
%EndExpansion
\left\{
\begin{array}
[c]{cc}%
2, & q\mid f\\
1+\left(  \dfrac{d}{q}\right)  & q\nmid f
\end{array}
\right\}  ,
\]
where we write $t^{2}-4=dl^{2}$ with $d$ a fundamental discriminant and
$l\geq1$.
\end{corollary}

\begin{proof}
\bigskip Immediately from proposition \ref{prop}.
\end{proof}

\[
\]
We can continue the process%

\[
\beta_{Q}(t)=%
%TCIMACRO{\dsum \limits_{f\mid l}}%
%BeginExpansion
{\displaystyle\sum\limits_{f\mid l}}
%EndExpansion
h(df^{2})\dfrac{\ln\varepsilon_{d}^{\left[  \mathfrak{r}[1]^{1}:\mathfrak{r}%
[f]^{1}\right]  }}{\sqrt{t^{2}-4}}\cdot%
%TCIMACRO{\dprod \limits_{q\mid Q}}%
%BeginExpansion
{\displaystyle\prod\limits_{q\mid Q}}
%EndExpansion
\left\{
\begin{array}
[c]{cc}%
2, & q\mid f\\
1+\left(  \dfrac{d}{q}\right)  & q\nmid f
\end{array}
\right\}  =
\]

\[
=%
%TCIMACRO{\dsum \limits_{f\mid l}}%
%BeginExpansion
{\displaystyle\sum\limits_{f\mid l}}
%EndExpansion
\dfrac{h(df^{2})\ln\varepsilon_{df^{2}}}{l\sqrt{d}}\cdot%
%TCIMACRO{\dprod \limits_{q\mid Q}}%
%BeginExpansion
{\displaystyle\prod\limits_{q\mid Q}}
%EndExpansion
\left\{
\begin{array}
[c]{cc}%
2, & q\mid f\\
1+\left(  \dfrac{d}{q}\right)  & q\nmid f
\end{array}
\right\}  =
\]

\[
=%
%TCIMACRO{\dsum \limits_{\substack{f\mid l \\(D=df^{2})}}}%
%BeginExpansion
{\displaystyle\sum\limits_{\substack{f\mid l \\(D=df^{2})}}}
%EndExpansion
\dfrac{h(D)\ln\varepsilon_{D}}{l\cdot\frac{\sqrt{D}}{f}}\cdot%
%TCIMACRO{\dprod \limits_{q\mid Q}}%
%BeginExpansion
{\displaystyle\prod\limits_{q\mid Q}}
%EndExpansion
\left\{
\begin{array}
[c]{cc}%
2, & q\mid f\\
1+\left(  \dfrac{d}{q}\right)  & q\nmid f
\end{array}
\right\}  .
\]
Now using Dirichlet's class number formula \ $h(D)\ln\varepsilon_{D}=\sqrt
{D}L(1,\chi_{D}),$ we get$\ $%
\[
\ \beta_{Q}(t)=%
%TCIMACRO{\dsum \limits_{\substack{f\mid l \\(D=df^{2})}}}%
%BeginExpansion
{\displaystyle\sum\limits_{\substack{f\mid l \\(D=df^{2})}}}
%EndExpansion
L(1,\chi_{D})\dfrac{f}{l}\cdot%
%TCIMACRO{\dprod \limits_{q\mid Q}}%
%BeginExpansion
{\displaystyle\prod\limits_{q\mid Q}}
%EndExpansion
\left\{
\begin{array}
[c]{cc}%
2, & q\mid f\\
1+\left(  \dfrac{D}{q}\right)  & q\nmid f
\end{array}
\right\}  =
\]

\[
=%
%TCIMACRO{\dsum \limits_{\substack{D,v\geq1 \\Dv^{2}=t^{2}-4}}}%
%BeginExpansion
{\displaystyle\sum\limits_{\substack{D,v\geq1 \\Dv^{2}=t^{2}-4}}}
%EndExpansion
\dfrac{1}{v}L(1,\chi_{D})\cdot%
%TCIMACRO{\dprod \limits_{q\mid Q}}%
%BeginExpansion
{\displaystyle\prod\limits_{q\mid Q}}
%EndExpansion
\left\{
\begin{array}
[c]{cc}%
2, & q\mid f\\
1+\left(  \dfrac{D}{q}\right)  & q\nmid f
\end{array}
\right\}  =
\]

\[
=%
%TCIMACRO{\dsum \limits_{\substack{D,v\geq1\\Dv^{2}=t^{2}-4}}}%
%BeginExpansion
{\displaystyle\sum\limits_{\substack{D,v\geq1\\Dv^{2}=t^{2}-4}}}
%EndExpansion
\dfrac{1}{v}L(1,\chi_{D})\cdot%
%TCIMACRO{\dprod \limits_{q\mid Q}}%
%BeginExpansion
{\displaystyle\prod\limits_{q\mid Q}}
%EndExpansion
\left\{
\begin{array}
[c]{cc}%
2, & q^{2}\mid D\\
1+\left(  \dfrac{D}{q}\right)  & q^{2}\nmid D
\end{array}
\right\}  ,
\]
and the theorem \ref{theorem} follows.

\begin{remark}
Here $D$ is a discriminant, i.e. $D\equiv0,1(\operatorname{mod}4)$. We assume
it from now on.
\end{remark}

\section{Factorization of weighted multiplicities}

\bigskip Here we factorize (lemma \ref{mult}) the weighted multiplicities
function as a finite products of a local terms, defined below.

We have \ \ %

\[
\beta_{Q}(n)=%
%TCIMACRO{\dsum \limits_{\substack{D,v\geq1\\Dv^{2}=n^{2}-4}}}%
%BeginExpansion
{\displaystyle\sum\limits_{\substack{D,v\geq1\\Dv^{2}=n^{2}-4}}}
%EndExpansion
\dfrac{1}{v}L(1,\chi_{D})\cdot%
%TCIMACRO{\dprod \limits_{q\mid Q}}%
%BeginExpansion
{\displaystyle\prod\limits_{q\mid Q}}
%EndExpansion
\left\{
\begin{array}
[c]{cc}%
2, & q^{2}\mid D\\
1+\left(  \dfrac{D}{q}\right)  & q^{2}\nmid D
\end{array}
\right\}  .
\]
We use now the Euler product formula for $L(1,\chi_{D})$. For $n\geq3,P\geq
Q,$ define%

\[
\beta_{P,Q}(n):=%
%TCIMACRO{\dsum \limits_{\substack{_{\substack{D,v\geq1 \\Dv^{2}=n^{2}-4}}
%\\p\mid v\Longrightarrow p\leq P}}}%
%BeginExpansion
{\displaystyle\sum\limits_{\substack{_{\substack{D,v\geq1 \\Dv^{2}=n^{2}-4}}
\\p\mid v\Longrightarrow p\leq P}}}
%EndExpansion
\left(  \dfrac{1}{v}%
%TCIMACRO{\dprod \limits_{p\leq P}}%
%BeginExpansion
{\displaystyle\prod\limits_{p\leq P}}
%EndExpansion
\left(  1-\dfrac{\chi_{D}(p)}{p}\right)  ^{-1}\right)  \cdot%
%TCIMACRO{\dprod \limits_{q\mid Q}}%
%BeginExpansion
{\displaystyle\prod\limits_{q\mid Q}}
%EndExpansion
\left\{
\begin{array}
[c]{cc}%
2, & q^{2}\mid D\\
1+\left(  \dfrac{D}{q}\right)  & q^{2}\nmid D
\end{array}
\right\}  .
\]
Note, that \ \ %

\[
\ q^{2}\mid D\Longleftrightarrow q^{2}\mid(n^{2}-4)v^{-2}\Longleftrightarrow
q^{2}\mid(n^{2}-4)\prod\limits_{\substack{p\leq P \\p\mid v}}p^{-2b_{p}}%
\]

\[
\Longleftrightarrow%
%TCIMACRO{\QATOPD{\{}{\}}{q^{2}\mid n^{2}-4\text{ \ \ \ \ \ \ \ \ \ }%
%,\text{\ }b_{q}=0}{q^{2}\mid(n^{2}-4)q^{-2b_{q}},\text{ }b_{q}\neq0}}%
%BeginExpansion
\genfrac{\{}{\}}{0pt}{}{q^{2}\mid n^{2}-4\text{ \ \ \ \ \ \ \ \ \ }%
,\text{\ }b_{q}=0}{q^{2}\mid(n^{2}-4)q^{-2b_{q}},\text{ }b_{q}\neq0}%
%EndExpansion
\Longleftrightarrow q^{2}\mid(n^{2}-4)q^{-2ord_{q}v}.
\]
Hence\ $\ \ \ \ \ \ \ \ \ \ \ \ \ \ $%
\[
\ \ \beta_{P,Q}(n)=%
%TCIMACRO{\dsum \limits_{\substack{_{\substack{D,v\geq1 \\Dv^{2}=n^{2}-4}}
%\\p\mid v\Longrightarrow p\leq P}}}%
%BeginExpansion
{\displaystyle\sum\limits_{\substack{_{\substack{D,v\geq1 \\Dv^{2}=n^{2}-4}}
\\p\mid v\Longrightarrow p\leq P}}}
%EndExpansion
\left(  \dfrac{1}{v}%
%TCIMACRO{\dprod \limits_{p\leq P}}%
%BeginExpansion
{\displaystyle\prod\limits_{p\leq P}}
%EndExpansion
\left(  1-\dfrac{\chi_{D}(p)}{p}\right)  ^{-1}\right)  \times
\]

\[
\times%
%TCIMACRO{\dprod \limits_{q\mid Q}}%
%BeginExpansion
{\displaystyle\prod\limits_{q\mid Q}}
%EndExpansion
\left\{
\begin{array}
[c]{cc}%
2 & \text{\ },\text{\ }q^{2}\mid(n^{2}-4)q^{-2ord_{q}v}\\
1+\left(  \dfrac{(n^{2}-4)q^{-2ord_{q}v}}{q}\right)  & \text{\ },\text{ }%
q^{2}\nmid(n^{2}-4)q^{-2ord_{q}v}%
\end{array}
\right\}  \ .
\]
Define \ a function%

\[
\beta_{(p,Q)}(n):=\sum\limits_{b\geq0}\frac{1}{p^{b}}\left(  1-\frac{1}{p}%
\chi_{(n^{2}-4)p^{-2b}}(p)\right)  ^{-1}\cdot\mathbb{I}_{p^{b}}(n)\ ,
\]
where%

\[
\mathbb{I}_{p^{b}}(n):=\left\{
\begin{array}
[c]{cc}%
1 & \text{,}n^{2}=4(\operatorname{mod}p^{2b})\\
0 & \text{, else}%
\end{array}
\right\}  \ ,\text{ for }p\neq2,p\nmid Q
\]

\[
\mathbb{I}_{2^{b}}(n):=\left\{
\begin{array}
[c]{ccc}%
1 & \text{, }n^{2}=4(\operatorname{mod}\text{ }2^{2b}) & \text{, }%
(n^{2}-4)2^{-2b}\text{ is a discriminant}\\
0 & \text{, else} &
\end{array}
\right\}  ,
\]
and for $q\mid Q$%

\[
\mathbb{I}_{q^{b}}(n):=\left\{
\begin{array}
[c]{ccc}%
2 & \text{, }n^{2}=4(\operatorname{mod}\text{ }q^{2b}) & \text{, }q^{2}%
\mid(n^{2}-4)q^{-2b}\\
1+\left(  \frac{(n^{2}-4)q^{-2b}}{q}\right)  & \text{, }n^{2}%
=4(\operatorname{mod}\text{ }q^{2b}) & \text{, }q^{2}\nmid(n^{2}-4)q^{-2b}\\
0 & \text{, else} &
\end{array}
\right\}  .
\]

\begin{lemma}
\label{mult}For $Q$ odd squarefree and $P\geq Q$ we have
$\ \ \ \ \ \ \ \ \ \ \ \ \ \ \ \ \ \ \ \ \ \ $%
\[
\beta_{P,Q}(n)=\prod\limits_{\substack{p\leq P \\p\nmid Q}}\beta
_{(p,Q)}(n)\cdot\prod\limits_{q\mid Q}\beta_{(q,Q)}(n).
\]

\end{lemma}

\begin{proof}%
\[%
%TCIMACRO{\dprod \limits_{\substack{p\leq P \\p\nmid Q}}}%
%BeginExpansion
{\displaystyle\prod\limits_{\substack{p\leq P \\p\nmid Q}}}
%EndExpansion
\beta_{(p,Q)}(n)\cdot%
%TCIMACRO{\dprod \limits_{q\mid Q}}%
%BeginExpansion
{\displaystyle\prod\limits_{q\mid Q}}
%EndExpansion
\beta_{(q,Q)}(n)=%
%TCIMACRO{\dprod \limits_{\substack{p\leq P \\p\nmid Q}}}%
%BeginExpansion
{\displaystyle\prod\limits_{\substack{p\leq P \\p\nmid Q}}}
%EndExpansion
\left(
%TCIMACRO{\dsum \limits_{b\geq0}}%
%BeginExpansion
{\displaystyle\sum\limits_{b\geq0}}
%EndExpansion
\dfrac{1}{p^{b}}\left(  1-\dfrac{1}{p}\chi_{(n^{2}-4)p^{-2b}}(p)\right)
^{-1}\cdot\mathbb{I}_{p^{b}}(n)\right)  \times
\]%
\[%
%TCIMACRO{\dprod \limits_{q\mid Q}}%
%BeginExpansion
{\displaystyle\prod\limits_{q\mid Q}}
%EndExpansion
\left(
%TCIMACRO{\dsum \limits_{b\geq0}}%
%BeginExpansion
{\displaystyle\sum\limits_{b\geq0}}
%EndExpansion
\dfrac{1}{q^{b}}\left(  1-\dfrac{1}{q}\chi_{(n^{2}-4)q^{-2b}}(q)\right)
^{-1}\mathbb{I}_{q^{b}}(n)\right)  .
\]
After opening the brackets, we will get the sum of \ terms of the
form:\medskip\medskip\medskip%
\[
\frac{1}{p_{1}^{b_{1}}}\left(  1-\frac{1}{p_{1}}\chi_{(n^{2}-4)p_{1}^{-2b_{1}%
}}(p_{1})\right)  ^{-1}\mathbb{I}_{p_{1}^{b_{1}}}(n)\cdot\frac{1}{p_{2}%
^{b_{2}}}\left(  1-\frac{1}{p_{2}}\chi_{(n^{2}-4)p_{2}^{-2b_{2}}}%
(p_{2})\right)  ^{-1}\mathbb{I}_{p_{2}^{b_{2}}}(n)\cdot\ldots\cdot
\]%
\[
\frac{1}{p_{k}^{b_{k}}}\left(  1-\frac{1}{p_{k}}\chi_{(n^{2}-4)p_{k}^{-2b_{k}%
}}(p_{k})\right)  ^{-1}\mathbb{I}_{p_{k}^{b_{k}}}(n)\cdot\frac{1}%
{q_{1}^{b_{q_{1}}}}\left(  1-\frac{1}{q_{1}}\chi_{(n^{2}-4)q_{1}^{-2b_{q_{1}}%
}}(q_{1})\right)  ^{-1}\mathbb{I}_{q_{1}^{b_{q_{1}}}}(n)\cdot\ldots\cdot
\]%
\[
\frac{1}{q_{l}^{b_{q_{l}}}}\left(  1-\frac{1}{q_{l}}\chi_{(n^{2}%
-4)q_{l}^{-2b_{q_{l}}}}(q_{l})\right)  ^{-1}\mathbb{I}_{q_{l}^{b_{q_{l}}}%
}(n).
\]
Therefore, since $P\geq Q$ we have

\mathstrut

$%
%TCIMACRO{\dprod \limits_{\substack{p\leq P \\p\nmid Q}}}%
%BeginExpansion
{\displaystyle\prod\limits_{\substack{p\leq P \\p\nmid Q}}}
%EndExpansion
\beta_{(p,Q)}(n)\cdot%
%TCIMACRO{\dprod \limits_{q\mid Q}}%
%BeginExpansion
{\displaystyle\prod\limits_{q\mid Q}}
%EndExpansion
\beta_{(q,Q)}(n)=$

\mathstrut

$%
%TCIMACRO{\dsum \limits_{b\geq0}}%
%BeginExpansion
{\displaystyle\sum\limits_{b\geq0}}
%EndExpansion%
%TCIMACRO{\dprod \limits_{\substack{_{\substack{p\leq P \\n^{2}%
%=4(\operatorname{mod}\text{ }p^{2b})}} \\(n^{2}-4)2^{-2b}\text{ is } \\\text{a
%discriminant}}}}%
%BeginExpansion
{\displaystyle\prod\limits_{\substack{_{\substack{p\leq P \\n^{2}%
=4(\operatorname{mod}\text{ }p^{2b})}} \\(n^{2}-4)2^{-2b}\text{ is } \\\text{a
discriminant}}}}
%EndExpansion
\dfrac{1}{p^{b}}\left(  1-\dfrac{1}{p}\chi_{(n^{2}-4)p^{-2b}}(p)\right)
^{-1}\times$

$\ \ \ \ \ \ \ \ \ \ \ \ \ \ \ \ \ \ \ \ \ \ \ \ \times%
%TCIMACRO{\dprod \limits_{q\mid Q}}%
%BeginExpansion
{\displaystyle\prod\limits_{q\mid Q}}
%EndExpansion
\left\{
\begin{array}
[c]{cc}%
2 & \text{\ },\text{\ }q^{2}\mid(n^{2}-4)q^{-2ord_{q}v}\\
1+\left(  \dfrac{(n^{2}-4)q^{-2ord_{q}v}}{q}\right)  & \text{\ },\text{ }%
q^{2}\nmid(n^{2}-4)q^{-2ord_{q}v}%
\end{array}
\right\}  =$

\mathstrut\mathstrut

\mathstrut

$%
%TCIMACRO{\dsum \limits_{\substack{D,v\geq1 \\Dv^{2}=n^{2}-4 \\p\mid
%v\Longrightarrow p\leq P}}}%
%BeginExpansion
{\displaystyle\sum\limits_{\substack{D,v\geq1 \\Dv^{2}=n^{2}-4 \\p\mid
v\Longrightarrow p\leq P}}}
%EndExpansion
\left[  \dfrac{1}{v}%
%TCIMACRO{\dprod \limits_{p\leq P}}%
%BeginExpansion
{\displaystyle\prod\limits_{p\leq P}}
%EndExpansion
\left(  1-\dfrac{1}{p}\chi_{D}(p)\right)  ^{-1}\right]  \times$

\mathstrut

$\times%
%TCIMACRO{\dprod \limits_{q\mid Q}}%
%BeginExpansion
{\displaystyle\prod\limits_{q\mid Q}}
%EndExpansion
\left\{
\begin{array}
[c]{cc}%
2 & \text{\ },\text{\ }q^{2}\mid(n^{2}-4)q^{-2ord_{q}v}\\
1+\left(  \dfrac{(n^{2}-4)q^{-2ord_{q}v}}{q}\right)  & \text{\ },\text{ }%
q^{2}\nmid(n^{2}-4)q^{-2ord_{q}v}%
\end{array}
\right\}  =\beta_{P,Q}(n).$
\end{proof}

\section{\label{LimitPeriodBckgr}Limit periodic functions and Fourier
analysis}

Let $s\geq1$. For $f:\mathbb{N}\mathbf{\longrightarrow}\mathbb{C}$, define the seminorm%

\[
\left\Vert f\right\Vert _{s}:=\left(  \limsup\limits_{N\rightarrow\infty}%
\frac{1}{N}\sum\limits_{1\leq n\leq N}\left\vert f(n)\right\vert ^{s}\right)
^{1/s}\in\lbrack0,\infty).
\]
A function $f$ is called $s$-limit periodic if for every $\varepsilon>0$ there
is a periodic function $h$ with $\left\Vert f-h\right\Vert _{s}\leq
\varepsilon.$ The set $\mathcal{D}^{s}$ of all $s$-limit periodic functions
becomes a Banach space with \ \ norm $\left\Vert \cdot\right\Vert _{s}$ if
functions $f_{1}$, $f_{2}$ with $\left\Vert f_{1}-f_{2}\right\Vert _{s}=0$ are
identified. If 1$\leq s_{1}\leq s_{2}<\infty$, we have $\mathcal{D}%
^{1}\supseteq\mathcal{D}^{s_{1}}\supseteq\mathcal{D}^{s_{2}}$ as sets (but
they are endowed with different norms). For all $f\in\mathcal{D}^{1}$, the
mean value $\ $%
\[
M(f):=\lim\limits_{N\rightarrow\infty}\frac{1}{N}\sum\limits_{1\leq n\leq
N}f(n)
\]
\ exists. The space \bigskip$\mathcal{D}^{2}$ is a Hilbert space with inner product%

\[
\left\langle f,h\right\rangle :=M(f\overline{h}),\ \ \ \ \ f,h\in
\mathcal{D}^{2}\mathcal{\ }.
\]
For $u\in$ $\mathbb{R}$, define $e_{u}(n):=e^{2\pi iun}$, $n\in\mathbb{N}$. In
$\mathcal{D}^{2}$, we have canonical orthonormal base $\left\{  e_{a/b}%
\right\}  $, where $1\leq a\leq b$ and $\gcd(a,b)=1.$
\[
\]
For all $f$\bigskip\ $\in\mathcal{D}^{1}$, the Fourier coefficients
$\widehat{f}(u):=M(fe_{-u}),$ $u\in\mathbb{R}$, exist.

\begin{lemma}
For $f$\bigskip\ $\in\mathcal{D}^{1},u\notin\mathbb{Q}$, we have $\widehat
{f}(u)=0.$
\end{lemma}

\begin{proof}
Let $f$\bigskip\ $\in\mathcal{D}^{1}$. For any $\varepsilon>0$ there is a
linear combination $%
%TCIMACRO{\tsum _{1\leq v\leq V}}%
%BeginExpansion
{\textstyle\sum_{1\leq v\leq V}}
%EndExpansion
e_{v}(n)$, such that $\left\Vert f-%
%TCIMACRO{\tsum _{1\leq v\leq V}}%
%BeginExpansion
{\textstyle\sum_{1\leq v\leq V}}
%EndExpansion
e_{v}\right\Vert _{1}<\varepsilon$, where $v\in\mathbb{Q}$. So$\left\vert
f(n)-%
%TCIMACRO{\tsum _{1\leq v\leq V}}%
%BeginExpansion
{\textstyle\sum_{1\leq v\leq V}}
%EndExpansion
e_{v}(n)\right\vert <\varepsilon$, for all $1\leq n\leq N$. Therefore we have

\mathstrut

$\left\vert \widehat{f}(u)-\lim\limits_{N\rightarrow\infty}\dfrac{1}{N}%
%TCIMACRO{\dsum \limits_{1\leq n\leq N}}%
%BeginExpansion
{\displaystyle\sum\limits_{1\leq n\leq N}}
%EndExpansion%
%TCIMACRO{\dsum \limits_{1\leq v\leq V}}%
%BeginExpansion
{\displaystyle\sum\limits_{1\leq v\leq V}}
%EndExpansion
e_{v}(n)e_{-u}(n)\right\vert =$

\mathstrut

$\ \ \ \ \ \ \ \ =\lim\limits_{N\rightarrow\infty}\dfrac{1}{N}\left\vert
%TCIMACRO{\dsum \limits_{1\leq n\leq N}}%
%BeginExpansion
{\displaystyle\sum\limits_{1\leq n\leq N}}
%EndExpansion
\left(  f(n)-%
%TCIMACRO{\dsum \limits_{1\leq v\leq V}}%
%BeginExpansion
{\displaystyle\sum\limits_{1\leq v\leq V}}
%EndExpansion
e_{v}(n)\right)  e_{-u}(n)\right\vert \leq$

\mathstrut

$\ \ \ \ \ \ \ \ \ \ \ \ \ \ \ \ \ \ \ \ \ \ \ \ \ \ \ \ \ \ \ \ \ \ \ \ \ \ \ \ \ \ \ \ \ \ \leq
\lim\limits_{N\rightarrow\infty}\dfrac{1}{N}%
%TCIMACRO{\dsum \limits_{1\leq n\leq N}}%
%BeginExpansion
{\displaystyle\sum\limits_{1\leq n\leq N}}
%EndExpansion
\varepsilon\left\vert e_{-u}(n)\right\vert \leq\varepsilon.$

\mathstrut

Let $u\notin\mathbb{Q}$, since $v-u\notin\mathbb{Q}$ we have%
\[
\left\vert \lim\limits_{N\rightarrow\infty}\dfrac{1}{N}%
%TCIMACRO{\dsum \limits_{1\leq n\leq N}}%
%BeginExpansion
{\displaystyle\sum\limits_{1\leq n\leq N}}
%EndExpansion%
%TCIMACRO{\dsum \limits_{1\leq v\leq V}}%
%BeginExpansion
{\displaystyle\sum\limits_{1\leq v\leq V}}
%EndExpansion
e_{v}(n)e_{-u}(n)\right\vert =\left\vert \lim\limits_{N\rightarrow\infty
}\dfrac{1}{N}%
%TCIMACRO{\dsum \limits_{1\leq n\leq N}}%
%BeginExpansion
{\displaystyle\sum\limits_{1\leq n\leq N}}
%EndExpansion%
%TCIMACRO{\dsum \limits_{1\leq v\leq V}}%
%BeginExpansion
{\displaystyle\sum\limits_{1\leq v\leq V}}
%EndExpansion
e^{2\pi in(v-u)}\right\vert =
\]%
\[
=\left\vert \lim\limits_{N\rightarrow\infty}\dfrac{1}{N}%
%TCIMACRO{\dsum \limits_{1\leq v\leq V}}%
%BeginExpansion
{\displaystyle\sum\limits_{1\leq v\leq V}}
%EndExpansion
\dfrac{1-e^{2\pi iN(v-u)}}{1-e^{2\pi i(v-u)}}\right\vert \leq\lim
\limits_{N\rightarrow\infty}\dfrac{1}{N}%
%TCIMACRO{\dsum \limits_{1\leq v\leq V}}%
%BeginExpansion
{\displaystyle\sum\limits_{1\leq v\leq V}}
%EndExpansion
\dfrac{2}{const}=0,
\]
and the lemma follows.
\end{proof}

\section{\bigskip\label{Limit-Period-sect}Limit periodicity of weighted
multiplicities}

In this section we will prove that the weighted multiplicities function
$\beta_{Q}(n)$ is limit periodic (prop. \ref{prop(D1)}), and as consequence,
we will obtain the formula for calculating its mean square (end of the section).

\begin{proposition}
\label{prop(D1)}The functions $\beta_{Q}(n)\in\mathcal{D}^{1}$, and for $1\leq
s\leq2$,
\[
\lim\limits_{P\rightarrow\infty}\left\Vert \beta_{Q}-\beta_{P,Q}\right\Vert
_{s}=0\text{.}%
\]

\end{proposition}%

\[
\]
Write $\beta_{Q}(n)-\beta_{P,Q}(n)=\bigtriangleup_{P}^{(1)}(n)+\bigtriangleup
_{P}^{(2)}(n),$ where%

\[
\bigtriangleup_{P}^{(1)}(n):=\sum\limits_{\substack{D,v\geq1;Dv^{2}%
=n^{2}-4\\p\mid v\text{ for some }p>P}}\frac{1}{v}\prod\limits_{q\mid
Q}\left\{
\begin{array}
[c]{cc}%
2, & q^{2}\mid D\\
1+\left(  \dfrac{D}{q}\right)  & q^{2}\nmid D
\end{array}
\right\}  L(1,\chi_{D})
\]
and%

\begin{align*}
\bigtriangleup_{P}^{(2)}(n)  &  :=\sum\limits_{\substack{D,v\geq1;Dv^{2}%
=n^{2}-4\\p\mid v\Rightarrow p\leq P}}\frac{1}{v}\prod\limits_{q\mid
Q}\left\{
\begin{array}
[c]{cc}%
2, & q^{2}\mid D\\
1+\left(  \dfrac{D}{q}\right)  & q^{2}\nmid D
\end{array}
\right\}  \times\\
&  \text{ \ \ \ \ \ \ \ \ \ \ \ \ \ \ \ \ \ \ \ \ \ \ \ \ \ \ \ }\times\left(
L(1,\chi_{D})-\prod\limits_{p\leq P}\left(  1-\frac{\chi_{D}(p)}{p}\right)
^{-1}\right)  .
\end{align*}

\begin{lemma}
\label{delta_p_1}For $P\geq Q$ we have
\[
\frac{1}{x}\sum\limits_{2<n\leq x}\left\vert \bigtriangleup_{P}^{(1)}%
(n)\right\vert ^{2}\ll\sum\limits_{v>P}\frac{1}{v^{2}}.
\]

\end{lemma}

\begin{proof}
Note that%

\begin{equation}
\bigtriangleup_{P}^{(1)}(n)\leq\sum\limits_{\substack{D,v\geq1;Dv^{2}=n^{2}-4
\\p\mid v\text{ for some }p>P}}\frac{2^{\omega(Q)}}{v}L(1,\chi_{D}),
\label{newbeta}%
\end{equation}
where $\omega(Q)$ is the number of prime divisors of $Q.$ Cauchy's inequality
gives%
\[
\left\vert \bigtriangleup_{P}^{(1)}(n)\right\vert \leq\left(  \sum
\limits_{\substack{D,v\geq1;Dv^{2}=n^{2}-4 \\p\mid v\text{ for some }%
p>P}}\frac{2^{2\omega(Q)}}{v^{2}}\right)  ^{1/2}\left(  \sum
\limits_{\substack{D,v\geq1;Dv^{2}=n^{2}-4 \\p\mid v\text{ for some }%
p>P}}L(1,\chi_{D})^{2}\right)  ^{1/2}.
\]
For $x\geq1$, this gives%
\[
\sum\limits_{2<n\leq x}\left\vert \bigtriangleup_{P}^{(1)}(n)\right\vert
^{2}\leq\sum\limits_{v>P}\frac{2^{2\omega(Q)}}{v^{2}}\sum
\limits_{\substack{2<n\leq x \\D,v\geq1;Dv^{2}=n^{2}-4}}L(1,\chi_{D})^{2}.
\]
M.Peter shows \cite{man},\cite{man1} that the last sum is%
\[
\sum\limits_{\substack{2<n\leq x \\D,v\geq1;Dv^{2}=n^{2}-4}}L(1,\chi_{D}%
)^{2}\sim const\cdot x\text{,}%
\]
as $x\rightarrow\infty$. Therefore we have the claim of the lemma.
\end{proof}

\bigskip

In order to estimate $\bigtriangleup_{P}^{(2)}(n)$ we must compare
$L(1,\chi_{D})$ with a partial product of its Euler products. This is done by
comparing both terms with a smoothed version of the Dirichlet series for
$L(1,\chi_{D})$. Let $N\geq1$. Then%

\[
\bigtriangleup_{P}^{(2)}(n)=\bigtriangleup_{P,N}^{(2,1)}(n)+\bigtriangleup
_{P,N}^{(2,2)}(n)+\bigtriangleup_{P,N}^{(2,3)}(n),
\]
where%

\[
\bigtriangleup_{P,N}^{(2,1)}(n):=\sum\limits_{\substack{D,v\geq1;Dv^{2}%
=n^{2}-4 \\p\mid v\Rightarrow p\leq P}}\frac{1}{v}\prod\limits_{q\mid
Q}\left\{
\begin{array}
[c]{cc}%
2, & q^{2}\mid D\\
1+\left(  \dfrac{D}{q}\right)  & q^{2}\nmid D
\end{array}
\right\}  \left(  L(1,\chi_{D})-\sum\limits_{l\geq1}\frac{\chi_{D}(l)}%
{l}e^{-l/N}\right)  ,
\]

\[
\bigtriangleup_{P,N}^{(2,2)}(n):=\sum\limits_{\substack{D,v\geq1;Dv^{2}%
=n^{2}-4 \\p\mid v\Rightarrow p\leq P}}\frac{1}{v}\prod\limits_{q\mid
Q}\left\{
\begin{array}
[c]{cc}%
2, & q^{2}\mid D\\
1+\left(  \dfrac{D}{q}\right)  & q^{2}\nmid D
\end{array}
\right\}  \sum\limits_{l\geq1:p\mid l\text{ \ for some }p>P}\frac{\chi_{D}%
(l)}{l}e^{-l/N},
\]

\[
\bigtriangleup_{P,N}^{(2,3)}(n):=\sum\limits_{\substack{D,v\geq1;Dv^{2}%
=n^{2}-4\\p\mid v\Rightarrow p\leq P}}\frac{1}{v}\prod\limits_{q\mid
Q}\left\{
\begin{array}
[c]{cc}%
2, & q^{2}\mid D\\
1+\left(  \dfrac{D}{q}\right)  & q^{2}\nmid D
\end{array}
\right\}  \sum\limits_{l\geq1:p\mid l\Rightarrow p\leq P}\frac{\chi_{D}(l)}%
{l}\left(  e^{-l/N}-1\right)  ;
\]
For the following approximations we will use two lemmas:

\begin{lemma}
[Sarnak]\label{sarnak}%
\[
\sum\limits_{\substack{2<n\leq x \\d,v\geq1 \\dv^{2}=n^{2}-4}}1\sim const\cdot
x
\]

\end{lemma}

\begin{proof}
\cite[Lemma 4.2]{sar}
\end{proof}

\begin{lemma}
[ Peter]\label{peter}For $l,v\in\mathbb{N}$ and $x\geq3$, we have%
\[
\sum\limits_{\substack{2<n\leq x \\d\geq1 \\dv^{2}=n^{2}-4}}\chi_{d}%
(l)\ll\frac{x}{v^{2-\varepsilon}K(l)}+v^{\varepsilon}l,
\]

where $K(l)$\ is the squarefree kernel of $l$\ and $\varepsilon>0$ is arbitrary.
\end{lemma}

\begin{proof}
See \cite{man2}, estimate (2.7).
\end{proof}

\begin{lemma}
For $P\geq Q$ and $x,N\geq1$, we have%
\[
\frac{1}{x}\sum\limits_{2<n\leq x}\left\vert \triangle_{P,N}^{(2,3)}%
(n)\right\vert ^{2}\ll\left(  N^{-1/2}+\sum\limits_{l>\sqrt{N}:p\mid
l\Rightarrow p\leq P}\frac{1}{l}\right)  ^{2}.
\]

\end{lemma}

\begin{proof}
Since $\left\vert e^{-u}-1\right\vert \ll u$ for $0\leq u\leq1$, we see that
for $n>2$ the inner sum in $\triangle_{P,N}^{(2,3)}(n)$ is%
\[
\ll%
%TCIMACRO{\dsum \limits_{l\geq1:p\mid l\Rightarrow p\leq P}}%
%BeginExpansion
{\displaystyle\sum\limits_{l\geq1:p\mid l\Rightarrow p\leq P}}
%EndExpansion
\dfrac{1}{l}\left\vert e^{-l/N}-1\right\vert \ll%
%TCIMACRO{\dsum \limits_{l>\sqrt{N}:p\mid l\Rightarrow p\leq P}}%
%BeginExpansion
{\displaystyle\sum\limits_{l>\sqrt{N}:p\mid l\Rightarrow p\leq P}}
%EndExpansion
\dfrac{2}{l}+%
%TCIMACRO{\dsum \limits_{1<l\leq\sqrt{N}:p\mid l\Rightarrow p\leq P}}%
%BeginExpansion
{\displaystyle\sum\limits_{1<l\leq\sqrt{N}:p\mid l\Rightarrow p\leq P}}
%EndExpansion
\dfrac{1}{l}\dfrac{l}{N}%
\]%
\[
\ll%
%TCIMACRO{\dsum \limits_{l>\sqrt{N}:p\mid l\Rightarrow p\leq P}}%
%BeginExpansion
{\displaystyle\sum\limits_{l>\sqrt{N}:p\mid l\Rightarrow p\leq P}}
%EndExpansion
\dfrac{1}{l}+N^{-1/2}=:c_{1}(P,N).
\]
Cauchy's inequality and (\ref{newbeta}) give%
\[%
%TCIMACRO{\dsum \limits_{2<n\leq x}}%
%BeginExpansion
{\displaystyle\sum\limits_{2<n\leq x}}
%EndExpansion
\left\vert \triangle_{P,N}^{(2,3)}(n)\right\vert ^{2}\ll%
%TCIMACRO{\dsum \limits_{2<n\leq x}}%
%BeginExpansion
{\displaystyle\sum\limits_{2<n\leq x}}
%EndExpansion
\left(
%TCIMACRO{\dsum \limits_{D,v\geq1;Dv^{2}=n^{2}-4}}%
%BeginExpansion
{\displaystyle\sum\limits_{D,v\geq1;Dv^{2}=n^{2}-4}}
%EndExpansion
\dfrac{2^{\omega(Q)}}{v}\right)  ^{2}c_{1}(P,N)^{2}\ll
\]%
\[
\ll c_{1}(P,N)^{2}%
%TCIMACRO{\dsum \limits_{2<n\leq x}}%
%BeginExpansion
{\displaystyle\sum\limits_{2<n\leq x}}
%EndExpansion
\left(
%TCIMACRO{\dsum \limits_{D,v\geq1;Dv^{2}=n^{2}-4}}%
%BeginExpansion
{\displaystyle\sum\limits_{D,v\geq1;Dv^{2}=n^{2}-4}}
%EndExpansion
\dfrac{2^{2\omega(Q)}}{v^{2}}\right)  \left(
%TCIMACRO{\dsum \limits_{D,v\geq1;Dv^{2}=n^{2}-4}}%
%BeginExpansion
{\displaystyle\sum\limits_{D,v\geq1;Dv^{2}=n^{2}-4}}
%EndExpansion
1\right)  \ll
\]%
\[
\ll c_{1}(P,N)^{2}%
%TCIMACRO{\dsum \limits_{\substack{2<n\leq x \\D,v\geq1;Dv^{2}=n^{2}-4}}}%
%BeginExpansion
{\displaystyle\sum\limits_{\substack{2<n\leq x \\D,v\geq1;Dv^{2}=n^{2}-4}}}
%EndExpansion
1.
\]
By using \ref{sarnak} the result follows.
\end{proof}

\begin{lemma}
For $P\geq Q$ and $x,N\geq1$, we have%
\[
\frac{1}{x}\sum\limits_{2<n\leq x}\left\vert \triangle_{P,N}^{(2,2)}%
(n)\right\vert ^{2}\ll\sum\limits_{l>P^{2}}\frac{\tau(l)}{lK(l)}+\frac
{1}{x^{1/3-\varepsilon}}N^{2}%
\]

where $\tau(l)$ is the number of positive divisors of $l.$
\end{lemma}

\begin{proof}
Let $\alpha>\dfrac{1}{2}$. We write%
\[
\triangle_{P,N}^{(2,2)}(n)=%
%TCIMACRO{\dsum \limits_{\substack{D,v\geq1 \\Dv^{2}=n^{2}-4 \\v\leq n^{\alpha
%}}}}%
%BeginExpansion
{\displaystyle\sum\limits_{\substack{D,v\geq1 \\Dv^{2}=n^{2}-4 \\v\leq
n^{\alpha}}}}
%EndExpansion
\dfrac{1}{v}%
%TCIMACRO{\dprod \limits_{q\mid Q}}%
%BeginExpansion
{\displaystyle\prod\limits_{q\mid Q}}
%EndExpansion
\left\{
\begin{array}
[c]{cc}%
2, & q^{2}\mid D\\
1+\left(  \dfrac{D}{q}\right)  & q^{2}\nmid D
\end{array}
\right\}
%TCIMACRO{\dsum \limits_{\substack{l\geq1:p\mid l\text{ for} \\\text{some }%
%p>P}}}%
%BeginExpansion
{\displaystyle\sum\limits_{\substack{l\geq1:p\mid l\text{ for} \\\text{some
}p>P}}}
%EndExpansion
\dfrac{\chi_{D}(l)}{l}e^{-l/N}+
\]%
\[
+%
%TCIMACRO{\dsum \limits_{\substack{D,v\geq1 \\Dv^{2}=n^{2}-4 \\v>n^{\alpha}}}}%
%BeginExpansion
{\displaystyle\sum\limits_{\substack{D,v\geq1 \\Dv^{2}=n^{2}-4 \\v>n^{\alpha}%
}}}
%EndExpansion
\dfrac{1}{v}%
%TCIMACRO{\dprod \limits_{q\mid Q}}%
%BeginExpansion
{\displaystyle\prod\limits_{q\mid Q}}
%EndExpansion
\left\{
\begin{array}
[c]{cc}%
2, & q^{2}\mid D\\
1+\left(  \dfrac{D}{q}\right)  & q^{2}\nmid D
\end{array}
\right\}
%TCIMACRO{\dsum \limits_{\substack{l\geq1:p\mid l\text{ for} \\\text{some }%
%p>P}}}%
%BeginExpansion
{\displaystyle\sum\limits_{\substack{l\geq1:p\mid l\text{ for} \\\text{some
}p>P}}}
%EndExpansion
\dfrac{\chi_{D}(l)}{l}e^{-l/N}=
\]%
\[
=\triangle_{P,N}^{(2,2,1)}(n)+\triangle_{P,N}^{(2,2,2)}(n).
\]
A trivial estimate gives%
\[
\triangle_{P,N}^{(2,2,2)}(n)\leq%
%TCIMACRO{\dsum \limits_{\substack{D,v\geq1 \\Dv^{2}=n^{2}-4 \\v>n^{\alpha}}}}%
%BeginExpansion
{\displaystyle\sum\limits_{\substack{D,v\geq1 \\Dv^{2}=n^{2}-4 \\v>n^{\alpha}%
}}}
%EndExpansion
\dfrac{2^{\omega(Q)}}{v}%
%TCIMACRO{\dsum \limits_{\substack{l\geq1:p\mid l\text{ for} \\\text{some }%
%p>P}}}%
%BeginExpansion
{\displaystyle\sum\limits_{\substack{l\geq1:p\mid l\text{ for} \\\text{some
}p>P}}}
%EndExpansion
\dfrac{\chi_{D}(l)}{l}e^{-l/N}\ll%
%TCIMACRO{\dsum \limits_{\substack{D,v\geq1 \\Dv^{2}=n^{2}-4 \\v>n^{\alpha}}}}%
%BeginExpansion
{\displaystyle\sum\limits_{\substack{D,v\geq1 \\Dv^{2}=n^{2}-4 \\v>n^{\alpha}%
}}}
%EndExpansion
\dfrac{1}{v}%
%TCIMACRO{\dsum \limits_{\substack{l\geq1:p\mid l\text{ for} \\\text{some }%
%p>P}}}%
%BeginExpansion
{\displaystyle\sum\limits_{\substack{l\geq1:p\mid l\text{ for} \\\text{some
}p>P}}}
%EndExpansion
\dfrac{1}{l}e^{-l/N}\ll
\]%
\[
\log N%
%TCIMACRO{\dsum \limits_{\substack{D,v\geq1 \\Dv^{2}=n^{2}-4 \\v>n^{\alpha}
%}}}%
%BeginExpansion
{\displaystyle\sum\limits_{\substack{D,v\geq1 \\Dv^{2}=n^{2}-4 \\v>n^{\alpha}
}}}
%EndExpansion
\dfrac{1}{v}\ll\log N\cdot\dfrac{1}{n^{\alpha}}\tau(n^{2}-4)\ll\log
N\cdot\dfrac{1}{n^{\alpha-\varepsilon}}.
\]
Thus , since $\alpha>\dfrac{1}{2}$, we have%
\begin{equation}
\sum\limits_{2<n\leq x}\left\vert \triangle_{P,N}^{(2,2,2)}(n)\right\vert
^{2}\ll\left(  \log N\right)  ^{2}\cdot\sum\limits_{2<n\leq x}\frac
{1}{n^{2(\alpha-\varepsilon)}}\ll\left(  \log N\right)  ^{2}. \label{delta222}%
\end{equation}
\ By Cauchy's inequality%
\[
\left\vert \triangle_{P,N}^{(2,2,1)}(n)\right\vert \leq%
%TCIMACRO{\dsum \limits_{\substack{D,v\geq1 \\Dv^{2}=n^{2}-4 \\v\leq n^{\alpha
%}}}}%
%BeginExpansion
{\displaystyle\sum\limits_{\substack{D,v\geq1 \\Dv^{2}=n^{2}-4 \\v\leq
n^{\alpha}}}}
%EndExpansion
\dfrac{2^{\omega(Q)}}{v}\left\vert
%TCIMACRO{\dsum \limits_{l\geq1:p\mid l\text{ for some }p>P}}%
%BeginExpansion
{\displaystyle\sum\limits_{l\geq1:p\mid l\text{ for some }p>P}}
%EndExpansion
\dfrac{\chi_{D}(l)}{l}e^{-l/N}\right\vert \leq
\]

\[
\left(
%TCIMACRO{\dsum \limits_{D,v\geq1;Dv^{2}=n^{2}-4}}%
%BeginExpansion
{\displaystyle\sum\limits_{D,v\geq1;Dv^{2}=n^{2}-4}}
%EndExpansion
\dfrac{2^{2\omega(Q)}}{v^{2}}\right)  ^{\frac{1}{2}}\cdot\left(
%TCIMACRO{\dsum \limits_{D,v\geq1;Dv^{2}=n^{2}-4;v\leq n^{\alpha}}}%
%BeginExpansion
{\displaystyle\sum\limits_{D,v\geq1;Dv^{2}=n^{2}-4;v\leq n^{\alpha}}}
%EndExpansion
\left(
%TCIMACRO{\dsum \limits_{l\geq1:p\mid l\text{ for some }p>P}}%
%BeginExpansion
{\displaystyle\sum\limits_{l\geq1:p\mid l\text{ for some }p>P}}
%EndExpansion
\dfrac{\chi_{D}(l)}{l}e^{-l/N}\right)  ^{2}\right)  ^{\frac{1}{2}}.
\]
Thus for $x\geq1,$%
\[%
%TCIMACRO{\dsum \limits_{2<n\leq x}}%
%BeginExpansion
{\displaystyle\sum\limits_{2<n\leq x}}
%EndExpansion
\left\vert \triangle_{P,N}^{(2,2,1)}(n)\right\vert ^{2}\ll%
%TCIMACRO{\dsum \limits_{2<n\leq x}}%
%BeginExpansion
{\displaystyle\sum\limits_{2<n\leq x}}
%EndExpansion%
%TCIMACRO{\dsum \limits_{D,v\geq1;Dv^{2}=n^{2}-4;v\leq x^{\alpha}}}%
%BeginExpansion
{\displaystyle\sum\limits_{D,v\geq1;Dv^{2}=n^{2}-4;v\leq x^{\alpha}}}
%EndExpansion
\left(
%TCIMACRO{\dsum \limits_{l\geq1:p\mid l\text{ for some }p>P}}%
%BeginExpansion
{\displaystyle\sum\limits_{l\geq1:p\mid l\text{ for some }p>P}}
%EndExpansion
\dfrac{\chi_{D}(l)}{l}e^{-l/N}\right)  ^{2}=
\]

\mathstrut%
\[
=%
%TCIMACRO{\dsum \limits_{\substack{l_{1},l_{2}:p_{i}\mid l_{i} \\\text{for some
%}p_{i}>P}}}%
%BeginExpansion
{\displaystyle\sum\limits_{\substack{l_{1},l_{2}:p_{i}\mid l_{i} \\\text{for
some }p_{i}>P}}}
%EndExpansion
\dfrac{1}{l_{1}l_{2}}e^{-(l_{1}+l_{2})/N}%
%TCIMACRO{\dsum \limits_{1\leq v\leq x^{\alpha}}}%
%BeginExpansion
{\displaystyle\sum\limits_{1\leq v\leq x^{\alpha}}}
%EndExpansion%
%TCIMACRO{\dsum \limits_{\substack{2<n\leq x \\D\geq1 \\Dv^{2}=n^{2}-4 }}}%
%BeginExpansion
{\displaystyle\sum\limits_{\substack{2<n\leq x \\D\geq1 \\Dv^{2}=n^{2}-4 }}}
%EndExpansion
\chi_{D}(l_{1}l_{2}).
\]
Applying Peter's lemma \ref{peter} to the innermost sum gives the estimate%
\[%
%TCIMACRO{\dsum \limits_{2<n\leq x}}%
%BeginExpansion
{\displaystyle\sum\limits_{2<n\leq x}}
%EndExpansion
\left\vert \triangle_{P,N}^{(2,2,1)}(n)\right\vert ^{2}\ll%
%TCIMACRO{\dsum \limits_{l>P^{2}}}%
%BeginExpansion
{\displaystyle\sum\limits_{l>P^{2}}}
%EndExpansion
\dfrac{1}{l}\tau(l)%
%TCIMACRO{\dsum \limits_{1\leq v\leq x^{\alpha}}}%
%BeginExpansion
{\displaystyle\sum\limits_{1\leq v\leq x^{\alpha}}}
%EndExpansion
\dfrac{x}{v^{2-\varepsilon}K(l)}+%
%TCIMACRO{\dsum \limits_{l_{1},l_{2}\geq1}}%
%BeginExpansion
{\displaystyle\sum\limits_{l_{1},l_{2}\geq1}}
%EndExpansion
\dfrac{l_{1}l_{2}}{l_{1}l_{2}}e^{-(l_{1}+l_{2})/N}%
%TCIMACRO{\dsum \limits_{1\leq v\leq x^{\alpha}}}%
%BeginExpansion
{\displaystyle\sum\limits_{1\leq v\leq x^{\alpha}}}
%EndExpansion
v^{\varepsilon}%
\]

\mathstrut%
\[
\ll x%
%TCIMACRO{\dsum \limits_{l>P^{2}}}%
%BeginExpansion
{\displaystyle\sum\limits_{l>P^{2}}}
%EndExpansion
\dfrac{\tau(l)}{lK(l)}+N^{2}x^{\alpha(1+\varepsilon)}.
\]
Thus together with (\ref{delta222}), for $\alpha=2/3>1/2$ we have%
\[
\dfrac{1}{x}%
%TCIMACRO{\dsum \limits_{2<n\leq x}}%
%BeginExpansion
{\displaystyle\sum\limits_{2<n\leq x}}
%EndExpansion
\left\vert \triangle_{P,N}^{(2,2)}(n)\right\vert ^{2}\ll%
%TCIMACRO{\dsum \limits_{l>P^{2}}}%
%BeginExpansion
{\displaystyle\sum\limits_{l>P^{2}}}
%EndExpansion
\dfrac{\tau(l)}{lK(l)}+\dfrac{1}{x}(\log N)^{2}+\dfrac{1}{x^{1/3-\varepsilon}%
}N^{2}.
\]

\end{proof}%

\[
\]
In order to estimate $\triangle_{P,N}^{(2,1)}(n)$ we must show that the error%

\[
I(D,N):=L(1,\chi_{D})-\sum\limits_{l\geq1}\frac{\chi_{D}(l)}{l}e^{-l/N},
\]
which comes from smoothing Dirichlet series expansion of $L(1,\chi_{D})$, is
small for large $N$.

\begin{lemma}
For $1/2<\sigma_{0}<1$ define the rectangle

\mathstrut

$R_{x}:=\{s\in\mathbb{C}\mid\sigma_{0}\leq\operatorname{Re}(s)\leq1,\left|
\operatorname{Im}(s)\right|  \leq\log^{2}x\}$.

\mathstrut

(a) If $L(s,\chi_{D})$ has no zeros in $R_{x}$ and $D\leq x^{2}$, then

\mathstrut

for $\operatorname{Re}(s)=\kappa,\left|  \operatorname{Im}(s)\right|
\leq(\log x)^{2}/2$ \ \ \ \ \ \ \ \ \ \ \ \ \ \ \ \ \ \ \ \ \ \ $I(D,N)\ll
x^{\varepsilon}N^{(\kappa-1)}$;

\mathstrut

(b) If $L(s,\chi_{D})$ has zeros in $R_{x}$, then

\mathstrut

$\#\{(n,v,D)\mid2<n\leq x,D,v\geq1,n^{2}-Dv^{2}=4$

$\ \ \ \ \ \ \ \ \ \ \ \ \ \ \ \ \ \ \ \ \ \ \ \ \ \ \ \ \ \ \ \ \ \ \ \ \ \ \ \ \ \ \ \ \ \ \ \ \ ,L(s,\chi
_{D})$ has zeros in $R_{x}\}\ll x^{\mu+\varepsilon}$,

\mathstrut

where $\mu:=8(1-\sigma_{0})/\sigma_{0}<1,$ $\sigma_{0}<\kappa<1$
\end{lemma}

\begin{proof}
See \cite[Lemma 3.6]{man}.
\end{proof}

\begin{lemma}
There are $0<\kappa,\mu<1$ such that for $P\geq Q,x,N\geq1$ and $\varepsilon
>0$ we have%
\[
\ \frac{1}{x}\sum\limits_{2<n\leq x}\left\vert \triangle_{P,N}^{(2,1)}%
(n)\right\vert ^{2}\ll x^{\varepsilon}N^{2(\kappa-1)}+x^{\mu-1+\varepsilon
}(\log(x^{2}N))^{2}.
\]

\end{lemma}

\begin{proof}
Note that a trivial estimation gives $I(D,N)\ll\log(DN)$. Cauchy's inequality,
previous lemma and (\ref{newbeta}) give%
\[%
%TCIMACRO{\dsum \limits_{2<n\leq x}}%
%BeginExpansion
{\displaystyle\sum\limits_{2<n\leq x}}
%EndExpansion
\left\vert \triangle_{P,N}^{(2,1)}(n)\right\vert ^{2}\ll%
%TCIMACRO{\dsum \limits_{2<n\leq x}}%
%BeginExpansion
{\displaystyle\sum\limits_{2<n\leq x}}
%EndExpansion
\left(
%TCIMACRO{\dsum \limits_{D,v\geq1;Dv^{2}=n^{2}-4}}%
%BeginExpansion
{\displaystyle\sum\limits_{D,v\geq1;Dv^{2}=n^{2}-4}}
%EndExpansion
\dfrac{2^{2\omega(Q)}}{v^{2}}\right)  \left(
%TCIMACRO{\dsum \limits_{D,v\geq1;Dv^{2}=n^{2}-4}}%
%BeginExpansion
{\displaystyle\sum\limits_{D,v\geq1;Dv^{2}=n^{2}-4}}
%EndExpansion
\left\vert I(D,N)\right\vert ^{2}\right)  \ll
\]

\mathstrut%
\[
\ll%
%TCIMACRO{\dsum \limits_{\substack{2<n\leq x;D,v\geq1:Dv^{2}=n^{2}-4
%\\L(s,\chi_{D})\text{ has no zeros in }R_{x}}}}%
%BeginExpansion
{\displaystyle\sum\limits_{\substack{2<n\leq x;D,v\geq1:Dv^{2}=n^{2}-4
\\L(s,\chi_{D})\text{ has no zeros in }R_{x}}}}
%EndExpansion
\left(  x^{\varepsilon}N^{(\kappa-1)}\right)  ^{2}+%
%TCIMACRO{\dsum \limits_{\substack{2<n\leq x;D,v\geq1:Dv^{2}=n^{2}-4
%\\L(s,\chi_{D})\text{ has a zero in }R_{x}}}}%
%BeginExpansion
{\displaystyle\sum\limits_{\substack{2<n\leq x;D,v\geq1:Dv^{2}=n^{2}-4
\\L(s,\chi_{D})\text{ has a zero in }R_{x}}}}
%EndExpansion
\log^{2}(DN)\ll
\]

\mathstrut%
\[
\ll x\left(  x^{\varepsilon}N^{(\kappa-1)}\right)  ^{2}+x^{\mu+\varepsilon
}(\log(x^{2}N))^{2},
\]
which proves the lemma.
\end{proof}

\bigskip

Now the results are collected.

\begin{lemma}
For $P\geq Q$, we have$\ \ \ \ \ \ \ \ \ \ \ \ \ \ \ \ $%
\[
\ \left\Vert \beta_{Q}-\beta_{P,Q}\right\Vert _{2}\ll\left(  \sum
\limits_{v>P}\frac{1}{v^{2}}\right)  ^{1/2}+\left(  \sum\limits_{l>P^{2}}%
\frac{\tau(l)}{lK(l)}\right)  ^{1/2}.
\]

\end{lemma}

\begin{proof}
For $\ x\geq1$ choose $N:=x^{1/8}$. Then previous lemmas show that%
\[
\frac{1}{x}\sum\limits_{2<n\leq x}\left\vert \bigtriangleup_{P}^{(2)}%
(n)\right\vert ^{2}\ll\left(  x^{-1/16}+\sum\limits_{l>x^{1/16}:p\mid
l\Rightarrow p\leq P}\frac{1}{l}\right)  ^{2}+\sum\limits_{l>P^{2}}\frac
{\tau(l)}{lK(l)}+
\]%
\[
+\frac{1}{x}(\log x)^{2}+\frac{1}{x^{1/12-\varepsilon}}+x^{(\kappa
-1)/4+\varepsilon}+x^{\mu-1+\varepsilon}(\log x)^{2}.
\]
Since the series%
\[
\sum\limits_{l\geq1:p\mid l\Rightarrow p\leq P}\frac{1}{l}%
\]
converges, we have for $P\geq Q$
fixed$\ \ \ \ \ \ \ \ \ \ \ \ \ \ \ \ \ \ \ \ \ \ \ \ \ \ \ \ \ \ \ \ \ \ \ $%
\[
\left\Vert \bigtriangleup_{P}^{(2)}(n)\right\Vert _{2}^{2}\ll\sum
\limits_{l>P^{2}}\frac{\tau(l)}{lK(l)}.
\]
Together with Lemma \ref{delta_p_1} this proves the claim.
\end{proof}

\[
\]
Now we are able to prove the proposition \ref{prop(D1)}

\begin{proof}
[ Proof of the proposition \ref{prop(D1)}]By previous lemma we
have$\ \ \ \ \ \ \ \ \ \ \ \ \ \ \ \ \ \ \ \ \ \ $%
\[
\ \ \left\Vert \beta_{Q}-\beta_{P,Q}\right\Vert _{2}\ll\left(  \sum
\limits_{v>P}\frac{1}{v^{2}}\right)  ^{1/2}+\left(  \sum\limits_{l>P^{2}}%
\frac{\tau(l)}{lK(l)}\right)  ^{1/2};
\]
here$\ \ \ \ \ \ \ \ \ \ \ \ \ \ \ \ \ \ \ \ \ \ \ \ \ \ \ \ \ \ \ \ \ \ \ \ \ \ \ \ $%
\[
\ \sum\limits_{v>P}\frac{1}{v^{2}}\longrightarrow0\text{,}%
\]
as $P\rightarrow\infty,$ since the series $\
%TCIMACRO{\tsum _{v\geq1}}%
%BeginExpansion
{\textstyle\sum_{v\geq1}}
%EndExpansion
\frac{1}{v^{2}}$ converges. Furthermore,%
\[
\ \ \sum\limits_{l>P^{2}}\frac{\tau(l)}{lK(l)}\longrightarrow0\text{,}%
\]
as $P\rightarrow\infty,$ since\ $\ \ \ \ \ \ \ \ \ \ \ \ \ \ \ $%
\[
\ \sum\limits_{l\geq1}\frac{\tau(l)}{lK(l)}=\sum\limits_{a,b\geq1:a\text{
squarefree}}\frac{\tau(ab^{2})}{ab^{2}\cdot a}\ll\sum\limits_{a\geq1}%
\frac{a^{\varepsilon}}{a^{2}}\sum\limits_{b\geq1}\frac{b^{2\varepsilon}}%
{b^{2}}<\infty.
\]
Thus $\lim\limits_{P\rightarrow\infty}\left\Vert \beta_{Q}-\beta
_{P,Q}\right\Vert _{2}=0$. For $f:\mathbb{N\longrightarrow C}$ arbitrary and
$1\leq s\leq2$ we have $\left\Vert f\right\Vert _{s}\leq\left\Vert
f\right\Vert _{2}$ by H\"{o}lder's inequality. Thus $\lim\limits_{P\rightarrow
\infty}\left\Vert \beta_{Q}-\beta_{P,Q}\right\Vert _{s}=0$, for all $1\leq
s\leq2$ and, in particular $\lim\limits_{P\rightarrow\infty}\left\Vert
\beta_{Q}-\beta_{P,Q}\right\Vert _{1}=0$.

Since the $b$-th summand of $\beta_{(p,Q)}$ is $p^{2b+1}$-periodic for $p\nmid
Q,$ $2^{2b+3}$-periodic in case $p=2$, and $p^{2b+2}$-periodic in case $p\mid
Q,$ and the series representing $\beta_{(p,Q)}$ is uniformly convergent, the
function $\beta_{(p,Q)}$ is uniformly limit periodic, i.e. $\beta_{(p,Q)}%
\in\mathcal{D}^{u}$; here $\mathcal{D}^{u}$ is the set of all functions which
can be approximated to an arbitrary accuracy by periodic functions with
respect to the supremum norm. Since $\mathcal{D}^{u}$ is closed under
multiplication it follows from Lemma \ref{mult}, that $\beta_{P,Q}%
\in\mathcal{D}^{u}$ for all $P\geq Q$. This gives $\beta_{Q}\in\mathcal{D}%
^{s}$ for all $1\leq s\leq2$.
\end{proof}

\[
\]
So we have now%

\[
\widehat{\beta_{Q}}(0):=M(\beta_{Q}):=\lim\limits_{N\rightarrow\infty}\frac
{1}{N}\sum\limits_{1\leq n\leq N}\beta_{Q}(n).
\]
One can prove the

\begin{lemma}
\label{fucomult}For $b\in\mathbb{N}$, $a\in\mathbb{Z}$, $\gcd(a,b)=1$, choose
$a_{p}\in\mathbb{Z}$ for all $p\mid b$ such that $\sum\limits_{p\mid b}%
a_{p}p^{-ord_{p}b}\equiv ab^{-1}(\operatorname{mod}1).$ Then
\begin{equation}
\widehat{\beta_{Q}}(\frac{a}{b})=\prod\limits_{p\mid b}\widehat{\beta_{(p,Q)}%
}(\frac{a_{p}}{p^{ord_{p}b}}). \label{man11}%
\end{equation}

\end{lemma}

\begin{proof}
Word by word the proof of the same fact in \cite[Lemma 4.3]{man}
\end{proof}

\begin{corollary}%
\[
\lim\limits_{N\rightarrow\infty}\dfrac{1}{N}%
%TCIMACRO{\dsum \limits_{1\leq n\leq N}}%
%BeginExpansion
{\displaystyle\sum\limits_{1\leq n\leq N}}
%EndExpansion
\beta_{Q}(n)=\widehat{\beta_{Q}}(0)=1
\]

\end{corollary}%

\[
\]
By Parseval's equality and by previous lemma and corollary%

\[
M(\beta\overline{\beta}):=\lim\limits_{N\rightarrow\infty}\frac{1}{N}%
\sum\limits_{2\leq n\leq N}\beta_{Q}^{2}(n)=\sum\limits_{b\geq1}%
\sum\limits_{\substack{1\leq a\leq b \\\gcd(a,b)=1}}\left\vert \widehat
{\beta_{Q}}\left(  \frac{a}{b}\right)  \right\vert ^{2}=
\]

\begin{equation}
=\prod\limits_{p\text{ - prime}}\left(  1+\sum\limits_{c\geq1}\sum
\limits_{\substack{1\leq a\leq p^{c} \\a\neq0(\operatorname{mod}p)}}\left\vert
\widehat{\beta_{(p,Q)}}\left(  \frac{a}{p^{c}}\right)  \right\vert
^{2}\right)  . \label{bet-sq-mean}%
\end{equation}
Here the term $1$ in a brackets is a contribution of $c=0,$ that is
$\left\vert \widehat{\beta_{(p,Q)}}\left(  0\right)  \right\vert ^{2}.$

\section{Calculating the mean square of weighted multiplicities $\beta_{Q}%
(n)$}

\bigskip In this section we will prove Theorem \ref{mean-sq-theorem(congr)}.%
\[
\]
Define the functions
\[
\beta_{(p,Q,b)}(n):=\left(  1-\frac{1}{p}\chi_{(n^{2}-4)p^{-2b}}(p)\right)
^{-1}\cdot\mathbb{I}_{p^{b}}(n),
\]
and calculate the Fourier coefficients of the $\beta_{(p,Q)}(n)$ by \ the
Fourier coefficients of the $\beta_{(p,Q,b)}(n).$%

\begin{equation}
\widehat{\beta_{(p,Q)}}(r)=\sum\limits_{b\geq0}\frac{1}{p^{b}}\widehat
{\beta_{(p,Q,b)}}(r) \label{fuco}%
\end{equation}
In \cite{man} was proved that for all $p\nmid Q,$
\begin{equation}
\widehat{\beta_{(p,Q)}}(0)=1 \label{man2}%
\end{equation}
We will prove that $\widehat{\beta_{(q,Q)}}(0)=1$ holds as well, for all
$q\mid Q$.

\subsection{Calculation of the period \ of $\beta_{(q,Q,b)}(n)$}

Let us calculate now the minimal period of the function defined above:%
\[
\beta_{(q,Q,b)}(n):=\left(  1-\frac{1}{q}\chi_{(n^{2}-4)q^{-2b}}(q)\right)
^{-1}\cdot\mathbb{I}_{q^{b}}(n)
\]

\begin{lemma}
For $q\mid Q$ the minimal period of the $\beta_{(q,Q,b)}(n)$ is $q^{2b+2}$.
\end{lemma}

\begin{proof}
a) we will\ find a period of a $\chi.$ We will look for a minimal $k$ ,such
that $\chi(n+k)=\chi(n)$ for all $n$. That is%

\[
\bigskip\left(  \frac{\left(  n^{2}-4\right)  q^{-2b}}{q}\right)  =\left(
\frac{\left(  \left(  n+k\right)  ^{2}-4\right)  q^{-2b}}{q}\right)  =\left(
\frac{\left(  n^{2}-4\right)  q^{-2b}+\left(  2nk+k^{2}\right)  q^{-2b}}%
{q}\right)
\]
\bigskip and it's true for $k=q^{2b+1}$.

\bigskip b) we will find a period of the $\mathbb{I}_{q^{b}}(n)$. We will look
for a minimal $k$ ,such that $\mathbb{I}_{q^{b}}(n+k)=\mathbb{I}_{q^{b}}(n)$
for all $\ n$.%
\[
n^{2}-4\equiv0(\operatorname{mod}q^{2b})\Longleftrightarrow\left(  n+k\right)
^{2}-4\equiv0(\operatorname{mod}q^{2b})\Longleftrightarrow
\]%
\[
n^{2}-4+2nk+k^{2}\equiv0(\operatorname{mod}q^{2b})\Longleftrightarrow
2nk+k^{2}\equiv0(\operatorname{mod}q^{2b})\text{,}%
\]
and it's true for $k=q^{2b}$.%
\[
\left(  n^{2}-4\right)  q^{-2b}\equiv0(\operatorname{mod}q^{2}%
)\Longleftrightarrow\left(  \left(  n+k\right)  ^{2}-4\right)  q^{-2b}%
\equiv0(\operatorname{mod}q^{2})\Longleftrightarrow
\]%
\[
\left(  n^{2}-4\right)  q^{-2b}+\left(  2nk+k^{2}\right)  q^{-2b}%
\equiv0(\operatorname{mod}q^{2})\Longleftrightarrow\left(  2nk+k^{2}\right)
q^{-2b}\equiv0(\operatorname{mod}q^{2})\text{,}%
\]
and it's true for $k=q^{2b+2}$.

\bigskip So the minimal period of the $\beta_{(q,Q,b)}(n)$ is $q^{2b+2}$ that
is $\beta_{(q,Q,b)}(n+q^{2b+2})=$ $\beta_{(q,Q,b)}(n)$ for all $n$.\bigskip
\end{proof}

\subsection{\label{Fur-coef-comp}Calculation of the Fourier coefficients
$\widehat{\beta_{(q,Q,b)}}(r)$ and $\widehat{\beta_{(q,Q)}}(r)$}

\begin{theorem}
\label{Fur-coef-summary-th}For any prime $q\mid Q$\ \ the Fourier coefficients
$\widehat{\beta_{(q,Q,b)}}(\dfrac{a}{q^{c}})$ are:
\end{theorem}

$c=0,$ \ \ \ \ \ \ \ $\ \ b=0,$ $\ \ \ \ \ \ \ \ \ \ \ \ \ \widehat
{\beta_{(q,Q,0)}}(0)=1-\dfrac{2}{q^{2}(q-1)}$

\mathstrut$c=0,$ \ \ \ \ \ \ \ $\ \ b\neq0,$
$\ \ \ \ \ \ \ \ \ \ \ \ \ \widehat{\beta_{(q,Q,b)}}(0)=\dfrac{2(q^{2}%
+q+1)}{q^{2b+2}}$

\mathstrut$c=2b+2,$ $\ b\neq0,$ \ \ \ \ \ \ \ \ \ \ \ \ $\widehat
{\beta_{(q,Q,b)}}(\dfrac{a}{q^{c}})=\dfrac{2}{q^{2b+2}}\cos\left(  \dfrac{4\pi
a}{q^{c}}\right)  $

$\bigskip c=2b+1,$ $\ b\neq0,$ $\ \ \ \ \ \ \ \ \ \ \ \ \ \widehat
{\beta_{(q,Q,b)}}(\dfrac{a}{q^{c}})=\dfrac{1}{q^{2b+2}}(\left(  1-\dfrac{1}%
{q}\right)  ^{-1}q^{\tfrac{3}{2}}\epsilon_{q}\times$

\bigskip
$\ \ \ \ \ \ \ \ \ \ \ \ \ \ \ \ \ \ \ \ \ \ \ \ \ \ \ \ \ \ \ \ \ \ \times
\left[  e^{-4\pi i\tfrac{a}{q^{c}}}(\dfrac{-a}{q})+e^{4\pi i\tfrac{a}{q^{c}}%
}(\dfrac{a}{q})\right]  -\dfrac{2}{q-1}\cos$ \bigskip$\left(  \dfrac{4\pi
a}{q^{c}}\right)  ),$

$\ \ \ \ \ \ \ \ \ \ \ \ \ \ \ \ \ \ \ \ \ \ \ \ \ \ \ \ \ \ \ \ \ \ \ \ \ \ \ \ \ \ \ \ \ \ \ \ \ \ \ $%
where $\epsilon_{q}=\left\{
\begin{array}
[c]{cc}%
1, & q\equiv1(\operatorname{mod}4)\\
i, & q\equiv3(\operatorname{mod}4)
\end{array}
\right\}  $

\bigskip$c\leq2b,$ $\ \ \ \ \ \ \ \ b\neq0,$
\ \ \ \ \ \ \ \ \ \ \ \ \ $\widehat{\beta_{(q,Q,b)}}(\dfrac{a}{q^{c}}%
)=\dfrac{2}{q^{2b+2}}\cos\left(  \dfrac{4\pi a}{q^{c}}\right)  (q^{2}+q+1)$

\bigskip$c=1,$ \ \ \ \ \ \ \ \ \ $\ b=0,$ \ \ \ \ \ \ \ \ \ \ \ \ \ $\widehat
{\beta_{(q,Q,0)}}(\dfrac{a}{q})=-\dfrac{2}{q^{2}(q-1)}\cos\left(  \dfrac{4\pi
a}{q}\right)  +\ \ \ \ $

$\ \ \ \ \ \ \ \ \ \ \ \ \ \ \ \ \ \ \ \ \ \ \ \ \ \ \ \ \ \ \ \ \ \ \ \ \ \ \ \ \ \ \ \ \ \ \ \ \ \ \ \ \ \ \ \ \ \ \ \ \ \ \ \ \ \ \ \ \ \ \ \ \ \dfrac
{1}{q-1}%
%TCIMACRO{\dsum \limits_{n(\operatorname{mod}q)}}%
%BeginExpansion
{\displaystyle\sum\limits_{n(\operatorname{mod}q)}}
%EndExpansion
\left(  \dfrac{n^{2}-4}{q}\right)  e^{-2\pi in\tfrac{a}{q}}$

\bigskip$c=2,$ \ \ \ \ \ \ \ \ \ $\ b=0,$ \ \ \ \ \ \ \ \ \ \ \ \ $\widehat
{\beta_{(q,Q,0)}}(\dfrac{a}{q^{2}})=$\ $\dfrac{2}{q^{2}}\cos\left(
\dfrac{4\pi a}{q^{2}}\right)  $.%

\[
\]
First we calculate $\widehat{\beta_{(q,Q,b)}}(0).$

\begin{remark}
In all of the sums below we want to be sure that the function $\beta
_{(q,Q,b)}(n)$ is defined at $n$, that is $\ n$ is a trace of some element of
$\Gamma_{0}(Q)$. The necessary and sufficient condition for $n$ to be a trace
of some element of \ $\Gamma_{0}(Q)$ is the condition that $\left(
\frac{n^{2}-4}{q}\right)  \neq-1$ for all $q\mid Q$. We will easily see that
those $n$, for which $\left(  \frac{n^{2}-4}{q}\right)  \neq-1$ do not
contribute to the sum. That is why we can sum over all $n$ of the given range
without any restriction.
\end{remark}

\[
\]
a) $b=0,$

\begin{center}%
\[
\widehat{\beta_{(q,Q,0)}}(0)=\frac{1}{q^{2}}\sum
\limits_{\substack{n(\operatorname{mod}q^{2})}}\left(  1-\frac{1}{q}\left(
\frac{n^{2}-4}{q}\right)  \right)  ^{-1}\cdot\left\{
\begin{array}
[c]{cc}%
2, & q^{2}\mid n^{2}-4\\
\left(  \dfrac{n^{2}-4}{q}\right)  +1, & q^{2}\nmid n^{2}-4
\end{array}
\right\}
\]

\[
=\frac{1}{q^{2}}\sum\limits_{\substack{n(\operatorname{mod}q^{2})}}\left(
1-\frac{1}{q}\left(  \frac{n^{2}-4}{q}\right)  \right)  ^{-1}\cdot\left\{
\begin{array}
[c]{cc}%
2, & n=\pm2,\text{ or }\left(  \frac{n^{2}-4}{q}\right)  =1\\
1, & n\neq\pm2,\text{ }q\mid n^{2}-4
\end{array}
\right\}
\]

\begin{align*}
&  =\frac{1}{q^{2}}(2\cdot2+2\cdot\#\{n(\operatorname{mod}q^{2})\mid\left(
\frac{n^{2}-4}{q}\right)  =1\}\left(  1-\frac{1}{q}\right)  ^{-1}+\\
\#\{n  &  \neq\pm2(\operatorname{mod}q^{2})\mid q\mid n^{2}-4\})=\frac
{1}{q^{2}}(2\cdot2+2\frac{q-3}{2}\cdot q\left(  1-\frac{1}{q}\right)
^{-1}+2(q-1))\\
&  =\frac{1}{q^{2}}(4+q(q-3)\frac{q}{q-1}+2(q-1))=1-\frac{2}{q^{2}(q-1)}.
\end{align*}

\end{center}

\[
\]
b) $b\neq0,$%

\[
\widehat{\beta_{(q,Q,b)}}(0)=\frac{1}{q^{2b+2}}\sum
\limits_{\substack{n(\operatorname{mod}q^{2b+2})}}\left(  1-\frac{1}{q}\left(
\frac{(n^{2}-4)q^{-2b}}{q}\right)  \right)  ^{-1}\times
\]

\begin{center}%
\[
\times\left\{
\begin{array}
[c]{ccc}%
2 & \text{, }n^{2}=4(\operatorname{mod}\text{ }q^{2b}) & \text{, }q^{2}%
\mid(n^{2}-4)q^{-2b}\\
1+\left(  \dfrac{(n^{2}-4)q^{-2b}}{q}\right)  & \text{, }n^{2}%
=4(\operatorname{mod}\text{ }q^{2b}) & \text{, }q^{2}\nmid(n^{2}-4)q^{-2b}\\
0 & \text{, else} &
\end{array}
\right\}  =
\]

\[
=\frac{1}{q^{2b+2}}(2\cdot2+2\left(  1-\frac{1}{q}\right)  ^{-1}\#\left\{
n(\operatorname{mod}q^{2b+2})\mid\left\{
\begin{array}
[c]{c}%
n^{2}\equiv4(\operatorname{mod}q^{2b})\\
\left(  \dfrac{(n^{2}-4)q^{-2b}}{q}\right)  =1
\end{array}
\right\}  \right\}  +
\]

\[
\#\left\{  n\neq\pm2(\operatorname{mod}q^{2b+2})\mid\left\{
\begin{array}
[c]{c}%
n^{2}\equiv4(\operatorname{mod}q^{2b})\\
q^{2b+1}\mid n^{2}-4
\end{array}
\right\}  \right\}  ).
\]

\end{center}

\begin{lemma}
The cardinality of the set \
\[
\left\{  n(\operatorname{mod}q^{2b+2})\mid\left\{
\begin{array}
[c]{c}%
n^{2}\equiv4(\operatorname{mod}q^{2b})\\
\left(  \dfrac{(n^{2}-4)q^{-2b}}{q}\right)  =1
\end{array}
\right\}  \right\}  \text{ \ is }q(q-1),
\]
and the cardinality of the set \
\[
\left\{  n\neq\pm2(\operatorname{mod}q^{2b+2})\mid\left\{
\begin{array}
[c]{c}%
n^{2}\equiv4(\operatorname{mod}q^{2b})\\
q^{2b+1}\mid n^{2}-4
\end{array}
\right\}  \right\}  \text{ is }2q-2.
\]
$\ $
\end{lemma}

\begin{proof}
a) There are $2q^{2}$ numbers $n$ modulo $q^{2b+2}$ such that $n^{2}%
-4=0(\operatorname{mod}q^{2b})$. They are of the form $kq^{2b}$, where
$\ k=\pm1,\pm2,...,\pm q^{2}$.So it's need to check how many of the $k$'s are
squares modulo $q$. There are $(q-1)/2$ squares modulo $q$, hence there are
$2q(q-1)/2=q(q-1)$ numbers in the first set.

\smallskip

b) The number of $n\neq\pm2$ modulo $q^{2b+2}$ such that $n^{2}=4$ modulo
$q^{2b+1}$ is $2q-2$.
\end{proof}

\[
\]
So $\ $%
\[
\ \widehat{\beta_{(q,Q,b)}}(0)=\frac{1}{q^{2b+2}}(4+2\frac{q}{q-1}%
q(q-1)+2q-2)=\frac{2(q^{2}+q+1)}{q^{2b+2}}.
\]
From the relation (\ref{fuco}) it follows that

\begin{center}%
\[
\widehat{\beta_{(q,Q)}}(0)=\sum\limits_{b\geq0}\frac{1}{q^{b}}\widehat
{\beta_{(q,Q,b)}}(0)=1-\frac{2}{q^{2}(q-1)}+\sum\limits_{b\geq1}\frac{1}%
{q^{b}}\frac{2(q^{2}+q+1)}{q^{2b+2}}=
\]

\end{center}

\[
1-\frac{2}{q^{2}(q-1)}+\frac{2(q^{2}+q+1)}{q^{2}}\sum\limits_{b\geq1}\frac
{1}{q^{3b}}=1-\frac{2}{q^{2}(q-1)}+\frac{2(q^{2}+q+1)}{q^{2}}\frac{1}{q^{3}%
-1}=1.
\]
Now we will compute the Fourier coefficients $\widehat{\beta_{(q,Q,b)}}%
(\frac{a}{q^{c}})$.%

\[
\]
a) $b\neq0,$ $c=2b+2.$

We will need a lemmas before we will start .

\begin{lemma}
If $q\nmid a,$ then
\[
\sum\limits_{k=1}^{q^{2}}\left(  \frac{k}{q}\right)  e^{-2\pi ik\frac{a}%
{q^{2}}}=0\text{,}%
\]
for $q$ prime.
\end{lemma}

\begin{proof}
Note that
\[
\left(  \frac{m}{q}\right)  =\left(  \frac{m+lq}{q}\right)  \text{,}%
\]
for $l\in\mathbb{Z}$. Now we can write%
\[
\sum\limits_{k=1}^{q^{2}}\left(  \frac{k}{q}\right)  e^{-2\pi ik\frac{a}%
{q^{2}}}=\sum\limits_{m=1}^{q}\sum\limits_{l=0}^{q-1}\left(  \frac{m+lq}%
{q}\right)  e^{-2\pi i(m+lq)\frac{a}{q^{2}}}=
\]%
\[
\sum\limits_{m=1}^{q}\sum\limits_{l=0}^{q-1}\left(  \frac{m}{q}\right)
e^{-2\pi im\frac{a}{q^{2}}}\cdot e^{-2\pi il\frac{a}{q}}=\sum\limits_{m=1}%
^{q}\left(  \frac{m}{q}\right)  e^{-2\pi im\frac{a}{q^{2}}}\sum\limits_{l=0}%
^{q-1}e^{-2\pi il\frac{a}{q}}.
\]
But
\[
\ \sum\limits_{l=0}^{q-1}e^{-2\pi il\frac{a}{q}}=0
\]
and hence \
\[
\sum\limits_{k=1}^{q^{2}}\left(  \frac{k}{q}\right)  e^{-2\pi ik\frac{a}%
{q^{2}}}=0\text{,}%
\]
as desired.
\end{proof}

\begin{lemma}
For $c=2b+2$ $\ \ \ $

$\ $%
\[
\ \sum\limits_{\substack{n\neq\pm2(\operatorname{mod}q^{2b+2}) \\q^{2b+1}\nmid
n^{2}-4 \\q^{2b}\mid n^{2}-4\text{ }}}\left(  1-\frac{1}{q}\left(
\frac{(n^{2}-4)q^{-2b}}{q}\right)  \right)  ^{-1}\left(  1+\left(
\frac{(n^{2}-4)q^{-2b}}{q}\right)  \right)  e^{-2\pi in\frac{a}{q^{c}}}=0.
\]

\end{lemma}

\begin{proof}%
\[
\sum\limits_{\substack{n\neq\pm2(\operatorname{mod}q^{2b+2})\\q^{2b+1}\nmid
n^{2}-4\\q^{2b}\mid n^{2}-4\text{ }}}\left(  1-\frac{1}{q}\left(  \frac
{(n^{2}-4)q^{-2b}}{q}\right)  \right)  ^{-1}\left(  1+\left(  \frac
{(n^{2}-4)q^{-2b}}{q}\right)  \right)  e^{-2\pi in\frac{a}{q^{c}}}=
\]%
\[
e^{-4\pi i\frac{a}{q^{c}}}\sum\limits_{\substack{k=1\\q\nmid k}}^{q^{2}%
-1}\bigskip\left(  1+\left(  \frac{k}{q}\right)  \right)  \left(  1-\frac
{1}{q}\left(  \frac{k}{q}\right)  \right)  ^{-1}e^{-2\pi ikaq^{2b-c}}+
\]%
\[
e^{4\pi i\frac{a}{q^{c}}}\sum\limits_{\substack{k=1\\q\nmid k}}^{q^{2}%
-1}\left(  1+\left(  \frac{-k}{q}\right)  \right)  \left(  1-\frac{1}%
{q}\left(  \frac{-k}{q}\right)  \right)  ^{-1}e^{-2\pi ikaq^{2b-c}%
}=\{c=2b+2\}=
\]%
\begin{align*}
&  =e^{-4\pi i\frac{a}{q^{c}}}\sum\limits_{k=1,\left(  \frac{k}{q}\right)
=1}^{q^{2}-1}2\left(  1-\frac{1}{q}\right)  ^{-1}e^{-2\pi ik\frac{a}{q^{2}}%
}+\\
&  \text{
\ \ \ \ \ \ \ \ \ \ \ \ \ \ \ \ \ \ \ \ \ \ \ \ \ \ \ \ \ \ \ \ \ \ \ \ \ \ \ \ \ \ \ \ \ \ \ \ \ }%
e^{4\pi i\frac{a}{q^{c}}}\sum\limits_{k=1,\left(  \frac{-k}{q}\right)
=1}^{q^{2}-1}2\left(  1-\frac{1}{q}\right)  ^{-1}e^{-2\pi ik\frac{a}{q^{2}}}%
\end{align*}%
\begin{align*}
&  =e^{-4\pi i\frac{a}{q^{c}}}\cdot2\left(  1-\frac{1}{q}\right)  ^{-1}%
\cdot\frac{1}{2}\sum\limits_{\substack{k=1\\q\nmid k}}^{q^{2}-1}e^{-2\pi
ik\frac{a}{q^{2}}}\left(  1+\left(  \frac{k}{q}\right)  \right)  +\\
&  \text{ \ \ \ \ \ \ \ \ \ \ \ \ \ \ \ \ \ \ \ \ \ \ \ \ \ \ \ \ \ }e^{4\pi
i\frac{a}{q^{c}}}\cdot2\left(  1-\frac{1}{q}\right)  ^{-1}\cdot\frac{1}{2}%
\sum\limits_{\substack{k=1\\q\nmid k}}^{q^{2}-1}e^{-2\pi ik\frac{a}{q^{2}}%
}\left(  1+\left(  \frac{-k}{q}\right)  \right)
\end{align*}%
\[
=e^{-4\pi i\frac{a}{q^{c}}}\left(  1-\frac{1}{q}\right)  ^{-1}\sum
\limits_{\substack{k=1\\q\nmid k}}^{q^{2}-1}\left(  \frac{k}{q}\right)
e^{-2\pi ik\frac{a}{q^{2}}}+e^{-4\pi i\frac{a}{q^{c}}}\left(  1-\frac{1}%
{q}\right)  ^{-1}\sum\limits_{\substack{k=1\\q\nmid k}}^{q^{2}-1}e^{-2\pi
ik\frac{a}{q^{2}}}+
\]%
\[
e^{4\pi i\frac{a}{q^{c}}}\left(  1-\frac{1}{q}\right)  ^{-1}\sum
\limits_{\substack{k=1\\q\nmid k}}^{q^{2}-1}\left(  \frac{-k}{q}\right)
e^{-2\pi ik\frac{a}{q^{2}}}+e^{4\pi i\frac{a}{q^{c}}}\left(  1-\frac{1}%
{q}\right)  ^{-1}\sum\limits_{\substack{k=1\\q\nmid k}}^{q^{2}-1}e^{-2\pi
ik\frac{a}{q^{2}}}.
\]
By the previous lemma and the fact that
\[
\sum\limits_{\substack{k=1\\q\nmid k}}^{q^{2}-1}e^{-2\pi ik\frac{a}{q^{2}}}=0
\]
the expression we want to compute is equal to $0.$
\end{proof}

\[
\]
Now it will be much easier to compute what we want to compute:%

\[
\widehat{\beta_{(q,Q,b)}}(\frac{a}{q^{c}})=\frac{1}{q^{2b+2}}\sum
\limits_{\substack{n(\operatorname{mod}q^{2b+2}) \\\text{ }}}\left(
1-\frac{1}{q}\left(  \frac{(n^{2}-4)q^{-2b}}{q}\right)  \right)  ^{-1}\times
\]

\[
\times\left\{
\begin{array}
[c]{ccc}%
2 & \text{, }n^{2}=4(\operatorname{mod}\text{ }q^{2b}) & \text{, }q^{2}%
\mid(n^{2}-4)q^{-2b}\\
1+\left(  \dfrac{(n^{2}-4)q^{-2b}}{q}\right)  & \text{, }n^{2}%
=4(\operatorname{mod}\text{ }q^{2b}) & \text{, }q^{2}\nmid(n^{2}-4)q^{-2b}\\
0 & \text{, else} &
\end{array}
\right\}  \cdot e^{-2\pi in\tfrac{a}{q^{c}}}%
\]

\[
=\frac{1}{q^{2b+2}}\left(  \sum\limits_{\substack{n=2\\n=-2+q^{2b+2}%
}}2e^{-2\pi in\frac{a}{q^{c}}}\right.  +\text{
\ \ \ \ \ \ \ \ \ \ \ \ \ \ \ \ \ \ \ \ \ \ \ \ \ \ \ \ \ \ \ \ \ \ \ \ \ \ \ \ \ \ \ \ \ \ \ \ \ \ \ \ \ \ \ \ \ \ \ \ \ \ \ \ \ \ \ \ \ \ \ \ \ \ \ \ \ \ \ \ \ \ \ }%
\]%
\[
\left.  \sum\limits_{\substack{n(\operatorname{mod}q^{2b+2})\\n^{2}%
=4(\operatorname{mod}q^{2b})\\q^{2}\nmid(n^{2}-4)q^{-2b}\text{ }}}\left(
1-\frac{1}{q}\left(  \frac{(n^{2}-4)q^{-2b}}{q}\right)  \right)  ^{-1}\left(
1+\left(  \frac{(n^{2}-4)q^{-2b}}{q}\right)  \right)  e^{-2\pi in\frac
{a}{q^{c}}}\right)
\]

\[
=\frac{1}{q^{2b+2}}\left(  \sum\limits_{\substack{n=2\\n=-2+q^{2b+2}%
}}2e^{-2\pi in\frac{a}{q^{c}}}+\sum\limits_{\substack{n\neq\pm
2(\operatorname{mod}q^{2b+2})\\q^{2b+1}\mid n^{2}-4\text{ }}}e^{-2\pi
in\frac{a}{q^{c}}}+\right.  \text{
\ \ \ \ \ \ \ \ \ \ \ \ \ \ \ \ \ \ \ \ \ \ \ \ \ \ \ \ \ \ \ \ \ \ \ \ \ \ \ \ \ \ \ \ \ \ \ }%
\]

\[
\left.  \sum\limits_{\substack{n\neq\pm2(\operatorname{mod}q^{2b+2}%
)\\q^{2b+1}\nmid n^{2}-4\\q^{2b}\mid n^{2}-4\text{ }}}\left(  1-\frac{1}%
{q}\left(  \frac{(n^{2}-4)q^{-2b}}{q}\right)  \right)  ^{-1}\left(  1+\left(
\frac{(n^{2}-4)q^{-2b}}{q}\right)  \right)  e^{-2\pi in\frac{a}{q^{c}}%
}\right)  =
\]

\[
\frac{1}{q^{2b+2}}\left(  \sum\limits_{\substack{n=2\\n=-2+q^{2b+2}}}2e^{-2\pi
in\frac{a}{q^{c}}}+2\cos\left(  \frac{4\pi a}{q^{c}}\right)  \sum
\limits_{k=1}^{q-1}e^{-2\pi ikaq^{2b-c+1}}+0\right)  .
\]
by the previous lemma. Note that%
\[
\sum\limits_{k=1}^{q-1}e^{-2\pi ikaq^{2b-c+1}}=-1
\]
when $c=2b+2$. Therefore%
\[
\widehat{\beta_{(q,Q,b)}}(\frac{a}{q^{c}})=\frac{1}{q^{2b+2}}\left(  2e^{-4\pi
i\frac{a}{q^{c}}}+2e^{4\pi i\frac{a}{q^{c}}}+2\cos\left(  \frac{4\pi a}{q^{c}%
}\right)  \cdot(-1)\right)  =
\]

\[
\frac{1}{q^{2b+2}}\cos\left(  \frac{4\pi a}{q^{c}}\right)  (4-2)=\frac
{2}{q^{2b+2}}\cos\left(  \frac{4\pi a}{q^{c}}\right)  .
\]%
\[
\]
b) $b\neq0,$ $c=2b+1.$%

\[
\]
The same way as before we will need two lemmas:

\begin{lemma}%
\[
\sum\limits_{k=1}^{q^{2}}\left(  \frac{k}{q}\right)  e^{-2\pi ik\frac{a}{q}%
}=q^{\frac{3}{2}}(\frac{-a}{q})\epsilon_{q}\text{, where }\epsilon
_{q}=\left\{
\begin{array}
[c]{cc}%
1, & q\equiv1(\operatorname{mod}4)\\
i, & q\equiv3(\operatorname{mod}4)
\end{array}
\right\}  .
\]

\end{lemma}

\begin{proof}%
\[
\sum\limits_{k=1}^{q^{2}}\left(  \frac{k}{q}\right)  e^{-2\pi ik\frac{a}{q}%
}=\sum\limits_{m=1}^{q}\sum\limits_{l=0}^{q-1}\left(  \frac{m+lq}{q}\right)
e^{-2\pi i(m+lq)\frac{a}{q}}=
\]%
\[
=\sum\limits_{m=1}^{q}\sum\limits_{l=0}^{q-1}\left(  \frac{m}{q}\right)
e^{-2\pi im\frac{a}{q}}\cdot e^{-2\pi ila}=q\sum\limits_{m=1}^{q}\left(
\frac{m}{q}\right)  e^{-2\pi im\frac{a}{q}}=
\]%
\[
=\mathstrut q(\frac{-1}{q})\sum\limits_{m=1}^{q}\left(  \frac{-m}{q}\right)
e^{2\pi i(-m)\frac{a}{q}}=\mathstrut q(\frac{-1}{q})\sum\limits_{m=1}%
^{q}\left(  \frac{m}{q}\right)  e^{2\pi im\frac{a}{q}}=
\]%
\[
=\mathstrut q(\frac{-1}{q})(\frac{a}{q})\sum\limits_{m=1}^{q}\left(  \frac
{m}{q}\right)  e^{2\pi i\frac{m}{q}}.
\]
The equality holds by using of property of the Gaussian sum. The value of the
last sum is $\ \epsilon_{q}\sqrt{q},$ where $\epsilon_{q}$ defined as above.
For this see, for example,\cite[chapter 6]{ire}. Hence
\[
\sum\limits_{k=1}^{q^{2}}\left(  \frac{k}{q}\right)  e^{-2\pi ik\frac{a}{q}%
}=\mathstrut q(\frac{-1}{q})(\frac{a}{q})\epsilon_{q}\sqrt{q}=q^{\frac{3}{2}%
}(\frac{-a}{q})\epsilon_{q}\text{,}%
\]
and we have the claim of the lemma.
\end{proof}

\begin{lemma}
\ \ For $c=2b+1,$ and $\epsilon_{q}$ defined as before%
\[
\sum\limits_{\substack{n\neq\pm2(\operatorname{mod}q^{2b+2}) \\q^{2b+1}\nmid
n^{2}-4 \\q^{2b}\mid n^{2}-4\text{ }}}\left(  1-\frac{1}{q}\left(
\frac{(n^{2}-4)q^{-2b}}{q}\right)  \right)  ^{-1}\left(  1+\left(
\frac{(n^{2}-4)q^{-2b}}{q}\right)  \right)  e^{-2\pi in\frac{a}{q^{c}}}=
\]%
\[
=\left(  1-\frac{1}{q}\right)  ^{-1}q^{\frac{3}{2}}\epsilon_{q}\left[
e^{-4\pi i\frac{a}{q^{c}}}(\frac{-a}{q})+e^{4\pi i\frac{a}{q^{c}}}(\frac{a}%
{q})\right]  -2q\left(  1-\frac{1}{q}\right)  ^{-1}\cos\bigskip\left(
\frac{4\pi a}{q^{c}}\right)  .
\]

\end{lemma}

\begin{proof}%
\[
\bigskip\sum\limits_{\substack{n\neq\pm2(\operatorname{mod}q^{2b+2})
\\q^{2b+1}\nmid n^{2}-4 \\q^{2b}\mid n^{2}-4\text{ }}}\left(  1-\frac{1}%
{q}\left(  \frac{(n^{2}-4)q^{-2b}}{q}\right)  \right)  ^{-1}\left(  1+\left(
\frac{(n^{2}-4)q^{-2b}}{q}\right)  \right)  e^{-2\pi in\frac{a}{q^{c}}}=
\]%
\[
=e^{-4\pi i\frac{a}{q^{c}}}\sum\limits_{\substack{k=1 \\q\nmid k}}^{q^{2}%
-1}\bigskip\left(  1+\left(  \frac{k}{q}\right)  \right)  \left(  1-\frac
{1}{q}\left(  \frac{k}{q}\right)  \right)  ^{-1}e^{-2\pi ikaq^{2b-c}}+
\]%
\[
e^{4\pi i\frac{a}{q^{c}}}\sum\limits_{\substack{k=1 \\q\nmid k}}^{q^{2}%
-1}\left(  1+\left(  \frac{-k}{q}\right)  \right)  \left(  1-\frac{1}%
{q}\left(  \frac{-k}{q}\right)  \right)  ^{-1}e^{-2\pi ikaq^{2b-c}%
}=\{c=2b+1\}=
\]%
\[
=e^{-4\pi i\frac{a}{q^{c}}}\sum\limits_{k=1,\left(  \frac{k}{q}\right)
=1}^{q^{2}-1}2\left(  1-\frac{1}{q}\right)  ^{-1}e^{-2\pi ik\frac{a}{q}}+
\]%
\[
e^{4\pi i\frac{a}{q^{c}}}\sum\limits_{k=1,\left(  \frac{-k}{q}\right)
=1}^{q^{2}-1}2\left(  1-\frac{1}{q}\right)  ^{-1}e^{-2\pi ik\frac{a}{q}}=
\]%
\[
=e^{-4\pi i\frac{a}{q^{c}}}\cdot2\left(  1-\frac{1}{q}\right)  ^{-1}\cdot
\frac{1}{2}\sum\limits_{\substack{k=1 \\q\nmid k}}^{q^{2}-1}e^{-2\pi
ik\frac{a}{q}}\left(  1+\left(  \frac{k}{q}\right)  \right)  +
\]%
\[
e^{4\pi i\frac{a}{q^{c}}}\cdot2\left(  1-\frac{1}{q}\right)  ^{-1}\cdot
\frac{1}{2}\sum\limits_{\substack{k=1 \\q\nmid k}}^{q^{2}-1}e^{-2\pi
ik\frac{a}{q}}\left(  1+\left(  \frac{-k}{q}\right)  \right)  =
\]%
\[
=e^{-4\pi i\frac{a}{q^{c}}}\left(  1-\frac{1}{q}\right)  ^{-1}\sum
\limits_{\substack{k=1 \\q\nmid k}}^{q^{2}-1}\left(  \frac{k}{q}\right)
e^{-2\pi ik\frac{a}{q}}+e^{-4\pi i\frac{a}{q^{c}}}\left(  1-\frac{1}%
{q}\right)  ^{-1}\sum\limits_{\substack{k=1 \\q\nmid k}}^{q^{2}-1}e^{-2\pi
ik\frac{a}{q}}+
\]%
\[
e^{4\pi i\frac{a}{q^{c}}}\left(  1-\frac{1}{q}\right)  ^{-1}\sum
\limits_{\substack{k=1 \\q\nmid k}}^{q^{2}-1}\left(  \frac{-k}{q}\right)
e^{-2\pi ik\frac{a}{q}}+e^{4\pi i\frac{a}{q^{c}}}\left(  1-\frac{1}{q}\right)
^{-1}\sum\limits_{\substack{k=1 \\q\nmid k}}^{q^{2}-1}e^{-2\pi ik\frac{a}{q}%
}.
\]
By the previous lemma and the fact that
\[
\sum\limits_{\substack{k=1 \\q\nmid k}}^{q^{2}-1}e^{-2\pi ik\frac{a}{q}}=-q
\]
we can write that our expression is equal to%
\[
\left(  1-\frac{1}{q}\right)  ^{-1}q^{\frac{3}{2}}\epsilon_{q}\left[  e^{-4\pi
i\frac{a}{q^{c}}}(\frac{-a}{q})+e^{4\pi i\frac{a}{q^{c}}}(\frac{a}{q})\right]
-2q\left(  1-\frac{1}{q}\right)  ^{-1}\cos\bigskip\left(  \frac{4\pi a}{q^{c}%
}\right)  \text{,}%
\]
proving the claim.
\end{proof}

\[
\]
By using the same arguments as in the case a) we will get%

\[
\widehat{\beta_{(q,Q,b)}}(\dfrac{a}{q^{c}})=\dfrac{1}{q^{2b+2}}\left(
%TCIMACRO{\dsum \limits_{\substack{n=2\\n=-2+q^{2b+2}}}}%
%BeginExpansion
{\displaystyle\sum\limits_{\substack{n=2\\n=-2+q^{2b+2}}}}
%EndExpansion
2e^{-2\pi in\frac{a}{q^{c}}}+%
%TCIMACRO{\dsum \limits_{\substack{n\neq\pm2(\operatorname{mod}q^{2b+2}%
%)\\q^{2b+1}\mid n^{2}-4\text{ }}}}%
%BeginExpansion
{\displaystyle\sum\limits_{\substack{n\neq\pm2(\operatorname{mod}%
q^{2b+2})\\q^{2b+1}\mid n^{2}-4\text{ }}}}
%EndExpansion
e^{-2\pi in\frac{a}{q^{c}}}\right.  +\text{
\ \ \ \ \ \ \ \ \ \ \ \ \ \ \ \ \ \ \ \ \ \ \ \ \ \ \ \ \ \ \ \ \ \ \ \ \ \ \ }%
\]

\[
\left.
%TCIMACRO{\dsum \limits_{\substack{n\neq\pm2(\operatorname{mod}q^{2b+2}%
%)\\q^{2b+1}\nmid n^{2}-4\\q^{2b}\mid n^{2}-4\text{ }}}}%
%BeginExpansion
{\displaystyle\sum\limits_{\substack{n\neq\pm2(\operatorname{mod}%
q^{2b+2})\\q^{2b+1}\nmid n^{2}-4\\q^{2b}\mid n^{2}-4\text{ }}}}
%EndExpansion
\left(  1-\dfrac{1}{q}\left(  \dfrac{(n^{2}-4)q^{-2b}}{q}\right)  \right)
^{-1}\left(  1+\left(  \dfrac{(n^{2}-4)q^{-2b}}{q}\right)  \right)  e^{-2\pi
in\frac{a}{q^{c}}}\right)  =
\]

\[
=\dfrac{1}{q^{2b+2}}\left(
%TCIMACRO{\dsum \limits_{\substack{n=2\\n=-2+q^{2b+2}}}}%
%BeginExpansion
{\displaystyle\sum\limits_{\substack{n=2\\n=-2+q^{2b+2}}}}
%EndExpansion
2e^{-2\pi in\frac{a}{q^{c}}}+2\cos\left(  \dfrac{4\pi a}{q^{c}}\right)
%TCIMACRO{\dsum \limits_{k=1}^{q-1}}%
%BeginExpansion
{\displaystyle\sum\limits_{k=1}^{q-1}}
%EndExpansion
e^{-2\pi ika}\right.  +\text{
\ \ \ \ \ \ \ \ \ \ \ \ \ \ \ \ \ \ \ \ \ \ \ \ \ \ \ \ \ \ \ \ \ \ \ \ \ \ \ \ \ \ \ \ \ \ \ \ \ \ \ \ \ \ }%
\]

\[
\left.  \left(  1-\dfrac{1}{q}\right)  ^{-1}q^{\frac{3}{2}}\epsilon_{q}\left[
e^{-4\pi i\frac{a}{q^{c}}}(\dfrac{-a}{q})+e^{4\pi i\frac{a}{q^{c}}}(\dfrac
{a}{q})\right]  -2q\left(  1-\dfrac{1}{q}\right)  ^{-1}\cos\bigskip\left(
\dfrac{4\pi a}{q^{c}}\right)  \right)  .
\]
by using previous lemma. Note that
\[
\sum\limits_{k=1}^{q-1}e^{-2\pi ika}=q-1\text{,}%
\]
so%
\[
\widehat{\beta_{(q,Q,b)}}(\dfrac{a}{q^{c}})=\dfrac{1}{q^{2b+2}}\left(
4\cos\left(  \dfrac{4\pi a}{q^{c}}\right)  +2(q-1)\cos\left(  \dfrac{4\pi
a}{q^{c}}\right)  \right.  +\text{
\ \ \ \ \ \ \ \ \ \ \ \ \ \ \ \ \ \ \ \ \ \ \ \ \ \ \ \ \ \ \ \ \ \ \ \ \ \ \ \ \ }%
\]

\[
\text{ \ \ \ \ \ \ \ \ \ }\left.  \left(  1-\dfrac{1}{q}\right)  ^{-1}%
q^{\frac{3}{2}}\epsilon_{q}\left[  e^{-4\pi i\frac{a}{q^{c}}}(\dfrac{-a}%
{q})+e^{4\pi i\frac{a}{q^{c}}}(\dfrac{a}{q})\right]  -2q\left(  1-\dfrac{1}%
{q}\right)  ^{-1}\cos\bigskip\left(  \dfrac{4\pi a}{q^{c}}\right)  \right)
\]

\[
=\dfrac{1}{q^{2b+2}}\left(  \left(  1-\dfrac{1}{q}\right)  ^{-1}q^{\frac{3}%
{2}}\epsilon_{q}\left[  e^{-4\pi i\frac{a}{q^{c}}}(\dfrac{-a}{q})+e^{4\pi
i\frac{a}{q^{c}}}(\dfrac{a}{q})\right]  -\frac{2}{q-1}\cos\bigskip\left(
\dfrac{4\pi a}{q^{c}}\right)  \right)  .
\]

\[
\]
c) $b\neq0,$ $2b\geq c.$%
\[
\]
By doing the same steps as before we can write

\begin{lemma}
\ For $2b\geq c,$%
\[
\sum\limits_{\substack{n\neq\pm2(\operatorname{mod}q^{2b+2}) \\q^{2b+1}\nmid
n^{2}-4 \\q^{2b}\mid n^{2}-4\text{ }}}\left(  1-\frac{1}{q}\left(
\frac{(n^{2}-4)q^{-2b}}{q}\right)  \right)  ^{-1}\left(  1+\left(
\frac{(n^{2}-4)q^{-2b}}{q}\right)  \right)  e^{-2\pi in\frac{a}{q^{c}}}=
\]%
\[
=2q^{2}\cos\left(  \frac{4\pi a}{q^{c}}\right)  .
\]

\end{lemma}

\begin{proof}%
\[
\sum\limits_{\substack{n\neq\pm2(\operatorname{mod}q^{2b+2}) \\q^{2b+1}\nmid
n^{2}-4 \\q^{2b}\mid n^{2}-4\text{ }}}\left(  1-\frac{1}{q}\left(
\frac{(n^{2}-4)q^{-2b}}{q}\right)  \right)  ^{-1}\left(  1+\left(
\frac{(n^{2}-4)q^{-2b}}{q}\right)  \right)  e^{-2\pi in\frac{a}{q^{c}}}=
\]%
\[
=e^{-4\pi i\frac{a}{q^{c}}}\sum\limits_{\substack{k=1 \\q\nmid k}}^{q^{2}%
-1}\bigskip\left(  1+\left(  \frac{k}{q}\right)  \right)  \left(  1-\frac
{1}{q}\left(  \frac{k}{q}\right)  \right)  ^{-1}e^{-2\pi ikaq^{2b-c}}+
\]%
\[
e^{4\pi i\frac{a}{q^{c}}}\sum\limits_{\substack{k=1 \\q\nmid k}}^{q^{2}%
-1}\left(  1+\left(  \frac{-k}{q}\right)  \right)  \left(  1-\frac{1}%
{q}\left(  \frac{-k}{q}\right)  \right)  ^{-1}e^{-2\pi ikaq^{2b-c}}=\{2b\geq
c\}=
\]%
\[
=e^{-4\pi i\frac{a}{q^{c}}}\sum\limits_{k=1,\left(  \frac{k}{q}\right)
=1}^{q^{2}-1}2\left(  1-\frac{1}{q}\right)  ^{-1}+e^{4\pi i\frac{a}{q^{c}}%
}\sum\limits_{k=1,\left(  \frac{-k}{q}\right)  =1}^{q^{2}-1}2\left(
1-\frac{1}{q}\right)  ^{-1}=
\]
\mathstrut%
\[
=2\frac{q}{q-1}\#\{k(\operatorname{mod}q^{2})\mid(\frac{k}{q})=1\}2\cos\left(
\frac{4\pi a}{q^{c}}\right)  =4\frac{q}{q-1}\frac{q(q-1)}{2}\cos\left(
\frac{4\pi a}{q^{c}}\right)  =
\]%
\[
=2q^{2}\cos\left(  \frac{4\pi a}{q^{c}}\right)  .
\]
$\bigskip$
\end{proof}

$\bigskip$%
\[
\]
And gathering all together we have%
\[
\widehat{\beta_{(q,Q,b)}}(\frac{a}{q^{c}})=\frac{1}{q^{2b+2}}\left(
\sum\limits_{\substack{n=2\\n=-2+q^{2b+2}}}2e^{-2\pi in\frac{a}{q^{c}}}%
+2\cos\left(  \frac{4\pi a}{q^{c}}\right)  \sum\limits_{k=1}^{q-1}e^{-2\pi
ikaq^{2b-c+1}}\right.  +\text{
\ \ \ \ \ \ \ \ \ \ \ \ \ \ \ \ \ \ \ \ \ \ \ \ \ \ }%
\]%
\[
\left.  \sum\limits_{\substack{n\neq\pm2(\operatorname{mod}q^{2b+2}%
)\\q^{2b+1}\nmid n^{2}-4\\q^{2b}\mid n^{2}-4\text{ }}}\left(  1-\frac{1}%
{q}\left(  \frac{(n^{2}-4)q^{-2b}}{q}\right)  \right)  ^{-1}\left(  1+\left(
\frac{(n^{2}-4)q^{-2b}}{q}\right)  \right)  e^{-2\pi in\frac{a}{q^{c}}%
}\right)
\]

\[
=\frac{1}{q^{2b+2}}\left(  4\cos\left(  \frac{4\pi a}{q^{c}}\right)
+2(q-1)\cos\left(  \frac{4\pi a}{q^{c}}\right)  +2q^{2}\cos\left(  \frac{4\pi
a}{q^{c}}\right)  \right)
\]

\[
=\frac{1}{q^{2b+2}}\cos\left(  \frac{4\pi a}{q^{c}}\right)  [4+2(q-1)+2q^{2}%
]=\frac{2}{q^{2b+2}}\cos\left(  \frac{4\pi a}{q^{c}}\right)  (q^{2}+q+1).
\]%
\[
\]
d) $b=0,$ $c=1.$%

\[
\]
As always we need two lemmas:

\begin{lemma}%
\[
\bigskip\sum\limits_{n(\operatorname{mod}q^{2})}\left(  \frac{n^{2}-4}%
{q}\right)  e^{-2\pi in\frac{a}{q}}=q\sum\limits_{n(\operatorname{mod}%
q)}\left(  \frac{n^{2}-4}{q}\right)  e^{-2\pi in\frac{a}{q}}.
\]

\end{lemma}

\begin{proof}%
\[
\bigskip\sum\limits_{n(\operatorname{mod}q^{2})}\left(  \frac{n^{2}-4}%
{q}\right)  e^{-2\pi in\frac{a}{q}}=\sum\limits_{m=1}^{q}\sum\limits_{l=0}%
^{q-1}\left(  \frac{(m+lq)^{2}-4}{q}\right)  e^{-2\pi i(m+lq)\frac{a}{q}}=
\]%
\[
=\sum\limits_{m=1}^{q}\sum\limits_{l=0}^{q-1}\left(  \frac{m^{2}-4}{q}\right)
e^{-2\pi im\frac{a}{q}}=q\sum\limits_{m=1}^{q}\left(  \frac{m^{2}-4}%
{q}\right)  e^{-2\pi im\frac{a}{q}}=
\]%
\[
=q\sum\limits_{n(\operatorname{mod}q)}\left(  \frac{n^{2}-4}{q}\right)
e^{-2\pi in\frac{a}{q}}.
\]

\end{proof}

\begin{lemma}%
\[
\sum\limits_{\substack{n\neq\pm2(\operatorname{mod}q^{2}) \\q\nmid n^{2}%
-4}}\left(  1-\frac{1}{q}\left(  \frac{n^{2}-4}{q}\right)  \right)
^{-1}\left(  \left(  \frac{n^{2}-4}{q}\right)  +1\right)  e^{-2\pi in\frac
{a}{q}}=
\]%
\[
=q\left(  1-\frac{1}{q}\right)  ^{-1}\left[  -2\cos\left(  \frac{4\pi a}%
{q}\right)  +\sum\limits_{n(\operatorname{mod}q)}\left(  \frac{n^{2}-4}%
{q}\right)  e^{-2\pi in\frac{a}{q}}\right]  .
\]
$.$
\end{lemma}

\begin{proof}%
\[
\sum\limits_{\substack{n\neq\pm2(\operatorname{mod}q^{2})\\q\nmid n^{2}%
-4}}\left(  1-\frac{1}{q}\left(  \frac{n^{2}-4}{q}\right)  \right)
^{-1}\left(  \left(  \frac{n^{2}-4}{q}\right)  +1\right)  e^{-2\pi in\frac
{a}{q}}=
\]%
\[
=\sum\limits_{\substack{n(\operatorname{mod}q^{2})\\n\neq\pm
2(\operatorname{mod}q)}}\left(  1-\frac{1}{q}\left(  \frac{n^{2}-4}{q}\right)
\right)  ^{-1}\left(  \left(  \frac{n^{2}-4}{q}\right)  +1\right)  e^{-2\pi
in\frac{a}{q}}=
\]%
\[
=\sum\limits_{\substack{n(\operatorname{mod}q^{2})\\\left(  \frac{n^{2}-4}%
{q}\right)  =1}}2\left(  1-\frac{1}{q}\right)  ^{-1}e^{-2\pi in\frac{a}{q}}=
\]%
\[
=2\cdot\frac{1}{2}\left(  1-\frac{1}{q}\right)  ^{-1}\sum
\limits_{\substack{n(\operatorname{mod}q^{2})\\n\neq\pm2(\operatorname{mod}%
q)}}\left(  \left(  \frac{n^{2}-4}{q}\right)  +1\right)  e^{-2\pi in\frac
{a}{q}}=
\]%
\[
=\left(  1-\frac{1}{q}\right)  ^{-1}\left[  \sum
\limits_{\substack{n(\operatorname{mod}q^{2})\\n\neq\pm2(\operatorname{mod}%
q)}}e^{-2\pi in\frac{a}{q}}+\sum\limits_{\substack{n(\operatorname{mod}%
q^{2})\\n\neq\pm2(\operatorname{mod}q)}}\left(  \frac{n^{2}-4}{q}\right)
e^{-2\pi in\frac{a}{q}}\right]  =
\]%
\[
=\left(  1-\frac{1}{q}\right)  ^{-1}\left(  -2q\cos\left(  \frac{4\pi a}%
{q}\right)  +\sum\limits_{n(\operatorname{mod}q^{2})}\left(  \frac{n^{2}-4}%
{q}\right)  e^{-2\pi in\frac{a}{q}}\right.  -\text{ \ \ \ \ \ \ }%
\]%
\[
\text{\ \ \ \ \ \ \ \ \ \ \ \ \ \ \ \ }\left.  \sum
\limits_{\substack{n(\operatorname{mod}q^{2})\\n=\pm2(\operatorname{mod}%
q)}}\left(  \frac{n^{2}-4}{q}\right)  e^{-2\pi in\frac{a}{q}}\right)
\]%
\[
=\left(  1-\frac{1}{q}\right)  ^{-1}\left(  -2q\cos\left(  \frac{4\pi a}%
{q}\right)  +\sum\limits_{n(\operatorname{mod}q^{2})}\left(  \frac{n^{2}-4}%
{q}\right)  e^{-2\pi in\frac{a}{q}}-0\right)  =
\]%
\[
=\left(  1-\frac{1}{q}\right)  ^{-1}\left[  -2q\cos\left(  \frac{4\pi a}%
{q}\right)  +q\sum\limits_{n(\operatorname{mod}q)}\left(  \frac{n^{2}-4}%
{q}\right)  e^{-2\pi in\frac{a}{q}}\right]  .
\]
by previous lemma, and that's it.
\end{proof}

\[
\]
Now we can gather the results%
\[
\widehat{\beta_{(q,Q,0)}}(\frac{a}{q})=\frac{1}{q^{2}}\sum
\limits_{n(\operatorname{mod}q^{2})}\left(  1-\frac{1}{q}\left(  \frac
{n^{2}-4}{q}\right)  \right)  ^{-1}\cdot\left\{
\begin{array}
[c]{cc}%
2, & q^{2}\mid n^{2}-4\\
\left(  \dfrac{n^{2}-4}{q}\right)  +1, & q^{2}\nmid n^{2}-4
\end{array}
\right\}  \cdot e^{-2\pi in\frac{a}{q}}=
\]

\[
\frac{1}{q^{2}}\left(  \sum\limits_{\substack{n=2\\n=-2+q^{2}}}2e^{-2\pi
in\frac{a}{q}}+\sum\limits_{\substack{n\neq\pm2(\operatorname{mod}%
q^{2})\\q\mid n^{2}-4}}e^{-2\pi in\frac{a}{q}}+\right.  \text{
\ \ \ \ \ \ \ \ \ \ \ \ \ \ \ \ \ \ \ \ \ \ \ \ \ \ \ \ \ \ \ \ \ \ \ \ \ \ \ \ \ \ \ \ \ \ \ \ \ \ \ \ \ \ \ \ \ \ \ \ \ \ \ \ \ \ \ \ \ \ \ \ \ \ \ \ \ \ \ \ \ \ \ \ \ \ \ \ \ \ \ }%
\]

\[
\text{ \ \ \ \ \ \ \ \ \ \ \ }\left.  \sum\limits_{\substack{n\neq
\pm2(\operatorname{mod}q^{2})\\q\nmid n^{2}-4}}\left(  1-\frac{1}{q}\left(
\frac{n^{2}-4}{q}\right)  \right)  ^{-1}\left(  \left(  \frac{n^{2}-4}%
{q}\right)  +1\right)  e^{-2\pi in\frac{a}{q}}\right)  =
\]

\mathstrut%

\[
\frac{1}{q^{2}}\left(  4\cos\left(  \frac{4\pi a}{q}\right)  +2(q-1)\cos
\left(  \frac{4\pi a}{q}\right)  \right.  +\text{
\ \ \ \ \ \ \ \ \ \ \ \ \ \ \ \ \ \ \ \ \ \ \ \ \ \ \ \ \ \ \ \ \ \ \ \ \ \ \ \ \ \ \ \ \ \ \ \ \ \ \ \ \ \ \ \ \ \ \ \ \ \ \ \ \ \ \ \ \ \ \ \ \ \ \ \ \ \ \ \ \ \ \ \ \ \ \ \ \ \ \ \ \ \ \ }%
\]%
\[
\text{ \ \ \ \ \ \ \ \ \ \ \ \ \ \ \ }\left.  q\left(  1-\frac{1}{q}\right)
^{-1}\left[  -2\cos\left(  \frac{4\pi a}{q}\right)  +\sum
\limits_{n(\operatorname{mod}q)}\left(  \frac{n^{2}-4}{q}\right)  e^{-2\pi
in\frac{a}{q}}\right]  \right)  =
\]

\mathstrut%

\[
-\frac{2}{q^{2}(q-1)}\cos\left(  \frac{4\pi a}{q}\right)  +\frac{1}{q-1}%
\sum\limits_{n(\operatorname{mod}q)}\left(  \frac{n^{2}-4}{q}\right)  e^{-2\pi
in\frac{a}{q}}.
\]

\[
\]
e) $b=0,$ $c=2.$

\begin{lemma}%
\[
\ \sum\limits_{\substack{n\neq\pm2(\operatorname{mod}q^{2}) \\q\mid n^{2}%
-4}}e^{-2\pi in\frac{a}{q^{2}}}=-2\cos\left(  \frac{4\pi a}{q^{2}}\right)  .
\]

\end{lemma}

\begin{proof}
$\bigskip$

The condition \{ $n\neq\pm2(\operatorname{mod}q^{2}),$ $q\mid n^{2}-4\}$ we
can rewrite in the \ form \ 

\{ $n=\pm2+kq,$ $k=1,...,q-1$ \}. Thus%
\[
\sum\limits_{\substack{n\neq\pm2(\operatorname{mod}q^{2}) \\q\mid n^{2}%
-4}}e^{-2\pi in\frac{a}{q^{2}}}=\sum\limits_{k=1}^{q-1}e^{-2\pi i(\pm
2+kq)\frac{a}{q^{2}}}=
\]
$\ \ \ \ \ \ \ \ \ \ \ \ \ \ \ \ \ \ \ \ \ \ \ \ \ $%
\[
=2\cos\left(  \frac{4\pi a}{q^{2}}\right)  \sum\limits_{k=1}^{q-1}e^{-2\pi
ik\frac{a}{q}}=-2\cos\left(  \frac{4\pi a}{q^{2}}\right)  \text{,}%
\]
and we are done.
\end{proof}

\begin{lemma}%
\[
\sum\limits_{\substack{n\neq\pm2(\operatorname{mod}q^{2}) \\q\nmid n^{2}%
-4}}\left(  1-\frac{1}{q}\left(  \frac{n^{2}-4}{q}\right)  \right)
^{-1}\left(  \left(  \frac{n^{2}-4}{q}\right)  +1\right)  e^{-2\pi in\frac
{a}{q^{2}}}=0.
\]

\end{lemma}

\begin{proof}%
\[
\sum\limits_{\substack{n\neq\pm2(\operatorname{mod}q^{2}) \\q\nmid n^{2}%
-4}}\left(  1-\frac{1}{q}\left(  \frac{n^{2}-4}{q}\right)  \right)
^{-1}\left(  \left(  \frac{n^{2}-4}{q}\right)  +1\right)  e^{-2\pi in\frac
{a}{q^{2}}}=
\]%
\[
=\sum\limits_{\substack{n(\operatorname{mod}q^{2}) \\n\neq\pm
2(\operatorname{mod}q)}}\left(  1-\frac{1}{q}\left(  \frac{n^{2}-4}{q}\right)
\right)  ^{-1}\left(  \left(  \frac{n^{2}-4}{q}\right)  +1\right)  e^{-2\pi
in\frac{a}{q^{2}}}=
\]%
\[
=\sum\limits_{\substack{n(\operatorname{mod}q^{2}) \\\left(  \frac{n^{2}-4}%
{q}\right)  =1}}2\left(  1-\frac{1}{q}\right)  ^{-1}e^{-2\pi in\frac{a}{q^{2}%
}}=
\]%
\[
=2\cdot\frac{1}{2}\left(  1-\frac{1}{q}\right)  ^{-1}\sum
\limits_{\substack{n(\operatorname{mod}q^{2}) \\n\neq\pm2(\operatorname{mod}%
q)}}\left(  \left(  \frac{n^{2}-4}{q}\right)  +1\right)  e^{-2\pi in\frac
{a}{q^{2}}}=
\]%
\[
=\left(  1-\frac{1}{q}\right)  ^{-1}\left[  \sum
\limits_{\substack{n(\operatorname{mod}q^{2}) \\n\neq\pm2(\operatorname{mod}%
q)}}e^{-2\pi in\frac{a}{q^{2}}}+\sum\limits_{\substack{n(\operatorname{mod}%
q^{2}) \\n\neq\pm2(\operatorname{mod}q)}}\left(  \frac{n^{2}-4}{q}\right)
e^{-2\pi in\frac{a}{q^{2}}}\right]  .
\]
The first sum in the brackets is $0$, and so we need to compute an expression%
\[
\left(  1-\frac{1}{q}\right)  ^{-1}\sum\limits_{\substack{n(\operatorname{mod}%
q^{2}) \\n\neq\pm2(\operatorname{mod}q)}}\left(  \frac{n^{2}-4}{q}\right)
e^{-2\pi in\frac{a}{q^{2}}}=
\]
$\ \ \ \ \ \ \ \ \ \ \ \ \ \ \ \ \ \ \ \ \ \ \ \ \ \ \ \ \ \ \ \ \ \ \ \ $%
\[
=\left(  1-\frac{1}{q}\right)  ^{-1}\sum\limits_{n(\operatorname{mod}q^{2}%
)}\left(  \frac{n^{2}-4}{q}\right)  e^{-2\pi in\frac{a}{q^{2}}}=
\]%
\[
=\left(  1-\frac{1}{q}\right)  ^{-1}\sum\limits_{m=1}^{q}\sum\limits_{l=0}%
^{q-1}\left(  \frac{\left(  m+lq\right)  ^{2}-4}{q}\right)  e^{-2\pi
i(m+lq)\frac{a}{q^{2}}}=
\]
$\ \ \ \ \ \ \ \ \ \ \ \ \ \ \ \ \ \ \ \ \ \ \ $%
\[
=\left(  1-\frac{1}{q}\right)  ^{-1}\sum\limits_{m=1}^{q}\sum\limits_{l=0}%
^{q-1}\left(  \frac{m^{2}-4}{q}\right)  e^{-2\pi im\frac{a}{q^{2}}}\cdot
e^{-2\pi il\frac{a}{q}}=
\]%
\[
=\left(  1-\frac{1}{q}\right)  ^{-1}\sum\limits_{m=1}^{q}\left(  \frac
{m^{2}-4}{q}\right)  e^{-2\pi im\frac{a}{q^{2}}}\sum\limits_{l=0}%
^{q-1}e^{-2\pi il\frac{a}{q}}=0\text{,}%
\]
since the last sum is 0.
\end{proof}

$\bigskip$%

\[
\widehat{\beta_{(q,Q,0)}}(\frac{a}{q^{2}})=\frac{1}{q^{2}}\sum
\limits_{n(\operatorname{mod}q^{2})}\left(  1-\frac{1}{q}\left(  \frac
{n^{2}-4}{q}\right)  \right)  ^{-1}\cdot\left\{
\begin{array}
[c]{cc}%
2, & q^{2}\mid n^{2}-4\\
\left(  \frac{n^{2}-4}{q}\right)  +1, & q^{2}\nmid n^{2}-4
\end{array}
\right\}  \cdot e^{-2\pi in\frac{a}{q^{2}}}=
\]

\[
=\frac{1}{q^{2}}\left(  \sum\limits_{\substack{n=2\\n=-2+q^{2}}}2e^{-2\pi
in\frac{a}{q^{2}}}+\sum\limits_{\substack{n\neq\pm2(\operatorname{mod}%
q^{2})\\q\mid n^{2}-4}}e^{-2\pi in\frac{a}{q^{2}}}\right.  +\text{
\ \ \ \ \ \ \ \ \ \ \ \ \ \ \ \ \ \ \ \ \ \ \ \ \ \ \ \ \ \ \ \ \ \ \ \ \ \ \ \ \ \ \ \ \ \ \ \ \ \ \ \ \ \ \ \ \ \ \ \ \ \ \ }%
\]

\begin{center}%
\[
\text{ \ \ \ \ \ \ \ \ \ \ }\left.  \sum\limits_{\substack{n\neq
\pm2(\operatorname{mod}q^{2})\\q\nmid n^{2}-4}}\left(  1-\frac{1}{q}\left(
\frac{n^{2}-4}{q}\right)  \right)  ^{-1}\left(  \left(  \frac{n^{2}-4}%
{q}\right)  +1\right)  e^{-2\pi in\frac{a}{q^{2}}}\right)  =
\]

\end{center}

\[
=\frac{1}{q^{2}}\left(  4\cos\left(  \frac{4\pi a}{q^{2}}\right)
-2\cos\left(  \frac{4\pi a}{q^{2}}\right)  +0\right)  =\frac{2}{q^{2}}%
\cos\left(  \frac{4\pi a}{q^{2}}\right)
\]

\[
\]
Summarize the above results:

\bigskip

$c=0,$ \ \ \ \ \ \ \ $b=0,$ $\ \ \ \ \ \ \ \ \ \ \ \ \ \widehat{\beta
_{(q,Q,0)}}(0)=1-\dfrac{2}{q^{2}(q-1)}$

$\bigskip$

$c=0,$ \ \ \ \ \ \ \ $b\neq0,$ $\ \ \ \ \ \ \ \ \ \ \ \ \ \widehat
{\beta_{(q,Q,b)}}(0)=\dfrac{2(q^{2}+q+1)}{q^{2b+2}}$

$\bigskip$

$c=2b+2,$ $b\neq0,$ \ \ \ \ \ \ \ \ \ \ \ \ $\widehat{\beta_{(q,Q,b)}}%
(\dfrac{a}{q^{c}})=\dfrac{2}{q^{2b+2}}\cos\left(  \dfrac{4\pi a}{q^{c}%
}\right)  $

$\bigskip$

$c=2b+1,$ $b\neq0,$ $\ \ \ \ \ \ \ \ \ \ \ \ \ \widehat{\beta_{(q,Q,b)}%
}(\dfrac{a}{q^{c}})=\dfrac{1}{q^{2b+2}}\left(  \left(  1-\dfrac{1}{q}\right)
^{-1}q^{\frac{3}{2}}\epsilon_{q}\right.  \times$

\bigskip

$\ \ \ \ \ \ \ \ \ \ \ \ \ \ \ \ \ \ \ \ \ \ \ \ \ \ \ \ \ \ \ \ \ \ \ \ \ \ \ \ \ \ \ \ \ \ \ \times
\left.  \left[  e^{-4\pi i\frac{a}{q^{c}}}(\dfrac{-a}{q})+e^{4\pi i\frac
{a}{q^{c}}}(\dfrac{a}{q})\right]  -\dfrac{2}{q-1}\cos\bigskip\left(
\dfrac{4\pi a}{q^{c}}\right)  \right)  $

\bigskip

$c\leq2b,$ $\ \ \ \ \ \ \ b\neq0,$ \ \ \ \ \ \ \ \ \ \ \ \ \ $\widehat
{\beta_{(q,Q,b)}}(\dfrac{a}{q^{c}})=\dfrac{2}{q^{2b+2}}\cos\left(  \dfrac{4\pi
a}{q^{c}}\right)  (q^{2}+q+1)$

\bigskip

\bigskip$c=1,$ \ \ \ \ \ \ \ \ \ $b=0,$ \ \ \ \ \ \ \ \ \ \ \ \ \ $\widehat
{\beta_{(q,Q,0)}}(\dfrac{a}{q})=-\dfrac{2}{q^{2}(q-1)}\cos\left(  \dfrac{4\pi
a}{q}\right)  +$

$\ \ \ \ $

$\ \ \ \ \ \ \ \ \ \ \ \ \ \ \ \ \ \ \ \ \ \ \ \ \ \ \ \ \ \ \ \ \ \ \ \ \ \ \ \ \ \ \ \ \ \ \ \ \ \ \ \ \ \ \ \ \ \ \ \ \ \ \ \ \ \ \ \ \ \ \ \ \ \dfrac
{1}{q-1}%
%TCIMACRO{\dsum \limits_{n(\operatorname{mod}q)}}%
%BeginExpansion
{\displaystyle\sum\limits_{n(\operatorname{mod}q)}}
%EndExpansion
\left(  \dfrac{n^{2}-4}{q}\right)  e^{-2\pi in\frac{a}{q}}$

\bigskip

$c=2,$ \ \ \ \ \ \ \ \ \ $\ b=0,$ \ \ \ \ \ \ \ \ \ \ \ \ $\widehat
{\beta_{(q,Q,0)}}(\dfrac{a}{q^{2}})=$\ $\dfrac{2}{q^{2}}\cos\left(
\dfrac{4\pi a}{q^{2}}\right)  $.%

\[
\]
The theorem is completely proved.%
\[
\]
Now we can compute the $\widehat{\beta_{(q,Q)}}(r)$ .

\begin{theorem}
\bigskip\ The Fourier coefficients $\widehat{\beta_{(q,Q)}}(\dfrac{a}{q^{c}})
$ are:%
\begin{align*}
\widehat{\beta_{(q,Q)}}\left(  \frac{a}{q}\right)   &  =\frac{1}{q-1}%
\sum\limits_{n(\operatorname{mod}q)}\left(  \frac{n^{2}-4}{q}\right)  e^{-2\pi
in\frac{a}{q}}\\
\widehat{\beta_{(q,Q)}}\left(  \frac{a}{q^{2}}\right)   &  =\frac{2}%
{q(q-1)}\cos\left(  \frac{4\pi a}{q^{2}}\right) \\
\widehat{\beta_{(q,Q)}}\left(  \frac{a}{q^{c}}\right)   &  =\frac{2}{q-1}%
\cos\left(  \frac{4\pi a}{q^{c}}\right)  \frac{1}{q^{\frac{3c-4}{2}}},\text{
for }c>2,c\text{ even}\\
\widehat{\beta_{(q,Q)}}\left(  \frac{a}{q^{c}}\right)   &  =\frac{1}{q-1}%
\frac{1}{q^{\frac{3c-4}{2}}}\epsilon_{q}\left(  \frac{a}{q}\right)  \left[
e^{-4\pi i\frac{a}{q^{c}}}(\frac{-1}{q})+e^{4\pi i\frac{a}{q^{c}}}\right]
,\text{ for }c>2,c\text{ odd}%
\end{align*}

\end{theorem}

\begin{proof}
By using (\ref{fuco}) we get

\bigskip a) for $c=1,$

\mathstrut

$\widehat{\beta_{(q,Q)}}\left(  \dfrac{a}{q}\right)  =%
%TCIMACRO{\dsum \limits_{b\geq0}}%
%BeginExpansion
{\displaystyle\sum\limits_{b\geq0}}
%EndExpansion
\dfrac{1}{q^{b}}\widehat{\beta_{(q,Q,b)}}\left(  \dfrac{a}{q}\right)
=-\dfrac{2}{q^{2}(q-1)}\cos\left(  \dfrac{4\pi a}{q}\right)  +\ \dfrac{1}{q-1}%
%TCIMACRO{\dsum \limits_{n(\operatorname{mod}q)}}%
%BeginExpansion
{\displaystyle\sum\limits_{n(\operatorname{mod}q)}}
%EndExpansion
\left(  \dfrac{n^{2}-4}{q}\right)  e^{-2\pi in\frac{a}{q}}+$

\mathstrut

$%
%TCIMACRO{\dsum \limits_{b\geq1}}%
%BeginExpansion
{\displaystyle\sum\limits_{b\geq1}}
%EndExpansion
\dfrac{1}{q^{b}}\widehat{\beta_{(q,Q,b)}}\left(  \dfrac{a}{q}\right)
=-\dfrac{2}{q^{2}(q-1)}\cos\left(  \dfrac{4\pi a}{q}\right)  +\ \dfrac{1}{q-1}%
%TCIMACRO{\dsum \limits_{n(\operatorname{mod}q)}}%
%BeginExpansion
{\displaystyle\sum\limits_{n(\operatorname{mod}q)}}
%EndExpansion
\left(  \dfrac{n^{2}-4}{q}\right)  e^{-2\pi in\frac{a}{q}}+$

\mathstrut

$%
%TCIMACRO{\dsum \limits_{b\geq1}}%
%BeginExpansion
{\displaystyle\sum\limits_{b\geq1}}
%EndExpansion
\dfrac{1}{q^{b}}\dfrac{2}{q^{2b+2}}\cos\left(  \dfrac{4\pi a}{q}\right)
(q^{2}+q+1)=-\dfrac{2}{q^{2}(q-1)}\cos\left(  \dfrac{4\pi a}{q}\right)
+\ \dfrac{1}{q-1}%
%TCIMACRO{\dsum \limits_{n(\operatorname{mod}q)}}%
%BeginExpansion
{\displaystyle\sum\limits_{n(\operatorname{mod}q)}}
%EndExpansion
\left(  \dfrac{n^{2}-4}{q}\right)  e^{-2\pi in\frac{a}{q}}+$

\mathstrut

$2\cos\left(  \dfrac{4\pi a}{q}\right)  (q^{2}+q+1)\dfrac{1}{q^{2}(q^{3}%
-1)}=\dfrac{1}{q-1}%
%TCIMACRO{\dsum \limits_{n(\operatorname{mod}q)}}%
%BeginExpansion
{\displaystyle\sum\limits_{n(\operatorname{mod}q)}}
%EndExpansion
\left(  \dfrac{n^{2}-4}{q}\right)  e^{-2\pi in\frac{a}{q}}.$

\bigskip b) for $c=2,$

\mathstrut

$\widehat{\beta_{(q,Q)}}\left(  \dfrac{a}{q^{2}}\right)  =%
%TCIMACRO{\dsum \limits_{b\geq0}}%
%BeginExpansion
{\displaystyle\sum\limits_{b\geq0}}
%EndExpansion
\dfrac{1}{q^{b}}\widehat{\beta_{(q,Q,b)}}\left(  \dfrac{a}{q^{2}}\right)
=\dfrac{2}{q^{2}}\cos\left(  \dfrac{4\pi a}{q^{2}}\right)  +%
%TCIMACRO{\dsum \limits_{b\geq1}}%
%BeginExpansion
{\displaystyle\sum\limits_{b\geq1}}
%EndExpansion
\dfrac{1}{q^{b}}\widehat{\beta_{(q,Q,b)}}\left(  \dfrac{a}{q^{2}}\right)  =$

\mathstrut

$=\dfrac{2}{q^{2}}\cos\left(  \dfrac{4\pi a}{q^{2}}\right)  +%
%TCIMACRO{\dsum \limits_{b\geq1}}%
%BeginExpansion
{\displaystyle\sum\limits_{b\geq1}}
%EndExpansion
\dfrac{1}{q^{b}}\dfrac{2}{q^{2b+2}}\cos\left(  \dfrac{4\pi a}{q^{2}}\right)
(q^{2}+q+1)=$

\mathstrut

$=\dfrac{2}{q^{2}}\cos\left(  \dfrac{4\pi a}{q^{2}}\right)  +2\cos\left(
\dfrac{4\pi a}{q^{2}}\right)  (q^{2}+q+1)\dfrac{1}{q^{2}(q^{3}-1)}=\dfrac
{2}{q(q-1)}\cos\left(  \dfrac{4\pi a}{q^{2}}\right)  .$

\bigskip c) for $c>2,$

\mathstrut

$\widehat{\beta_{(q,Q)}}\left(  \dfrac{a}{q^{c}}\right)  =%
%TCIMACRO{\dsum \limits_{0\leq b<\frac{c-2}{2}}}%
%BeginExpansion
{\displaystyle\sum\limits_{0\leq b<\frac{c-2}{2}}}
%EndExpansion
\dfrac{1}{q^{b}}\cdot0+%
%TCIMACRO{\dsum \limits_{b=\frac{c-2}{2}}}%
%BeginExpansion
{\displaystyle\sum\limits_{b=\frac{c-2}{2}}}
%EndExpansion
\dfrac{1}{q^{b}}\dfrac{2}{q^{2b+2}}\cos\left(  \dfrac{4\pi a}{q^{c}}\right)
+$

\mathstrut

$%
%TCIMACRO{\dsum \limits_{b=\frac{c-1}{2}}}%
%BeginExpansion
{\displaystyle\sum\limits_{b=\frac{c-1}{2}}}
%EndExpansion
\dfrac{1}{q^{b}}\left[  \dfrac{1}{q^{2b+2}}\left(  \left(  1-\dfrac{1}%
{q}\right)  ^{-1}q^{\frac{3}{2}}\epsilon_{q}\left[  e^{-4\pi i\frac{a}{q^{c}}%
}(\dfrac{-a}{q})+e^{4\pi i\frac{a}{q^{c}}}(\dfrac{a}{q})\right]  -\dfrac
{2}{q-1}\cos\bigskip\left(  \dfrac{4\pi a}{q^{c}}\right)  \right)  \right]  +$

\mathstrut

$%
%TCIMACRO{\dsum \limits_{b\geq\frac{c}{2}}}%
%BeginExpansion
{\displaystyle\sum\limits_{b\geq\frac{c}{2}}}
%EndExpansion
\dfrac{1}{q^{b}}\dfrac{2}{q^{2b+2}}\cos\left(  \dfrac{4\pi a}{q^{c}}\right)
(q^{2}+q+1).$

\mathstrut

For $c$ even we have

\mathstrut\mathstrut

\bigskip$\widehat{\beta_{(q,Q)}}\left(  \dfrac{a}{q^{c}}\right)  =\dfrac
{1}{q^{\frac{c-2}{2}}}\dfrac{2}{q^{c}}\cos\left(  \dfrac{4\pi a}{q^{c}%
}\right)  +\dfrac{2}{q^{2}}\cos\left(  \dfrac{4\pi a}{q^{c}}\right)
(q^{2}+q+1)%
%TCIMACRO{\dsum \limits_{b\geq\frac{c}{2}}}%
%BeginExpansion
{\displaystyle\sum\limits_{b\geq\frac{c}{2}}}
%EndExpansion
\dfrac{1}{q^{3b}}=$

\mathstrut

\bigskip$=2\cos\left(  \dfrac{4\pi a}{q^{c}}\right)  \left[  \dfrac
{1}{q^{\frac{3c-2}{2}}}+\dfrac{1}{q^{2}}(q^{2}+q+1)\dfrac{1}{q^{\frac{3c-6}%
{2}}}\dfrac{1}{q^{3}-1}\right]  =\dfrac{2}{q-1}\cos\left(  \dfrac{4\pi
a}{q^{c}}\right)  \dfrac{1}{q^{\frac{3c-4}{2}}}.$

\bigskip And for $c$ odd

\mathstrut

$\dfrac{1}{q^{\frac{c-1}{2}}}\left[  \dfrac{1}{q^{c+1}}\left(  \left(
1-\dfrac{1}{q}\right)  ^{-1}q^{\frac{3}{2}}\epsilon_{q}\left[  e^{-4\pi
i\frac{a}{q^{c}}}(\dfrac{-a}{q})+e^{4\pi i\frac{a}{q^{c}}}(\dfrac{a}%
{q})\right]  -\dfrac{2}{q-1}\cos\bigskip\left(  \dfrac{4\pi a}{q^{c}}\right)
\right)  \right]  +$

\mathstrut

$\dfrac{2}{q^{2}}\cos\left(  \dfrac{4\pi a}{q^{c}}\right)  (q^{2}+q+1)%
%TCIMACRO{\dsum \limits_{b\geq\frac{c+1}{2}}}%
%BeginExpansion
{\displaystyle\sum\limits_{b\geq\frac{c+1}{2}}}
%EndExpansion
\dfrac{1}{q^{3b}}=$

\mathstrut

$=\dfrac{1}{q^{\frac{3c+1}{2}}}\left(  1-\dfrac{1}{q}\right)  ^{-1}q^{\frac
{3}{2}}\epsilon_{q}\left[  e^{-4\pi i\frac{a}{q^{c}}}(\dfrac{-a}{q})+e^{4\pi
i\frac{a}{q^{c}}}(\dfrac{a}{q})\right]  -\dfrac{1}{q^{\frac{3c+1}{2}}}%
\dfrac{2}{q-1}\cos\bigskip\left(  \dfrac{4\pi a}{q^{c}}\right)  +$

\mathstrut

$\dfrac{2}{q^{2}}\cos\left(  \dfrac{4\pi a}{q^{c}}\right)  (q^{2}%
+q+1)\dfrac{1}{q^{\frac{3c-3}{2}}}\dfrac{1}{q^{3}-1}=$

\mathstrut

$=\dfrac{1}{q^{\frac{3c-4}{2}}}\dfrac{1}{q-1}\epsilon_{q}\left[  e^{-4\pi
i\frac{a}{q^{c}}}(\dfrac{-a}{q})+e^{4\pi i\frac{a}{q^{c}}}(\dfrac{a}%
{q})\right]  -\dfrac{1}{q^{\frac{3c+1}{2}}}\dfrac{2}{q-1}\cos\bigskip\left(
\dfrac{4\pi a}{q^{c}}\right)  +\dfrac{1}{q^{\frac{3c+1}{2}}}\dfrac{2}{q-1}%
\cos\bigskip\left(  \dfrac{4\pi a}{q^{c}}\right)  =$

$=\dfrac{1}{q-1}\dfrac{1}{q^{\frac{3c-4}{2}}}\epsilon_{q}\left(  \dfrac{a}%
{q}\right)  \left[  e^{-4\pi i\frac{a}{q^{c}}}(\dfrac{-1}{q})+e^{4\pi
i\frac{a}{q^{c}}}\right]  .$
\end{proof}

\subsection{\label{mean-qs-calc-sect}Calculating the mean square of weighted
multiplicities function}

In this subsection we will calculate the mean-square of the weighted
multiplicities $\beta_{Q}(n).$%

\[
\]
To calculate the limit
\[
\lim\limits_{N\rightarrow\infty}\frac{1}{N}\sum\limits_{2\leq n\leq N}%
\beta_{Q}^{2}(n)
\]
we will use (\ref{bet-sq-mean}). Let us define the function $\ \ \ $%
\begin{equation}
A_{Q}(p^{c}):=\sum\limits_{\substack{1\leq a\leq p^{c}\\p\nmid a}}\left\vert
\widehat{\beta_{(p,Q)}}(\frac{a}{p^{c}})\right\vert ^{2}. \label{A-Q}%
\end{equation}
\ So we have that%

\[
\lim\limits_{N\rightarrow\infty}\frac{1}{N}\sum\limits_{2\leq n\leq N}%
\beta_{Q}^{2}(n)=\prod\limits_{p\text{ - prime}}\left(  1+\sum\limits_{c\geq
1}A_{Q}(p^{c})\right)  .
\]
The values of the $A_{Q}(p^{c})$ for $p\nmid Q$ were calculate by M.Peter
\cite{man}. They are:%
\[%
\begin{array}
[c]{ccc}%
p\neq2,p\nmid Q; & A_{Q}(p)=\dfrac{p^{2}-2p-1}{(p^{2}-1)^{2}}, & A_{Q}%
(p^{c})=\dfrac{2(p-1)}{(p^{2}-1)^{2}p^{2c-3}}\\
&  & \\
p=2; & A_{Q}(2)=\dfrac{1}{9}, & A_{Q}(4)=\dfrac{1}{18},\\
&  & \\
& A_{Q}(8)=0, & \ A_{Q}(16)=\dfrac{1}{9\cdot16},\\
&  & \\
& A_{Q}(32)=0, & A_{Q}(2^{c})=\dfrac{1}{9\cdot2^{2c-5}},c\geq6.
\end{array}
\]
We just need to complete his work by adding the case $q\mid Q.$

a) $c=1;$%

\[
A_{Q}(q)=\sum\limits_{\substack{1\leq a\leq q \\q\nmid a}}\left\vert
\widehat{\beta_{(q,Q)}}(\frac{a}{q})\right\vert ^{2}=\sum
\limits_{\substack{1\leq a\leq q \\q\nmid a}}\left\vert \frac{1}{q-1}%
\sum\limits_{n(\operatorname{mod}q)}\left(  \frac{n^{2}-4}{q}\right)  e^{-2\pi
in\frac{a}{q}}\right\vert ^{2}=
\]

\[
=\frac{1}{(q-1)^{2}}\sum\limits_{\substack{1\leq a\leq q \\q\nmid a}%
}\sum\limits_{n_{1},n_{2}(\operatorname{mod}q)}\left(  \frac{n_{1}^{2}-4}%
{q}\right)  \left(  \frac{n_{2}^{2}-4}{q}\right)  e^{2\pi i(n_{1}-n_{2}%
)\frac{a}{q}}=
\]

\[
=\frac{1}{(q-1)^{2}}\sum\limits_{n_{1},n_{2}(\operatorname{mod}q)}\left(
\frac{n_{1}^{2}-4}{q}\right)  \left(  \frac{n_{2}^{2}-4}{q}\right)
\sum\limits_{a=1}^{q-1}e^{2\pi i(n_{1}-n_{2})\frac{a}{q}}.
\]
The sum
\[
\sum\limits_{a=1}^{q-1}e^{2\pi i(n_{1}-n_{2})\frac{a}{q}}=\left\{
\begin{array}
[c]{cc}%
q-1, & n_{1}=n_{2}=n(\operatorname{mod}q)\\
-1, & \text{else}%
\end{array}
\right.  \text{.}%
\]
Note that%

\[
\sum\limits_{n(\operatorname{mod}q)}\left(  \frac{n^{2}-4}{q}\right)
=\frac{q-3}{2}-\frac{q-1}{2}=-1,
\]
hence%

\[
A_{Q}(q)=\frac{1}{(q-1)^{2}}\left(  (q-1)\sum\limits_{n(\operatorname{mod}%
q)}\left(  \frac{n^{2}-4}{q}\right)  ^{2}-\sum\limits_{n_{1}%
(\operatorname{mod}q)}\sum\limits_{n_{2}\neq n_{1}(\operatorname{mod}%
q)}\left(  \frac{n_{1}^{2}-4}{q}\right)  \left(  \frac{n_{2}^{2}-4}{q}\right)
\right)  =
\]

\[
=\frac{1}{(q-1)^{2}}\left(  (q-1)(q-2)-\sum\limits_{n_{1}(\operatorname{mod}%
q)}\left(  \sum\limits_{n_{2}(\operatorname{mod}q)}\left(  \frac{n_{1}^{2}%
-4}{q}\right)  \left(  \frac{n_{2}^{2}-4}{q}\right)  -\left(  \frac{n_{1}%
^{2}-4}{q}\right)  ^{2}\right)  \right)  =
\]

\[
=\frac{1}{(q-1)^{2}}\left(  (q-1)(q-2)-\left(  1-\sum\limits_{n_{1}%
(\operatorname{mod}q)}\left(  \frac{n_{1}^{2}-4}{q}\right)  ^{2}\right)
\right)  =
\]

\[
=\frac{1}{(q-1)^{2}}\left(  q^{2}-3q+2-(1-(q-2))\right)  =\frac{q^{2}%
-2q-1}{(q-1)^{2}}.
\]

\[
\]
b) $c=2;$%

\[
A_{Q}(q^{2})=\sum\limits_{\substack{1\leq a\leq q^{2} \\q\nmid a}}\left\vert
\widehat{\beta_{(q,Q)}}(\frac{a}{q^{2}})\right\vert ^{2}=\sum
\limits_{\substack{1\leq a\leq q^{2} \\q\nmid a}}\left\vert \frac{2}%
{q(q-1)}\cos\left(  \frac{4\pi a}{q^{2}}\right)  \right\vert ^{2}=
\]

\[
=\frac{4}{q^{2}(q-1)^{2}}\sum\limits_{\substack{1\leq a\leq q^{2} \\q\nmid
a}}\left\vert \cos\left(  \frac{4\pi a}{q^{2}}\right)  \right\vert ^{2}%
=\frac{1}{q^{2}(q-1)^{2}}\sum\limits_{\substack{1\leq a\leq q^{2} \\q\nmid
a}}\left(  e^{8\pi i\frac{a}{q^{2}}}+e^{-8\pi i\frac{a}{q^{2}}}+2\right)  =
\]

\[
=\frac{1}{q^{2}(q-1)^{2}}\left(  2(q^{2}-q)+\sum\limits_{\substack{1\leq a\leq
q^{2} \\q\nmid a}}\left(  e^{8\pi i\frac{a}{q^{2}}}+e^{-8\pi i\frac{a}{q^{2}}%
}\right)  \right)  =\frac{2(q^{2}-q)}{q^{2}(q-1)^{2}}=\frac{2}{q(q-1)}.
\]

\[
\]
c1) $c>2,$ $c$ even;%

\[
A_{Q}(q^{c})=\sum\limits_{\substack{1\leq a\leq q^{c} \\q\nmid a}}\left\vert
\widehat{\beta_{(q,Q)}}(\frac{a}{q^{c}})\right\vert ^{2}=\sum
\limits_{\substack{1\leq a\leq q^{c} \\q\nmid a}}\left\vert \frac{2}{q-1}%
\cos\left(  \frac{4\pi a}{q^{c}}\right)  \frac{1}{q^{\frac{3c-4}{2}}%
}\right\vert ^{2}=
\]

\[
=\frac{1}{q^{3c-4}}\frac{4}{(q-1)^{2}}\sum\limits_{\substack{1\leq a\leq q^{c}
\\q\nmid a}}\left\vert \cos\left(  \frac{4\pi a}{q^{c}}\right)  \right\vert
^{2}=\frac{1}{q^{3c-4}}\frac{1}{(q-1)^{2}}\sum\limits_{\substack{1\leq a\leq
q^{c} \\q\nmid a}}\left(  e^{8\pi i\frac{a}{q^{c}}}+e^{-8\pi i\frac{a}{q^{c}}%
}+2\right)  =
\]

\[
=\frac{1}{q^{3c-4}}\frac{1}{(q-1)^{2}}\left(  2(q^{c}-q^{c-1})+\sum
\limits_{\substack{1\leq a\leq q^{c} \\q\nmid a}}\left(  e^{8\pi i\frac
{a}{q^{c}}}+e^{-8\pi i\frac{a}{q^{c}}}\right)  \right)  =\frac{2q^{c-1}%
(q-1)}{q^{3c-4}(q-1)^{2}}=\frac{2}{q^{2c-3}(q-1)}.
\]

\[
\]
c2) $c>2,$ $c$ odd;%

\[
A_{Q}(q^{c})=\sum\limits_{\substack{1\leq a\leq q^{c} \\q\nmid a}}\left\vert
\widehat{\beta_{(q,Q)}}(\frac{a}{q^{c}})\right\vert ^{2}=\sum
\limits_{\substack{1\leq a\leq q^{c} \\q\nmid a}}\left\vert \frac{1}{q-1}%
\frac{1}{q^{\frac{3c-4}{2}}}\epsilon_{q}\left(  \frac{a}{q}\right)  \left[
e^{-4\pi i\frac{a}{q^{c}}}(\frac{-1}{q})+e^{4\pi i\frac{a}{q^{c}}}\right]
\right\vert ^{2}=
\]

\[
=\frac{1}{(q-1)^{2}}\frac{1}{q^{3c-4}}\sum\limits_{\substack{1\leq a\leq q^{c}
\\q\nmid a}}\epsilon_{q}\left[  e^{-4\pi i\frac{a}{q^{c}}}(\frac{-1}%
{q})+e^{4\pi i\frac{a}{q^{c}}}\right]  \cdot\overline{\epsilon_{q}\left[
e^{-4\pi i\frac{a}{q^{c}}}(\frac{-1}{q})+e^{4\pi i\frac{a}{q^{c}}}\right]  }=
\]

\[
=\frac{1}{(q-1)^{2}}\frac{1}{q^{3c-4}}\sum\limits_{\substack{1\leq a\leq
q^{c}\\q\nmid a}}\left(  2+\left(  \frac{-1}{q}\right)  \left[  e^{8\pi
i\frac{a}{q^{c}}}+e^{-8\pi i\frac{a}{q^{c}}}\right]  \right)  =\frac
{2q^{c-1}(q-1)}{q^{3c-4}(q-1)^{2}}=\frac{2}{q^{2c-3}(q-1)}.
\]
And we can see that for $c>2,$ $A_{Q}(q^{c})$ does not depend on parity of
$c.$ Now we have by \cite{man}%

\[%
\begin{array}
[c]{ccc}%
p\neq2,p\nmid Q; & A_{Q}(p)=\dfrac{p^{2}-2p-1}{(p^{2}-1)^{2}}, & A_{Q}%
(p^{c})=\dfrac{2(p-1)}{(p^{2}-1)^{2}p^{2c-3}}\\
&  & \\
p=2; & A_{Q}(2)=\dfrac{1}{9}, & A_{Q}(4)=\dfrac{1}{18},\\
&  & \\
& A_{Q}(8)=0, & \ A_{Q}(16)=\dfrac{1}{9\cdot16},\\
&  & \\
& A_{Q}(32)=0, & A_{Q}(2^{c})=\dfrac{1}{9\cdot2^{2c-5}},c\geq6.
\end{array}
\]
and what we have find, for $q\mid Q,$%
\[%
\begin{array}
[c]{ccc}%
A_{Q}(q)=\dfrac{q^{2}-2q-1}{(q-1)^{2}}, & \ A_{Q}(q^{2})=\dfrac{2}{q(q-1)}, &
A_{Q}(q^{c})=\dfrac{2}{q^{2c-3}(q-1)}.
\end{array}
\]
Now we can calculate%

\[
\lim\limits_{N\rightarrow\infty}\frac{1}{N}\sum\limits_{2\leq n\leq N}%
\beta_{Q}^{2}(n)=\prod\limits_{p\text{ - prime}}\left(  1+\sum\limits_{c\geq
1}A_{Q}(p^{c})\right)  =\left(  1+\frac{1}{9}+\frac{1}{18}+\frac{1}{9\cdot
16}+\sum\limits_{c\geq6}\frac{1}{9\cdot2^{2c-5}}\right)  \times
\]

\[
\times\prod\limits_{q\mid Q}\left(  1+\frac{q^{2}-2q-1}{(q-1)^{2}}+\frac
{2}{q(q-1)}+\sum\limits_{c>2}\frac{2}{q^{2c-3}(q-1)}\right)  \cdot
\prod\limits_{\substack{p\neq2 \\p\nmid Q}}\left(  1+\frac{p^{2}-2p-1}%
{(p^{2}-1)^{2}}+\sum\limits_{c\geq2}\frac{2(p-1)}{(p^{2}-1)^{2}p^{2c-3}%
}\right)  =
\]

\[
=\frac{1015}{864}\prod\limits_{q\mid Q}\frac{2q(q^{2}-q-1)}{(q+1)(q-1)^{2}%
}\prod\limits_{\substack{p\neq2 \\p\nmid Q}}\frac{p^{2}(p^{3}+p^{2}%
-p-3)}{(p^{2}-1)^{2}(p+1)}.
\]
That is \ \ \ \ \ \ \ %

\[
\ \lim\limits_{N\rightarrow\infty}\frac{1}{N}\sum\limits_{2\leq n\leq N}%
\beta_{Q}^{2}(n)=\frac{1015}{864}\prod\limits_{q\mid Q}\frac{2q(q^{2}%
-q-1)}{(q+1)(q-1)^{2}}\prod\limits_{\substack{p\neq2 \\p\nmid Q}}\frac
{p^{2}(p^{3}+p^{2}-p-3)}{(p^{2}-1)^{2}(p+1)}=
\]

\[
=\ \frac{1015}{864}\prod\limits_{q\mid Q}\frac{2q(q^{2}-q-1)}{(q+1)(q-1)^{2}%
}\frac{(q^{2}-1)^{2}(q+1)}{q^{2}(q^{3}+q^{2}-q-3)}\ \prod\limits_{p\neq2}%
\frac{p^{2}(p^{3}+p^{2}-p-3)}{(p^{2}-1)^{2}(p+1)}=
\]

\[
=\frac{1015}{864}\prod\limits_{q\mid Q}\allowbreak\frac{\allowbreak2\left(
q+1\right)  ^{2}(q^{2}-q-1)}{q\left(  q^{3}+q^{2}-q-3\right)  }\ \prod
\limits_{p\neq2}\frac{p^{2}(p^{3}+p^{2}-p-3)}{(p^{2}-1)^{2}(p+1)}%
\]%
\[
=\left(  \prod\limits_{q\mid Q}\frac{2(q^{2}-q-1)(q+1)^{2}}{q(q^{3}%
+q^{2}-q-3)}\right)  \cdot1.328...=C_{1}\prod\limits_{q\mid Q}\frac
{2(q^{2}-q-1)(q+1)^{2}}{q(q^{3}+q^{2}-q-3)},
\]
\ proving the result pointed out in introduction.

\section{\label{Quater-sect}Proof of the Theorem \ref{mean-sq-theorem(quater)}%
. Main points}

Let $B$ be an indefinite division quaternion algebra over $\mathbb{Q}$ with
discriminant $d_{B}$, and $R$ be the \textit{maximal} order in $B$. Let
\[
\Gamma_{R}=\{\alpha\in R\mid N_{B}(\alpha)=1\}
\]
be the unit group of $R$. Define the weighted multiplicities function in the
similar way as for $\Gamma_{0}(Q)$:
\[
\beta_{R}(n)=\dfrac{1}{4}%
%TCIMACRO{\dsum \limits_{\substack{\{T\}~\subset\Gamma_{R}\\T=T_{0}%
%^{k}~is~hyperbolic\\\left\vert trT\right\vert =n}}}%
%BeginExpansion
{\displaystyle\sum\limits_{\substack{\{T\}~\subset\Gamma_{R}\\T=T_{0}%
^{k}~is~hyperbolic\\\left\vert trT\right\vert =n}}}
%EndExpansion
\dfrac{\ln\mathcal{N}(T_{0})}{\mathcal{N}(T)^{1/2}-\mathcal{N}(T)^{-1/2}%
}\text{.}%
\]
The first step is to express weighted multiplicities in the terms of
Dirichlet's $L$-functions. The analog of the theorem \ref{theorem} is

\begin{theorem}
\bigskip Let $B$ be indefinite division quaternion algebra over $\mathbb{Q}$,
and $d_{B}$ be its reduced discriminant. Then
\[
\beta_{R}(n)=%
%TCIMACRO{\dsum \limits_{\substack{D,v\geq1 \\Dv^{2}=n^{2}-4}}}%
%BeginExpansion
{\displaystyle\sum\limits_{\substack{D,v\geq1 \\Dv^{2}=n^{2}-4}}}
%EndExpansion
\dfrac{1}{v}L(1,\chi_{D})%
%TCIMACRO{\dprod \limits_{p\mid d_{B}}}%
%BeginExpansion
{\displaystyle\prod\limits_{p\mid d_{B}}}
%EndExpansion
\left\{
\begin{array}
[c]{cc}%
1-\left(  \dfrac{D}{p}\right)  , & p^{2}\nmid D\\
0, & p^{2}\mid D
\end{array}
\right\}  \text{,}%
\]
where $D$ is a discriminant, i.e. $D\equiv0,1(\operatorname{mod}4)$.
\end{theorem}

\begin{proof}
One uses exactly the same techniques as for the $\Gamma_{0}(Q)$ case.
\end{proof}

\[
\]
In the second step we change $L$-function in the last formula by the Euler's
product formula, and define%
\[
\beta_{P,R}(n):=%
%TCIMACRO{\dsum \limits_{\substack{D,v\geq1 \\Dv^{2}=n^{2}-4 \\p\mid
%v\Rightarrow p\leq P}}}%
%BeginExpansion
{\displaystyle\sum\limits_{\substack{D,v\geq1 \\Dv^{2}=n^{2}-4 \\p\mid
v\Rightarrow p\leq P}}}
%EndExpansion
\left(  \dfrac{1}{v}%
%TCIMACRO{\dprod \limits_{p\leq P}}%
%BeginExpansion
{\displaystyle\prod\limits_{p\leq P}}
%EndExpansion
\left(  1-\dfrac{\chi_{D}(p)}{p}\right)  ^{-1}\right)
%TCIMACRO{\dprod \limits_{p\mid d_{B}}}%
%BeginExpansion
{\displaystyle\prod\limits_{p\mid d_{B}}}
%EndExpansion
\left\{
\begin{array}
[c]{cc}%
1-\left(  \dfrac{D}{p}\right)  , & p^{2}\nmid D\\
0, & p^{2}\mid D
\end{array}
\right\}  \text{.}%
\]
We have the local factor decomposition for this function \cite{man},
\cite{luk}:

\begin{lemma}
Let $P\geq d_{B}$, then
\[
\beta_{P,R}(n)=%
%TCIMACRO{\dprod \limits_{p\leq P}}%
%BeginExpansion
{\displaystyle\prod\limits_{p\leq P}}
%EndExpansion
\beta_{(p,R)}(n),
\]
where%
\[
\beta_{(p,R)}(n):=%
%TCIMACRO{\dsum \limits_{b\geq0}}%
%BeginExpansion
{\displaystyle\sum\limits_{b\geq0}}
%EndExpansion
\dfrac{1}{p^{b}}\left(  1-\dfrac{1}{p}\chi_{(n^{2}-4)p^{-2b}}(p)\right)
^{-1}\cdot\mathbb{I}_{p^{b}}(n).
\]

\end{lemma}

\begin{remark}
\bigskip Here the definition of $\mathbb{I}_{p^{b}}(n)$ is
\begin{align*}
\mathbb{I}_{2^{b}}(n)  &  =\left\{
\begin{array}
[c]{cc}%
1, & n^{2}=4(\operatorname{mod}2^{2b})\text{, }(n^{2}-4)2^{-2b}\text{ is a
discriminant}\\
0, & \text{else}%
\end{array}
\right.  \text{;}\\
& \\
\underset{p\nmid d_{B}}{\mathbb{I}_{p^{b}}(n)}  &  =\left\{
\begin{array}
[c]{cc}%
1, & n^{2}=4(\operatorname{mod}p^{2b})\text{, \ }p\neq2\\
0, & \text{else}%
\end{array}
\right.  \text{;}\\
& \\
\underset{p\mid d_{B}}{\mathbb{I}_{p^{b}}(n)}  &  =\left\{
\begin{array}
[c]{cc}%
1-\left(  \dfrac{(n^{2}-4)p^{-2b}}{p}\right)  , & n^{2}=4(\operatorname{mod}%
p^{2b})\text{, \ }p^{2}\nmid(n^{2}-4)p^{-2b}\\
0, & \text{else}%
\end{array}
\right.  \text{.}%
\end{align*}

\end{remark}%

\[
\]
For the next result we use the uniform limit periodicity of $\beta_{(p,R)}$
and once again (see lemma \ref{fucomult}) get the formula
\[
\widehat{\beta_{R}}(\frac{a}{b})=\prod\limits_{p\mid b}\widehat{\beta_{(p,R)}%
}(\frac{a_{p}}{p^{ord_{p}b}})
\]
for the Fourier coefficients of $\beta_{R}$. \ It implies%
\begin{align}
\lim\limits_{N\rightarrow\infty}\frac{1}{N}\sum\limits_{2\leq n\leq N}%
\beta_{R}^{2}(n)  &  =\sum\limits_{b\geq1}\sum\limits_{\substack{1\leq a\leq b
\\\gcd(a,b)=1}}\left\vert \widehat{\beta_{R}}\left(  \frac{a}{b}\right)
\right\vert ^{2}=\label{betameansquare(quater)}\\
& \nonumber\\
&  =\prod\limits_{p\text{ - prime}}\left(  1+\sum\limits_{c\geq1}%
\sum\limits_{\substack{1\leq a\leq p^{c} \\a\neq0(\operatorname{mod}%
p)}}\left\vert \widehat{\beta_{(p,R)}}\left(  \frac{a}{p^{c}}\right)
\right\vert ^{2}\right)  \text{,}\nonumber
\end{align}
by using Parseval's equality.%
\[
\]
Defining%
\[
\beta_{(p,R,b)}(n):=\left(  1-\frac{1}{p}\chi_{(n^{2}-4)p^{-2b}}(p)\right)
^{-1}\cdot\mathbb{I}_{p^{b}}(n),
\]
one obtain \cite{man}
\[
\widehat{\beta_{(p,R)}}(r)=\sum\limits_{b\geq0}\frac{1}{p^{b}}\widehat
{\beta_{(p,R,b)}}(r)\text{.}%
\]
For the Fourier coefficients in the right-hand side we have

\begin{theorem}
For any prime $p\mid d_{B}$\ \ the Fourier coefficients $\widehat
{\beta_{(p,R,b)}}(\dfrac{a}{p^{c}})$ are:

\mathstrut

$c=0,$ \ \ \ \ \ \ \ $\ \ b=0,$ $\ \ \ \ \ \ \ \ \ \ \ \ \ \widehat
{\beta_{(p,R,0)}}(0)=\dfrac{(p-1)(p^{2}+2p+2)}{p^{2}(p+1)}$\mathstrut

\mathstrut$c=0,$ \ \ \ \ \ \ \ $\ \ b\neq0,$
$\ \ \ \ \ \ \ \ \ \ \ \ \ \widehat{\beta_{(p,R,b)}}(0)=\dfrac{2(p^{3}%
-1)}{p^{2b+2}(p+1)}$

\mathstrut$c=2b+2,$ $\ \ \ b\neq0,$ \ \ \ \ \ \ \ \ \ \ \ \ $\ \widehat
{\beta_{(p,R,b)}}(\dfrac{a}{p^{c}})=-\dfrac{2}{p^{2b+2}}\cos\left(
\dfrac{4\pi a}{p^{c}}\right)  $

$\bigskip c=2b+1,$ $\ \ \ b\neq0,$ $\ \ \ \ \ \ \ \ \ \ \ \ \ \widehat
{\beta_{(p,R,b)}}(\dfrac{a}{p^{c}})=-\dfrac{2}{p^{2b+2}(p+1)}\cos
\bigskip\left(  \dfrac{4\pi a}{p^{c}}\right)  -$

$\ \ \ \ \ \ \ \ \ \ \ \ \ \ \ \ \ \ \ \ \ \ \ \ \ \ \ \ \ \ \ \ \ \ \ \ \ \ \ \ \ \ \ \ \ \ \ \ \ \ \ \ \ \ \ \ \ \ \ \ \ \ \ \ \ \ \ \ \ \ \ -\ \dfrac
{p^{1/2}\epsilon_{p}}{p^{2b}(p+1)}\left[  e^{-4\pi i\tfrac{a}{p^{c}}}%
(\dfrac{-a}{p})+e^{4\pi i\tfrac{a}{p^{c}}}(\dfrac{a}{p})\right]  ,$

\mathstrut

$\ \ \ \ \ \ \ \ \ \ \ \ \ \ \ \ \ \ \ \ \ \ \ \ \ \ \ \ \ \ \ \ \ \ \ \ \ \ \ \ \ \ \ \ \ \ \ \ \ \ \ $%
where $\epsilon_{p}=\left\{
\begin{array}
[c]{cc}%
1, & p\equiv1(\operatorname{mod}4)\\
i, & p\equiv3(\operatorname{mod}4)
\end{array}
\right\}  $

\bigskip$c\leq2b,$ $\ \ \ \ \ \ \ \ b\neq0,$
\ \ \ \ \ \ \ \ \ \ \ \ \ \ $\widehat{\beta_{(p,R,b)}}(\dfrac{a}{p^{c}%
})=\dfrac{2(p^{3}-1)}{p^{2b+2}(p+1)}\cos\left(  \dfrac{4\pi a}{p^{c}}\right)
$

$c=1,$ \ \ \ \ \ \ \ \ $\ b=0,$ \ \ \ \ \ \ \ \ \ \ \ \ \ $\ \widehat
{\beta_{(p,R,0)}}(\dfrac{a}{p})=-\dfrac{2}{p^{2}(p+1)}\cos\left(  \dfrac{4\pi
a}{p}\right)  -\ \ $

$\ \ \mathstrut$

$\ \ \ \ \ \ \ \ \ \ \ \ \ \ \ \ \ \ \ \ \ \ \ \ \ \ \ \ \ \ \ \ \ \ \ \ \ \ \ \ \ \ \ \ \ \ \ \ \ \ \ \ \ \ \ \ \ \ \ \ \ \ \ \ \ \ \ \ \ \ \ \ \ -\dfrac
{1}{p+1}%
%TCIMACRO{\dsum \limits_{n(\operatorname{mod}p)}}%
%BeginExpansion
{\displaystyle\sum\limits_{n(\operatorname{mod}p)}}
%EndExpansion
\left(  \dfrac{n^{2}-4}{p}\right)  e^{-2\pi in\tfrac{a}{p}}$

\bigskip$c=2,$ \ \ \ \ \ \ \ \ \ $\ b=0,$ \ \ \ \ \ \ \ \ \ \ \ \ \ $\widehat
{\beta_{(p,R,0)}}(\dfrac{a}{p^{2}})=$\ $-\dfrac{2}{p^{2}}\cos\left(
\dfrac{4\pi a}{p^{2}}\right)  $.
\end{theorem}

\begin{remark}
\bigskip The case $p\nmid d_{B}$ is done by M.Peter \cite{man}. We will use
these results later.
\end{remark}

\begin{corollary}
For any prime $p\mid d_{B}$\ the Fourier coefficients\ $\ \widehat
{\beta_{(p,R)}}(\dfrac{a}{p^{c}})$ are:

$\widehat{\beta_{(p,R)}}\left(  0\right)  =1$

\mathstrut

$\widehat{\beta_{(p,R)}}\left(  \dfrac{a}{p}\right)  =\dfrac{-1}{p+1}%
\sum\limits_{n(\operatorname{mod}p)}\left(  \dfrac{n^{2}-4}{p}\right)
e^{-2\pi in\frac{a}{p}}$

\mathstrut

$\widehat{\beta_{(p,R)}}\left(  \dfrac{a}{p^{2}}\right)  =\dfrac{-2}%
{p(p+1)}\cos\left(  \dfrac{4\pi a}{p^{2}}\right)  $

\mathstrut

$\widehat{\beta_{(p,R)}}\left(  \dfrac{a}{p^{c}}\right)  =\dfrac
{-2}{(p+1)p^{\frac{3c-4}{2}}}\cos\left(  \dfrac{4\pi a}{p^{c}}\right)  ,$ for
$2<c$ even

\mathstrut

$\widehat{\beta_{(p,R)}}\left(  \dfrac{a}{p^{c}}\right)  =\dfrac{-\epsilon
_{p}}{(p+1)p^{\frac{3c-4}{2}}}\left(  \dfrac{a}{p}\right)  \left[  e^{-4\pi
i\frac{a}{p^{c}}}(\dfrac{-1}{p})+e^{4\pi i\frac{a}{p^{c}}}\right]  ,$ for
$2<c$ odd
\end{corollary}

\[
\]
Now we able to calculate%
\[
\lim\limits_{N\rightarrow\infty}\frac{1}{N}\sum\limits_{2\leq n\leq N}%
\beta_{R}^{2}(n)\text{.}%
\]
Let us define the function $\ \ \ $%
\[
A_{R}(p^{c}):=\sum\limits_{\substack{1\leq a\leq p^{c} \\p\nmid a}}\left\vert
\widehat{\beta_{(p,R)}}(\frac{a}{p^{c}})\right\vert ^{2}.
\]
\ Using (\ref{betameansquare(quater)}) we find that%

\[
\lim\limits_{N\rightarrow\infty}\frac{1}{N}\sum\limits_{2\leq n\leq N}%
\beta_{R}^{2}(n)=\prod\limits_{p\text{ - prime}}\left(  1+\sum\limits_{c\geq
1}A_{R}(p^{c})\right)  .
\]
The values of the $A_{R}(p^{c})$ for $p\nmid d_{B}$ were calculate by M.Peter
\cite{man}. They are:%
\[%
\begin{array}
[c]{ccc}%
p\neq2,p\nmid d_{B}; & A_{R}(p)=\dfrac{p^{2}-2p-1}{(p^{2}-1)^{2}}, &
A_{R}(p^{c})=\dfrac{2(p-1)}{(p^{2}-1)^{2}p^{2c-3}}\\
&  & \\
p=2; & A_{R}(2)=\dfrac{1}{9}, & A_{R}(4)=\dfrac{1}{18},\\
&  & \\
& A_{R}(8)=0, & \ A_{R}(16)=\dfrac{1}{9\cdot16},\\
&  & \\
& A_{R}(32)=0, & A_{R}(2^{c})=\dfrac{1}{9\cdot2^{2c-5}},c\geq6.
\end{array}
\]
We complete his results by adding the case $p\mid d_{B}$. By simple
calculations we have

\begin{theorem}
Let $p\mid d_{B}$, then
\[
A_{R}(p)=\dfrac{p^{2}-2p-1}{(p+1)^{2}};~~\ \ A_{R}(p^{2})=\dfrac
{2(p-1)}{p(p+1)^{2}};~~~~~A_{R}(p^{c})=\dfrac{2(p-1)}{(p+1)^{2}p^{2c-3}}.
\]

\end{theorem}%

\[
\]
Gathering the results we conclude%
\begin{align*}
\lim\limits_{N\rightarrow\infty}\frac{1}{N}\sum\limits_{2\leq n\leq N}%
\beta_{R}^{2}(n)  &  =\prod\limits_{p\text{ - prime}}\left(  1+\sum
\limits_{c\geq1}A_{R}(p^{c})\right)  =\\
&  \left(  1+\frac{1}{9}+\frac{1}{18}+\frac{1}{9\cdot16}+\sum\limits_{c\geq
6}\frac{1}{9\cdot2^{2c-5}}\right)  \times\\
&  \times%
%TCIMACRO{\dprod \limits_{p\mid d_{B}}}%
%BeginExpansion
{\displaystyle\prod\limits_{p\mid d_{B}}}
%EndExpansion
\left(  1+\dfrac{p^{2}-2p-1}{(p+1)^{2}}+\dfrac{2(p-1)}{p(p+1)^{2}}+%
%TCIMACRO{\dsum \limits_{c>2}}%
%BeginExpansion
{\displaystyle\sum\limits_{c>2}}
%EndExpansion
\dfrac{2(p-1)}{(p+1)^{2}p^{2c-3}}\right)  \times\\
&  \times\prod\limits_{\substack{p\neq2\\p\nmid d_{B}}}\left(  1+\frac
{p^{2}-2p-1}{(p^{2}-1)^{2}}+\sum\limits_{c\geq2}\frac{2(p-1)}{(p^{2}%
-1)^{2}p^{2c-3}}\right) \\
&  =\dfrac{1015}{864}%
%TCIMACRO{\dprod \limits_{p\mid d_{B}}}%
%BeginExpansion
{\displaystyle\prod\limits_{p\mid d_{B}}}
%EndExpansion
\dfrac{2p(p^{2}+p+1)}{(p+1)^{3}}\prod\limits_{\substack{p\neq2\\p\nmid d_{B}%
}}\frac{p^{2}(p^{3}+p^{2}-p-3)}{(p^{2}-1)^{2}(p+1)}.
\end{align*}
And finally we have
\begin{align*}
\lim\limits_{N\rightarrow\infty}\frac{1}{N}\sum\limits_{2\leq n\leq N}%
\beta_{R}^{2}(n)  &  =\dfrac{1015}{864}%
%TCIMACRO{\dprod \limits_{p\mid d_{B}}}%
%BeginExpansion
{\displaystyle\prod\limits_{p\mid d_{B}}}
%EndExpansion
\dfrac{2p(p^{2}+p+1)}{(p+1)^{3}}\prod\limits_{\substack{p\neq2\\p\nmid d_{B}%
}}\frac{p^{2}(p^{3}+p^{2}-p-3)}{(p^{2}-1)^{2}(p+1)}\\
& \\
&  =\dfrac{1015}{864}%
%TCIMACRO{\dprod \limits_{p\mid d_{B}}}%
%BeginExpansion
{\displaystyle\prod\limits_{p\mid d_{B}}}
%EndExpansion
\dfrac{2p(p^{2}+p+1)}{(p+1)^{3}}\cdot\dfrac{(p^{2}-1)^{2}(p+1)}{p^{2}%
(p^{3}+p^{2}-p-3)}\times\\
&  \text{ \ \ \ \ \ \ \ }\times\prod\limits_{p\neq2}\frac{p^{2}(p^{3}%
+p^{2}-p-3)}{(p^{2}-1)^{2}(p+1)}\\
& \\
&  =C_{1}\cdot%
%TCIMACRO{\dprod \limits_{p\mid d_{B}}}%
%BeginExpansion
{\displaystyle\prod\limits_{p\mid d_{B}}}
%EndExpansion
\dfrac{2(p^{3}-1)(p-1)}{p(p^{3}+p^{2}-p-3)}\text{, \ \ \ where }C_{1}=1.328...
\end{align*}
This calculation concludes the proof.

\end{document}